\documentclass[11pt,reqno,a4paper]{amsart}

\usepackage[left=25mm,right=25mm,top=25mm,bottom=25mm]{geometry}
\usepackage{amsmath,amssymb,amsthm,mathtools}
\usepackage{booktabs,array}
\usepackage{enumitem}
\usepackage{microtype}
\usepackage[unicode,hypertexnames=false,hidelinks]{hyperref}

\numberwithin{equation}{section}

\theoremstyle{plain}
\newtheorem{thm}{Theorem}[section]
\newtheorem{lem}[thm]{Lemma}
\newtheorem{cor}[thm]{Corollary}
\newtheorem{prop}[thm]{Proposition}

\theoremstyle{remark}
\newtheorem{rem}[thm]{Remark}

\newcommand{\T}{\mathbb{T}}
\newcommand{\D}{\mathbb{D}}
\newcommand{\R}{\mathbb{R}}
\newcommand{\C}{\mathbb{C}}
\newcommand{\Z}{\mathbb{Z}}
\newcommand{\E}{\mathcal{E}}
\newcommand{\Var}{\operatorname{Var}}
\newcommand{\Cov}{\operatorname{Cov}}
\newcommand{\av}[1]{\langle #1\rangle}
\newcommand{\dist}{\operatorname{dist}}
\newcommand{\Lip}{\operatorname{Lip}}

\newcommand{\tr}{\operatorname{tr}}
\newcommand{\arccot}{\operatorname{arccot}}
\newcommand{\pv}{\operatorname{p.v.}}

\allowdisplaybreaks[2]

\hypersetup{
  pdftitle={The polygonal Weinstock inequality and quantitative stability},
  pdfauthor={Zhuo Cheng, Yian Zhang}
}

\title[The polygonal Weinstock inequality and quantitative stability]{The polygonal Weinstock inequality and quantitative stability}
\author{Zhuo Cheng}
\address{School of Mathematics and Statistics, HNP-LAMA, Central South University, Changsha, Hunan 410083, P.~R.~China}
\email{zhuo@outlook.cz}

\author{Yian Zhang}
\address{School of Mathematics and Statistics, HNP-LAMA, Central South University, Changsha, Hunan 410083, P.~R.~China}
\email{yian93217@gmail.com}
\date{}


\begin{document}

\begin{abstract}
For every integer $N \ge 3$, the regular $N$-gon $R_N$ uniquely maximizes the perimeter-normalized first positive Steklov eigenvalue $J(P) = |\partial P|\sigma_1(P)$ among all simple planar polygons $P$ with at most $N$ sides, up to similarity. More generally, if $P$ has $p$ strictly convex vertices and $q$ reflex vertices after merging consecutive collinear sides, then $J(P) \le J(R_{p+\min(q,2)})$, with equality if and only if $q = 0$ and $P$ is similar to $R_p$. Among polygons with at most $N$ sides, restricting to those that are nonconvex or have fewer than $N$ genuine sides lowers the supremum of $J$ by $12\pi\zeta(3)N^{-4} + O(N^{-5})$ as $N \to \infty$, with a positive gap for every $N \ge 4$.

For convex $p$-gons, the spectral deficit $J(R_p) - J(P)$ controls the normalized squared deviations of support heights, midpoint offsets, and side lengths with weights $1$, $p^{-2}$, and $p^{-3}$, as well as $p^{-3}$ times the sum of squared exterior-angle deviations from $2\pi/p$; the length weight $p^{-3}$ is sharp. For triangles, the deficit is instead comparable to the square root of the relative area deficit at fixed perimeter, or equivalently to $1 - r$, where $r$ is the singular-value ratio of an affine map from an equilateral triangle. Independently of the number of reflex vertices, a boundary phase identity in arc-length balanced disk coordinates yields the sharp fourth-order bound $2\pi/J(P) - 1 \ge 2\zeta(3)p^{-3} + 8\zeta(4)p^{-4} + O(p^{-5})$. For fixed $p$, a degree-$p$ Blaschke transport links this identity to $H^{1/p}$ coercivity for the regularized square root of the conformal derivative, perimeter-concentration rigidity at equally spaced convex prevertices, and quadratic Hausdorff stability for $p \ge 4$.
\end{abstract}

\maketitle

\section{Introduction}\label{sec:intro}

\subsection{Spectral isoperimetry on polygons and the Weinstock problem}
For a bounded Lipschitz domain $\Omega \subset \R^2$, the Steklov eigenvalue problem
\[
  \Delta u = 0 \quad \text{in }\Omega, \qquad \partial_\nu u = \sigma u \quad \text{on }\partial\Omega,
\]
has a discrete spectrum $0 = \sigma_0(\Omega) < \sigma_1(\Omega) \le \sigma_2(\Omega) \le \cdots \to \infty$ coinciding with the eigenvalues of the Dirichlet-to-Neumann map on $\partial\Omega$ (see \cite{McLean2000, Grisvard2011}). Among all simply connected planar domains of finite perimeter, Weinstock's classical inequality \cite{Weinstock1954} asserts that the scale-invariant functional
\[
  J(\Omega) \coloneqq |\partial\Omega|\,\sigma_1(\Omega)
\]
satisfies $J(\Omega) \le 2\pi$, with equality if and only if $\Omega$ is a disk. Higher Steklov eigenvalues satisfy the sharp Hersch--Payne--Schiffer bounds $|\partial\Omega|\sigma_k(\Omega) \le 2\pi k$ \cite{HPS1975, GP2010}, while at fixed volume the ball maximizes $\sigma_1$ in every dimension \cite{Brock2001}; see \cite{GP2017, CGGS2024} for comprehensive surveys of the Steklov problem, \cite{GLS2016} for quasiconformal bounds, and \cite{LPPS2022} for Steklov asymptotics on curvilinear polygons. Weinstock's theorem occupies the same foundational position for the Dirichlet-to-Neumann operator as the Faber--Krahn inequality and the Szeg\H{o}--Weinberger inequality \cite{Szego1954, Weinberger1956} occupy for the Dirichlet and Neumann Laplacians.

Restricting an elliptic shape-optimization problem from general domains to polygons with at most $N$ sides fundamentally alters its mathematical character. Because Steiner symmetrization generally increases the number of vertices, rearrangement techniques break down past four sides. In the polygonal isoperimetric program initiated by P\'olya and Szeg\H{o} \cite{PolyaSzego1951} (see also \cite{Henrot2006}, and \cite{LS2009, LS2011} for Neumann eigenvalues of triangles and symmetric polygons), the polygonal Faber--Krahn conjecture has long been known unconditionally only for $N = 3$ and $N = 4$. Substantial recent progress has come from the validated computational and local minimality framework of Bogosel and Bucur \cite{BogoselBucur2024, BogoselBucur2026}, the resolution of the convex polygonal Faber--Krahn inequality for all sufficiently large $N$ by Cheng, Gui, Hu, and Li \cite{CGHL2026}, the quantitative P\'olya--Szeg\H{o} theorem for torsional rigidity in the tangential class by Gui, Hu, and Li \cite{GHL2026}, and recent advances on the polygonal Saint--Venant and P\'olya--Szeg\H{o} inequalities by Wang \cite{Wang2026} and by Bogosel, Bucur, and Fragal\`a \cite{BBF2026}.

The analogous question for the Steklov eigenvalue is the \emph{polygonal Weinstock problem} \cite[Open Problem 8]{Nigam2020}: Oudet originally asked which convex $N$-gon maximizes $J(P)$, and the general formulation allowing arbitrary simple polygons with at most $N$ sides was raised at the 2018 AIM workshop on Steklov eigenproblems (see \cite[Section 1]{regular}, as well as the numerical optimization study \cite{DNS2017}). Recently, Cheng, Gui, Hu, Li, and Yao \cite{regular} formulated the regular-polygon conjecture and established both the global strict monotonicity $J(R_N) < J(R_{N+1})$ for all $N \ge 3$ and the rotational-character identification of the first Steklov eigenspace of $R_N$. Independent proofs of the triangle case and the convex quadrilateral case have also been obtained recently by Davies \cite{Davies2026} and by Dong, Yang, and Pan \cite{DYP2026}, respectively.

The Steklov problem on polygons carries two further subtleties that distinguish it sharply from both the smooth Weinstock problem and the Dirichlet problem. First, while quantitative stability holds under a fixed-volume constraint \cite{BDR2012} or within convex sets \cite{GLPT2020} (see also \cite{FL2022}), Bucur and Nahon \cite{BucurNahon2021} showed that in the unrestricted class of smooth simply connected planar domains, Weinstock's perimeter inequality is genuinely unstable: microscopic boundary oscillations can decouple geometric convergence from the behavior of $J(\Omega)$. In dimensions three and higher, Fraser and Schoen \cite{FS2019} (building on their conformal and free-boundary minimal-surface theory \cite{FS2011, FS2016}) showed that the perimeter/surface-area Weinstock inequality fails in general without convexity, whereas Bucur, Ferone, Nitsch, and Trombetti \cite{BFNT2021} proved that convexity restores it. Second, within the class of simple planar $N$-gons, allowing \emph{nonconvex} polygons creates a severe geometric obstruction. Unlike classical area-perimeter isoperimetry, passing to the convex hull does not transparently increase $J(P)$, and inserting reflex vertices can redistribute arc length and deform the domain far away from its convex hull without adding convex corners.

The central principle established in this paper is that, from the conformal perspective, strictly convex vertices and reflex vertices play asymmetric roles. Under a Riemann map $\D \to P$, the Schwarz--Christoffel boundary arc-length density on the unit circle $\T$ acquires integrable power singularities $|1 - \bar{\xi}_i z|^{-2b_i}$ ($0 < b_i < 1/2$) exclusively at the preimages of strictly convex vertices, whereas reflex vertices carry negative turning weights ($-1/2 < b_i < 0$) and merely produce zeros of the conformal derivative. As a result, the Weinstock deficit
\[
  d(P) \coloneqq \frac{2\pi}{J(P)} - 1 > 0
\]
is governed primarily by the number of strictly convex vertices rather than by the total number of sides. Developing this observation into a quantitative theory resolves the polygonal Weinstock problem completely for every $N \ge 3$ and every simple polygon, convex or nonconvex, by purely analytic arguments without computer assistance. At the same time, it reveals that polygonal constraints restore strong quantitative stability to Weinstock's inequality, with sharp anisotropic scaling laws and an exact spectral gap separating nonconvex and fewer-sided polygons from the regular $N$-gon.

\subsection{Main results}
Throughout the paper, $P \subset \R^2$ denotes a simple planar polygon of positive area, and $R_n$ denotes a regular $n$-gon ($n \ge 3$). After merging consecutive collinear sides so that no interior angle equals $\pi$, let $p$ be the number of strictly convex vertices (interior angle strictly below $\pi$) and $q$ the number of reflex vertices (interior angle strictly above $\pi$). Because the exterior turning angles sum to $2\pi$ and each strictly convex vertex contributes less than $\pi$, one always has $p \ge 3$.

Our main theorem sharpens the polygonal Weinstock inequality by linking $J(P)$ directly to $p$ and $\min(q, 2)$.

\begin{thm}[Polygonal Weinstock inequality]\label{thm:main}
Let $P$ be a simple planar polygon with $p$ strictly convex vertices and $q$ reflex vertices after merging consecutive collinear sides, and set $n = p + \min(q, 2)$. Then
\[
  J(P) \le J(R_n),
\]
with equality if and only if $q = 0$ and $P$ is similar to $R_p$. Consequently, for every integer $N \ge 3$, the regular $N$-gon $R_N$ uniquely maximizes $J$ among all simple planar polygons with at most $N$ sides, up to similarity.
\end{thm}

Theorem~\ref{thm:main} exhibits a sharp saturation at $q = 2$: any simple polygon with $p$ strictly convex vertices and at least two reflex vertices satisfies $J(P) < J(R_{p+2})$, no matter how many additional reflex vertices it possesses. Combining this uniform bound with the strict monotonicity and asymptotic expansion of $J(R_n)$ established in \cite{regular} and refined in Section~\ref{sec:regular},
\[
  J(R_N) = 2\pi - 4\pi\zeta(3)N^{-3} + O(N^{-4}), \qquad
  J(R_N) - J(R_{N-1}) = 12\pi\zeta(3)N^{-4} + O(N^{-5}),
\]
shows that nonconvex $N$-gons cannot asymptotically approach $J(R_N)$ any more closely than the regular $(N - 1)$-gon does.

\begin{thm}[Spectral gap for nonconvex and degenerate polygons]\label{cor:defective-gap}
For $N \ge 4$, let $\mathcal{D}_N$ be the class of simple planar polygons with at most $N$ sides that are either nonconvex or have fewer than $N$ genuine sides after merging consecutive collinear edges. There exists an absolute constant $c > 0$ such that
\[
  J(R_N) - J(P) \ge c N^{-4} \qquad \text{for all } P \in \mathcal{D}_N.
\]
More precisely,
\[
  \inf_{P \in \mathcal{D}_N} \bigl(J(R_N) - J(P)\bigr)
  = 12\pi\zeta(3) N^{-4} + O(N^{-5}) \qquad (N \to \infty).
\]
\end{thm}

In contrast with the generic instability of Weinstock's inequality on smooth simply connected domains \cite{BucurNahon2021}, convex $p$-gons satisfy a uniform quantitative stability estimate with explicit, sharp dependence on $p$. Let $Q$ be a convex $p$-gon normalized to unit perimeter, translated so that its boundary centroid lies at the origin. Put $x = 1/p$ and let $H_* = (x/2)\cot(\pi x)$ be the apothem of the unit-perimeter regular $p$-gon $R_p$. On the $i$-th edge, let $y_i$ denote its length, $h_i > 0$ its support height from the origin, and $m_i$ the signed displacement of its midpoint from the foot of the perpendicular dropped from the origin. Define the normalized edge parameters
\[
  a_i = \frac{h_i}{H_*}, \qquad b_i = p y_i, \qquad v_i = \frac{m_i}{H_*},
\]
and write $\mathbb{E} f = \sum_{i=1}^p y_i f_i$ for the edge-length average. Denote by $2\theta_i \in (0, \pi)$ the exterior turning angle at the $i$-th vertex (so that $\sum_{i=1}^p \theta_i = \pi$), and let $\mathcal{I} \ge 0$ be the cubic support-function inertia remainder defined in \eqref{eq:support-remainder}.

\begin{thm}[Quantitative stability for convex polygons]\label{thm:convex}
There exists an absolute constant $c > 0$ such that every unit-perimeter convex $p$-gon $Q$ ($p \ge 3$) satisfies
\[
  J(R_p) - J(Q)
  \ge c\Bigl\{\mathcal{I} + \mathbb{E}\bigl[(a_i - 1)^2 + p^{-2} v_i^2 + p^{-3}(b_i - 1)^2\bigr]\Bigr\}.
\]
In particular, there exists an absolute constant $c' > 0$ such that
\[
  J(R_p) - J(Q) \ge \frac{c'}{p^3}\sum_{i=1}^p\Bigl(\theta_i - \frac{\pi}{p}\Bigr)^2.
\]
Moreover, the weight $p^{-3}$ multiplying the squared side-length deviations is optimal as $p \to \infty$.

For triangles ($p = 3$), the stability exponent is strictly stronger: setting $\mathcal{D} = \mathbb{E}[(a_i - 1)^2 + (b_i - 1)^2 + v_i^2]$ and writing $r \in (0, 1]$ for the ratio of the singular values of an affine map sending an equilateral triangle onto $Q$, one has
\[
  J(R_3) - J(Q) \asymp 1 - r \asymp \sqrt{1 - |Q|/|R_3|} \asymp \sqrt{\mathcal{D}},
\]
where $R_3$ has the same perimeter as $Q$.
\end{thm}

The three weights $(1, p^{-2}, p^{-3})$ in Theorem~\ref{thm:convex} reflect the distinct geometric sensitivities of $J$ near $R_p$: perturbing support heights alters the macroscopic shape at order $1$, shifting edge midpoints acts at order $p^{-2}$, and redistributing edge lengths along a nearly circular contour incurs a penalty of order $p^{-3}$. Splitting a single edge of $R_{p-1}$ and pushing its midpoint outward infinitesimally shows that $p^{-3}$ cannot be improved (Remark~\ref{rem:length-weight-sharp}). On the other hand, triangles exhibit a linear dependence on affine distortion $1 - r$, hence a square-root law in the relative area deficit, rather than the quadratic scaling appearing for $p \ge 4$. As a direct consequence of Theorem~\ref{thm:convex}, among all nonsingular affine images $T R_p$ ($p \ge 3$), $J$ is uniquely maximized when $T$ is a similarity (Corollary~\ref{prop:affine-hand}).

Our conformal machinery yields two further results that require no convexity assumption and place no restriction on the number of reflex vertices $q$. Pulling the normalized boundary arc-length measure $d\mu = ds/|\partial P|$ back to $\T$ via a Riemann map $\D \to P$ produces an absolutely continuous probability measure $\rho\,dm$ on $\T$, where $dm = d\theta/(2\pi)$. In \emph{arc-length balanced coordinates} (uniquely determined up to rotation by $\int_\T z\,\rho(z)\,dm(z) = 0$), the preimages of the $p$ strictly convex vertices divide $\T$ into $p$ open arcs of normalized lengths $s_1, \dots, s_p$ summing to $1$. Write $Z = \zeta(3)$ and $W = \zeta(4)$, and let $H$ be the analytic square root of the conformal derivative $F$ normalized by $|F| = \rho$ on $\T$ and $H(0) > 0$.

\begin{thm}[Sharp fourth-order bound independent of reflex vertices]\label{thm:unrestricted-reflex-asymptotic}
Among all simple planar polygons $P$ with $p$ strictly convex vertices and arbitrarily many reflex vertices,
\[
  \inf_P d(P) = 2Z p^{-3} + 8W p^{-4} + O(p^{-5}) = d(R_p) + O(p^{-5}) \qquad (p \to \infty).
\]
More precisely, there is an absolute constant $C > 0$ such that for all sufficiently large $p$,
\[
  d(P) \ge d(R_p) + 2p^{-4}\sum_{i=1}^p (p s_i - 1)^2 + \frac{1}{6}p^{-3}\|H - 1\|_2^2 - C p^{-5}.
\]
In particular, $J(P) < J(R_{p+1})$ holds for all sufficiently large $p$, uniformly in the number of reflex vertices, and $d(P) \le d(R_p) + O(p^{-5})$ implies both $\sum_{i=1}^p (ps_i - 1)^2 = O(p^{-1})$ and $\|\rho - 1\|_1 = O(p^{-1})$.
\end{thm}

Even without any restriction on the number or angles of reflex vertices, Theorem~\ref{thm:unrestricted-reflex-asymptotic} recovers not only the leading third-order regular term $2\zeta(3)p^{-3}$ but also the exact fourth-order coefficient $8\zeta(4)p^{-4}$, while quantitatively controlling both the spacing of convex prevertices and the $L^2$-deviation of $H$ from $1$.

For fixed $p \ge 3$, a canonical degree-$p$ Blaschke factorization adapted to the convex prevertices governs both perimeter distribution and Hausdorff proximity to $R_p$.

\begin{thm}[Equally spaced prevertices and quadratic Hausdorff stability]\label{thm:intro-blaschke}
Let $P$ be a simple planar polygon of unit perimeter with $p \ge 3$ strictly convex vertices and arbitrarily many reflex vertices, and work in arc-length balanced coordinates.
\begin{enumerate}[label=\textup{(\roman*)},leftmargin=*]
\item If the $p$ convex prevertices are equally spaced on $\T$ and either $p \ge 4$, or $p = 3$ and $P$ admits a nonidentity rotational symmetry, then $J(P) \le J(R_p)$, with equality if and only if $P$ is similar to $R_p$. Moreover, for every $t \in (0, \pi)$, the fraction of $|\partial P|$ whose preimage lies within angular distance $t/p$ of a convex prevertex is minimized by $R_p$, and its excess above the regular value controls the total variation distance between $(z \mapsto z^p)_*\mu$ and the corresponding measure for $R_p$.
\item Let $e^2 = \sum_{i=1}^p (s_i - 1/p)^2$ and $H_p(P) = \inf_g d_{\mathrm{H}}(\partial P, g\partial R_p)$, where $d_{\mathrm{H}}$ denotes Hausdorff distance, $R_p$ has perimeter $1$, and $g$ ranges over rigid motions. For each $p \ge 4$, there exist constants $c_p, K_p > 0$ such that
\[
  d(P) - d(R_p) \ge c_p H_p(P)^2 - K_p e^2.
\]
In particular, there exists $\gamma_p > 0$ such that $d(P) \ge d(R_{p+1}) + \gamma_p$ whenever $p \ge 4$ and $e$ is sufficiently small.
\end{enumerate}
\end{thm}

Stronger weighted versions of Theorem~\ref{thm:intro-blaschke} (including the optimal $H^{1/p}$ Sobolev coercivity of $(FS^{2/p})^{1/2}$, the $(p-3)$-dimensional marked logarithmic density model, and the linear Hausdorff strengthening at $p = 3$) are established in Section~\ref{sec:asymptotics-blaschke} (Lemma~\ref{lem:radial-marker}, Corollary~\ref{cor:arc-order}, Theorem~\ref{thm:symmetric}, and Corollary~\ref{cor:sharp-gap-stability}).

\subsection{Analytic and geometric mechanisms of the proof}
The variational characterization of $J(P)$ is
\begin{equation}\label{eq:variational-physical}
  J(P)\,\Var_\mu(F) \le \int_P |DF|^2, \qquad
  \Var_\mu(F) \coloneqq \int_{\partial P} \Bigl|F - \int_{\partial P} F\,d\mu\Bigr|^2 d\mu,
\end{equation}
valid for all nonconstant scalar- or vector-valued $F \in H^1(P)$, where $d\mu = ds/|\partial P|$ and $|DF|^2 = \sum_k |\nabla F_k|^2$. Equality for a nonconstant harmonic field holds if and only if every centered component is zero or a first Steklov eigenfunction, and $J(P)$ is invariant under similarities. The proofs of Theorems~\ref{thm:main}--\ref{thm:intro-blaschke} combine four distinct analytic and geometric mechanisms.

\medskip
\noindent\emph{1. The boundary phase identity and gap optimization ($q \ge 2$).}
By conformal invariance, \eqref{eq:variational-physical} translates to a spectral problem on $\T$ pairing the circular Dirichlet-to-Neumann energy $\E(f)$ with the $\mu$-variance of $f$. In arc-length balanced coordinates, optimizing radial-tangential perturbations $z(g + ik)$ yields a sharp duality inequality bounding $|\mu g - \bar{g}|^2$ by $d(P)\,\mathcal{N}(g)$ (Lemma~\ref{lem:circle-gram}). On the other hand, the Schwarz--Christoffel derivative $F$ satisfies a boundary phase identity (Lemma~\ref{lem:boundary-phase-product}): for any real $f$ vanishing at every strictly convex prevertex, $-\int_\T (\Lambda + 1)f\,\rho\,dm$ equals a strictly nonnegative double integral measuring boundary variations of the unit tangent phase $u = F/|F|$. Choosing a profile supported on each convex gap $I_i$ and optimizing its Fourier modes via half-strip harmonic extensions bounds $d(P)$ from below solely in terms of the balanced gap lengths $s_1, \dots, s_p$ (Lemma~\ref{lem:gap}). Convexity then collapses the comparison against $R_{p+2}$ to a scalar inequality in $\sum_{i=1}^p s_i^2$, settling all polygons with $q \ge 2$ at once (Proposition~\ref{prop:pure-two-finite}).

\medskip
\noindent\emph{2. Spectral analysis of the regular polygons.}
To separate $P$ from $R_n$ at finite $n$ as well as asymptotically, we refine the regular-polygon Fourier matrix analysis of \cite{regular}. Multiplying the singular weight $w_{1/n}$ by a reciprocal density $w_{-\lambda}$ produces an exact positive semidefinite quadratic form that bounds the compressed matrix norm $\|M_*\| \le (n + 1)/(n - 1)^2$ and delivers tight rational upper bounds for $d(R_5), \dots, d(R_{10})$ (Lemma~\ref{lem:regular-six-short}). Simultaneously, evaluating a two-dimensional half-plane Gaussian chord integral in polar and Cartesian coordinates yields a closed-form positive power series for the coupling energy $S(x)$, from which the Neumann run-length expansion and Euler's beta integral extract the exact fifth-order expansion $d(R_n) = 2Z n^{-3} + 8W n^{-4} + 26\zeta(5)n^{-5} + O(n^{-6})$.

\medskip
\noindent\emph{3. Physical-space inertia, half-strip trials, and equivariant transplantation ($q = 0$).}
On convex polygons, conformal testing alone does not easily reach the exact maximizer $R_p$, whereas physical-space geometry captures the exact cancellation: at fourth order in $1/p$, the coordinate-field deficit matches the regular Steklov deficit. Expressing the boundary moment of inertia through the support function exposes an exact cubic remainder $\mathcal{I} \ge 0$ controlling the exterior angles, and augmenting the coordinate field by an odd harmonic half-strip profile along each edge produces a global a priori budget confining any near-maximizing convex $p$-gon ($p \ge 10$) to an explicit neighborhood of $R_p$.

Inside this neighborhood, we transplant the regular first eigenfield linearly on each centroid sector. Controlling the resulting quadratic form requires a sharp trace bound on the energy complement of the coordinate field $X$ on $R_p$; we prove an equivariant trace inequality via Stokes' theorem and Legendre--Wirtinger estimates (Lemma~\ref{lem:sector-trace}) valid for every $p \ge 3$. For triangles, this trace bound directly yields both rigidity and the linear-in-$(1 - r)$ stability law (Theorem~\ref{thm:pure-triangles}). For $4 \le p \le 9$, we construct explicit low-mode vertex-fixing harmonic fields (exact on the square, three-mode polynomials for $5 \le p \le 9$) whose transplantation forces competitors into the exact-eigenfield stability box, unifying the proof of convex rigidity and Theorem~\ref{thm:convex} across all $p \ge 3$.

\medskip
\noindent\emph{4. Single-reflex comparisons ($q = 1$), cyclic phase averaging, and Blaschke transport.}
When $q = 1$, beating $J(R_{p+1})$ requires resolving the fourth-order term in $1/p$. For $p \ge 6$, we add an interior plateau $g$ to the gap profile $f$ so that $(I + \Lambda)f + kg$ extracts the full within-gap phase integral, combine it with a cross-gap chord bound (Lemma~\ref{lem:cross-chord}), and pass to the nonlinear radial-tangential test $z(\sqrt{1 + 2h} + ik)$ via a circle chain rule obtained by modifying $\Lambda - I$ on the first harmonics to render its difference kernel pointwise nonnegative (Lemma~\ref{lem:markov-square-chain}). Together with pointwise half-cylinder barriers for $\Lambda f$ (Lemma~\ref{lem:strip-phase-pointwise}), a single analytic estimate settles all $p \ge 6$ (Proposition~\ref{prop:phase64}).

The remaining low convex counts divide cleanly: $p = 3$ is settled by a convex-hull pocket decomposition showing $J(P) < 3\sqrt{3} < J(R_4)$ (Proposition~\ref{prop:hull}), while for $p = 4$ and $p = 5$ gap dispersion localizes the convex prevertices near a regular tuple, where pulling back the reciprocal regular density $w_{-1/p}$ through a degree-$p$ Blaschke product $B$ (and controlling its energy via an Euler integral and a carr\'e du champ pushforward bound) contradicts $d(P) \le d(R_{p+1})$. In parallel, averaging boundary phase differences across cyclic gap shifts yields the unrestricted fourth-order bound of Theorem~\ref{thm:unrestricted-reflex-asymptotic}, while the Blaschke transport identity yields the sharp $H^{1/p}$ coercivity, arc-mass rigidity, and quadratic Hausdorff stability of Theorem~\ref{thm:intro-blaschke}.

\subsection{Organization of the paper and conventions}
Section~\ref{sec:gap} develops the balanced conformal coordinates, radial-tangential duality, arc-length discrepancy, boundary phase identity, and gap optimization lemmas. Section~\ref{sec:regular} analyzes the regular-polygon Fourier matrix, proves the strict monotonicity, fifth-order expansion, and quantitative bounds for $d(R_n)$, and establishes the comparison with two additional sides. Section~\ref{sec:convex} treats convex $p$-gons for every $p \ge 3$ and proves Theorem~\ref{thm:convex} together with its affine-image corollary. Section~\ref{sec:asymptotics-blaschke} proves Theorem~\ref{thm:unrestricted-reflex-asymptotic}, develops the Blaschke-transported moment theory, and establishes Theorem~\ref{thm:intro-blaschke} and its weighted refinements. Section~\ref{sec:one-large} proves the single-reflex comparison for $p \ge 6$, and Section~\ref{sec:low-reflex} resolves the remaining single-reflex cases $p \in \{3, 4, 5\}$. Section~\ref{sec:assembly} assembles the proofs of Theorems~\ref{thm:main} and~\ref{cor:defective-gap}.

Throughout all explicit numerical estimates, terminating decimals denote exact rational numbers.

\section{Balanced conformal coordinates and the boundary phase identity}\label{sec:gap}

\subsection{Variational formulation and balanced disk coordinates}
Let $\D = \{z \in \C : |z| < 1\}$ be the open unit disk and $\T = \partial\D$ the unit circle equipped with the normalized Lebesgue probability measure $dm = d\theta/(2\pi)$ at $z = e^{i\theta}$. Pulling the normalized boundary arc-length measure $d\mu = ds/|\partial P|$ back to $\T$ under a Riemann map $\D \to P$ yields an absolutely continuous probability measure $\rho\,dm$ on $\T$, which we continue to denote by $\mu$.

For an integrable function $f$ on $\T$, write $\widehat{f}(k) = \int_\T f(z)\bar{z}^k\,dm(z)$ for its $k$-th Fourier coefficient and $\bar{f} = \int_\T f\,dm$, $\mu f = \int_\T f\,d\mu$. Let $\Lambda$ be the circular Dirichlet-to-Neumann operator defined by $\Lambda(z^k) = |k|z^k$, and set
\[
  \E(f) \coloneqq \sum_{k \in \Z} |k|\,|\widehat{f}(k)|^2,
  \qquad
  \E(f,g) \coloneqq \sum_{k \in \Z} |k|\,\widehat{f}(k)\overline{\widehat{g}(k)}.
\]
Since the harmonic extension of $f$ to $\D$ has Dirichlet energy $2\pi\E(f)$, conformal invariance of the Dirichlet integral rewrites \eqref{eq:variational-physical} in terms of $d(P) = 2\pi/J(P) - 1$ as
\begin{equation}\label{eq:variational-disk}
  1 + d(P) = \sup_{\E(f) > 0} \frac{\Var_\mu(f)}{\E(f)},
\end{equation}
where $\Var_\mu(f) = \int_\T |f - \mu f|^2\,d\mu$.

Let $\alpha_i \in (0,2\pi)$ denote the interior angles at the genuine vertices of $P$, and define the signed exterior turning weights $b_i = (\pi - \alpha_i)/(2\pi)$, so that $\sum_i b_i = 1$. By the Schwarz--Christoffel formula (see \cite{DT2002}), the boundary density has the form
\begin{equation}\label{eq:sc-density}
  \rho(z) = C \prod_i |1 - \bar{\xi}_i z|^{-2b_i},
\end{equation}
where the product runs over the distinct prevertices $\xi_i \in \T$ of all genuine vertices of $P$. Because $\alpha_i > 0$, every vertex satisfies $b_i < 1/2$; hence $p \ge 3$ and $\rho \in L^s(dm)$ for every $s > 1$ satisfying $2s \max_i (b_i)_+ < 1$.

On the quotient Hilbert space $H^{1/2}(\T)/\C$ equipped with the inner product $\E(\cdot,\cdot)$, define the bounded positive self-adjoint operator $T_\mu$ by
\[
  \E(T_\mu f, g) = \Cov_\mu(f,g) \coloneqq \int_\T (f - \mu f)\overline{(g - \mu g)}\,d\mu.
\]
By the compact embedding of mean-zero $H^{1/2}(\T)$ into $L^{2s/(s-1)}(dm)$ and H\"older's inequality, $T_\mu$ is compact, and its operator norm is $\lambda_P \coloneqq 1 + d(P)$.

We say that disk coordinates are \emph{balanced for a probability measure $\nu$ on $\T$} if $\int_\T z\,d\nu(z) = 0$. We will use this normalization both for arc-length measures and for finite sets of prevertices.

\begin{lem}[Existence, uniqueness, and continuity of balanced coordinates]\label{lem:balanced-coords}
Let $\nu$ be a probability measure on $\T$ whose atoms each have mass strictly less than $1/2$. Then there exists a unique point $a \in \D$ minimizing the functional
\[
  \Phi_\nu(a) \coloneqq \int_\T \log\frac{|z - a|^2}{1 - |a|^2}\,d\nu(z).
\]
This minimizer is the unique solution in $\D$ of the balancing equation
\begin{equation}\label{eq:balance-eq}
  \int_\T \frac{z - a}{1 - \bar{a} z}\,d\nu(z) = 0,
\end{equation}
and depends continuously on $\nu$ under weak convergence to any probability measure with all atomic masses below $1/2$. In particular, arc-length balanced coordinates for $P$ exist and are unique up to rotation.
\end{lem}

\begin{proof}
To verify coercivity near $\partial\D$, suppose $a_k \to \zeta \in \T$. Choose a closed arc $I \subset \T$ centered at $\zeta$ with $\nu(I) = m < 1/2$. For $z \in I$ we have $|z - a_k| \ge 1 - |a_k|$, whereas $|z - a_k|$ stays uniformly bounded below on $\T \setminus I$ for large $k$. Hence
\[
  \Phi_\nu(a_k) \ge (2m - 1)\log(1 - |a_k|) - O(1) \longrightarrow +\infty,
\]
so $\Phi_\nu$ attains its global minimum in $\D$.

For strict convexity along hyperbolic geodesics, apply a disk automorphism mapping a given geodesic to a diameter; this shifts each integrand by an additive constant independent of the hyperbolic distance $\ell$ along the geodesic, transforming the integrand into $\log(\cosh\ell - \cos\theta\,\sinh\ell)$. Its second derivative with respect to $\ell$ is $\sin^2\theta / (\cosh\ell - \cos\theta\,\sinh\ell)^2$, which is strictly positive away from the two endpoints $\theta \in \{0,\pi\}$. Because these two endpoints together carry $\nu$-mass strictly less than $1$, $\Phi_\nu$ is strictly convex along every hyperbolic geodesic, yielding uniqueness of the minimizer and of solutions to the critical equation \eqref{eq:balance-eq}.

Finally, if $\nu_k \rightharpoonup \nu$ weakly and every atom of $\nu$ has mass below $1/2$, a finite covering of $\T$ by closed arcs of $\nu$-mass strictly less than $1/2$ yields the same uniform lower bound on $\Phi_{\nu_k}$, confining the minimizers $a_k$ to a compact subdisk of $\D$. Passing to the limit in \eqref{eq:balance-eq} along any convergent subsequence gives $a_k \to a$ by uniqueness.
\end{proof}

Retaining in counterclockwise cyclic order only the preimages of the $p$ strictly convex vertices of $P$, we denote the $p$ open circular arcs (\emph{gaps}) separating adjacent convex prevertices by $I_1, \dots, I_p$, and set $s_i = m(I_i)$, so that $\sum_{i=1}^p s_i = 1$. Throughout the paper we set $Z = \zeta(3)$, $W = \zeta(4)$, and define the constants and functions
\begin{equation}\label{eq:gap-constants}
  c = \frac{1}{2Z}, \qquad b = \frac{W}{Z^2}, \qquad e = \frac{\zeta(6)}{Z^2}, \qquad
  D(s) = c + bs - es^3, \qquad \mathcal{H}(s) = \frac{s^4}{D(s)}.
\end{equation}

\subsection{Radial-tangential duality and arc-length discrepancy}
Perturbing the balanced coordinate map $z$ by a radial-tangential vector field produces an energy-optimal duality inequality controlling linear functionals of $\mu - m$ in terms of the Weinstock deficit.

\begin{lem}[Radial amplitude and tangential correction]\label{lem:circle-gram}
For real mean-zero functions $g, k \in H^{1/2}(\T)$, set $V = z(g + ik)$ and define $Cg \in H^{1/2}(\T)$ by $\widehat{Cg}(0) = 0$ and $\widehat{Cg}(n) = i\widehat{g}(n)/n$ for $n \neq 0$. In arc-length balanced coordinates,
\begin{align*}
  |\mu(g + ik)|^2 + d(P)\Var_\mu(V) &\le d(P)\bigl(1 + d(P)\bigr)\E(V), \\
  \E(V) &= \mathcal{N}(g) + \E(k - Cg),
\end{align*}
where
\[
  \mathcal{N}(g) \coloneqq \sum_{n \neq 0} \bigl(|n| - |n|^{-1}\bigr)|\widehat{g}(n)|^2.
\]
In particular, $(\mu g)^2 + d(P)\Var_\mu(zg) \le d(P)(1 + d(P))\E(g)$.
\end{lem}

\begin{proof}
By \eqref{eq:variational-disk}, the Hermitian form $(1+d)\E - \Cov_\mu$ (with $d = d(P)$) is positive semidefinite. On the pair $(z, V)$, its $(1,1)$-entry is $(1+d) - 1 = d$, while its $(1,2)$-entry is $-\overline{\mu(g+ik)}$ because $\mu z = 0$ and $\widehat{V}(1) = \widehat{g}(0) + i\widehat{k}(0) = 0$. Taking the determinant of the $2 \times 2$ Gram matrix yields the first inequality.

To evaluate $\E(V)$, pair the Fourier modes at frequencies $1 + n$ and $1 - n$ for each $n \ge 1$. Using $\widehat{g}(-n) = \overline{\widehat{g}(n)}$ and $\widehat{k}(-n) = \overline{\widehat{k}(n)}$, the sum of these two terms in $\E(V)$ equals
\[
  (n+1)|\widehat{g}(n) + i\widehat{k}(n)|^2 + (n-1)|\overline{\widehat{g}(n)} + i\overline{\widehat{k}(n)}|^2
  = 2\Bigl(n - \frac{1}{n}\Bigr)|\widehat{g}(n)|^2 + 2n\Bigl|\widehat{k}(n) - \frac{i\widehat{g}(n)}{n}\Bigr|^2.
\]
Summing over $n \ge 1$ gives the energy identity. Setting $k = 0$ recovers $\E(zg) = \E(g)$ and hence the final inequality.
\end{proof}

Choosing $k = Cg$ minimizes $\E(V)$ for a given radial profile $g$. Since $1$ and $z$ are orthonormal in $L^2(\mu)$, decomposing $V = z(g + iCg)$ along $\operatorname{span}\{1,z\}$ in Lemma~\ref{lem:circle-gram} gives
\[
  |\mu(g + iCg)|^2 + \frac{d(P)}{1 + d(P)}\dist_{L^2(\mu)}\bigl(z(g + iCg), \operatorname{span}\{1,z\}\bigr)^2
  \le d(P)\,\mathcal{N}(g).
\]
This estimate has an exact spectral interpretation in terms of the operator $T = T_\mu$. Writing $\lambda = 1 + d(P)$ and $Sg = z(g + iCg)$ for real mean-zero $g$, define $g_\rho$ by $\widehat{g}_\rho(n) = \widehat{\rho}(n)/(|n| - |n|^{-1})$ for $|n| \ge 2$ and $\widehat{g}_\rho(n) = 0$ for $|n| \le 1$. Balance gives $\E(Tz, z) = 1$ and $\widehat{Tz}(n) = \widehat{\rho}(n-1)/|n|$ for $n \neq 0$. Because $\widehat{Sg}(1+n) = (1 - 1/n)\widehat{g}(n)$ for $n \neq 0$, we have $(T - I)z = Sg_\rho$ modulo constants. Defining the spectral moment
\[
  B_\mu \coloneqq \E\bigl((T - I)z\bigr) = \mathcal{N}(g_\rho)
  = \sum_{n \ge 2} \frac{2n}{n^2 - 1}|\widehat{\rho}(n)|^2,
\]
we obtain the identity
\begin{equation}\label{eq:spectralmom}
  d(P) - B_\mu = \E\bigl(T(\lambda I - T)z, z\bigr) \ge 0.
\end{equation}
Here $\E(Tz) = 1 + B_\mu$. Because $\rho > 0$ almost everywhere, $T$ is injective on $H^{1/2}(\T)/\C$, so equality holds in \eqref{eq:spectralmom} if and only if $Tz = \lambda z$, which is equivalent to $\rho \equiv 1$. Moreover, since $-2\log|z - w|$ is the convolution kernel of $\Lambda^{-1}$ on mean-zero functions, pairing $z\rho$ with $\Lambda^{-1}(z\rho)$ (justified via Abel limits in $L^{s/(s-1)}$ since $\rho \in L^s$) yields the chord representation
\[
  1 + B_\mu = \E(Tz) = -2\iint_{\T \times \T} \Re(z\bar{w})\log|z - w|\,d\mu(z)\,d\mu(w).
\]
For any real mean-zero $g \in H^{1/2}(\T)$, one has $\E(Sg, (T - I)z) = \mu(g + iCg)$, so Cauchy--Schwarz yields the sharp duality bound
\begin{equation}\label{eq:duality}
  |\mu g - \bar{g}|^2 \le B_\mu\,\mathcal{N}(g) \le d(P)\,\mathcal{N}(g) \le d(P)\,\E(g),
\end{equation}
valid for every real $g \in H^{1/2}(\T)$. The leading constant $B_\mu$ is attained at $g = g_\rho$, for which $\mu g_\rho = B_\mu$.

As an immediate application of \eqref{eq:duality}, $B_\mu$ controls the discrepancy between $\mu$ and $m$ uniformly over all arcs of $\T$.

\begin{lem}[Arc-length discrepancy in balanced coordinates]\label{lem:arc-smoothing}
For $0 < h \le 1/24$, set $a = \sin(\pi h/2)$ and
\[
  C_h = \pi^{-2}\Bigl[\log\frac{1}{a} - \frac{5}{8} + 2a^2 - a^4\Bigr].
\]
Then every arc $I \subset \T$ satisfies $|\mu(I) - m(I)| \le h + \sqrt{C_h B_\mu}$. In particular, $C_{1/48} < 0.294$ and $C_{1/72} < 0.334$. As $B_\mu \downarrow 0$, one has the uniform bound
\[
  \sup_I |\mu(I) - m(I)| \le \Bigl(\frac{1}{\pi\sqrt{2}} + o(1)\Bigr)\sqrt{B_\mu \log(1/B_\mu)}.
\]
\end{lem}

\begin{proof}
Work in turn coordinates on $\R/\Z$. Define a random variable $U \in [-h/2, h/2]$ by $\sin(\pi U) = a\cos\theta$, where $\theta$ is uniform on $[0,\pi]$. Its circular logarithmic potential is constant on $[-h/2, h/2]$: by symmetry and the arcsine integral (obtained by factoring the difference of cosines and using $\int_0^\pi \log\sin t\,dt = -\pi\log 2$), for every $|v| \le h/2$,
\[
  \mathbb{E}\log|\sin\pi(v - U)|
  = \frac{1}{2}\mathbb{E}\log|\sin^2\pi v - \sin^2\pi U|
  = \mathbb{E}\log|\sin\pi v - a\cos\theta| = \log\frac{a}{2}.
\]
The first two cosine moments of $U$ are $k_1 = \mathbb{E}\cos(2\pi U) = 1 - a^2$ and $k_2 = \mathbb{E}\cos(4\pi U) = (1 - a^2)(1 - 3a^2)$.

Let $I$ be an arc of length $s = m(I) > h$, and let $g$ be the convolution of the indicator of the concentric arc of length $\ell = s - h$ with the law of $U$. Then $0 \le g \le \mathbf{1}_I$ almost everywhere and $\bar{g} = \ell$. Let $U'$ be an independent copy of $U$, set $X = U - U'$, and define
\[
  K(t) = \log|\sin\pi t| + \cos(2\pi t) + \frac{1}{8}\cos(4\pi t).
\]
In $\mathcal{N}(g) = \E(g) - \|g - \bar{g}\|_{\dot{H}^{-1/2}}^2$, retain only frequencies $\pm 1$ and $\pm 2$ in the subtracted norm. Using the Fourier coefficients of the interval indicator (justified by Abel limits since the potential is continuous; for $t = |\sin\pi v| > a$ the same arcsine integral gives $\log((t + \sqrt{t^2 - a^2})/2)$), we obtain
\[
  \pi^2 \mathcal{N}(g) \le \mathbb{E}\bigl[K(\ell + X) - K(X)\bigr]
  \le \max K - \log\frac{a}{2} - k_1^2 - \frac{k_2^2}{8}.
\]
In terms of $u = \sin^2\pi t$, we have $K(t) = \frac{1}{2}\log u + \frac{9}{8} - 3u + u^2$, which achieves its maximum at $u_* = (3 - \sqrt{5})/4$. From $4u_*^2 - 6u_* + 1 = 0$ and $\log(4u_*) \le 4u_* - 1$, it follows that $\max K \le \frac{1}{2} - \log 2 - (\sqrt{5} - 2)/8$. Since $a < 1/15$, we have $k_2^2 \ge 1 - 8a^2 > 3 - \sqrt{5}$, and hence $\mathcal{N}(g) \le C_h$.

Applying \eqref{eq:duality} yields $\mu(I) \ge \mu g \ge s - h - \sqrt{C_h B_\mu}$, a bound that holds trivially when $s \le h$. Applying the same estimate to $\T \setminus I$ gives the two-sided discrepancy bound. The numerical bounds at $h = 1/48, 1/72$ follow from $a \ge 3h/2$, $2a^2 < 0.01$, $\log 2 < 0.7$, $\log 3 < 1.1$, and $\pi^2 > 9.85$. Choosing $h = \sqrt{B_\mu}$ yields the asymptotic statement.
\end{proof}

\subsection{The boundary phase identity and gap optimization}
We next connect the modulus of the Schwarz--Christoffel derivative to its boundary phase. Consider a zero-free holomorphic product on $\D$ of the form
\[
  F(z) = C \prod_j (1 - \bar{\xi}_j z)^{-2b_j}, \qquad C \neq 0, \quad \xi_j \in \T\ \text{distinct}, \quad b_j < \frac{1}{2},
\]
with analytic branches in $\D$. Set $\beta = \sum_j b_j$, and define the boundary phase $u = F/|F|$ on $\T$ away from the finitely many points $\xi_j$ (a constant global phase factor is immaterial).

\begin{lem}[Boundary phase identity]\label{lem:boundary-phase-product}
Choose $s \in (1,\infty)$ such that $2s\max_j (b_j)_+ < 1$, and let $s' = s/(s-1)$ be its conjugate exponent. Set $H = F^{1/2}$ and $D_z H(w) = (H(w) - H(z))/(w - z)$. For every real $f \in W^{1,s'}(\T)$ vanishing at every point $\xi_j$ with $b_j > 0$, one has the absolutely convergent identities
\begin{equation}\label{eq:phase-product}
\begin{aligned}
  -\int_\T (\Lambda + \beta)f\,|F|\,dm
  &= 2\int_\T f(z)\|D_z H\|_2^2\,dm(z) \\
  &= \iint_{\T \times \T} f(z)|F(w)|\left|\frac{u(w) - u(z)}{w - z}\right|^2 dm(w)\,dm(z).
\end{aligned}
\end{equation}
\end{lem}

\begin{proof}
Set $\rho = |H|^2 = |F|$, let $P_r$ denote the Poisson extension to radius $r \in (0,1)$, and assume first that $f \ge 0$. Because $2s\max_j (b_j)_+ < 1$, we have $H \in H^{2s}(\D)$ and $P_r H(z) = H(rz)$. Expanding $|H(w) - H(z)|^2$ and using the self-adjointness of $P_r$ gives
\begin{multline*}
  \iint_{\T \times \T} f(z)\frac{r|H(w) - H(z)|^2}{|w - rz|^2}\,dm(w)\,dm(z) \\
  = \frac{r}{1 - r^2}\left[\int_\T \rho(P_r f - f)\,dm
    - 2\Re\int_\T f \bar{H}\bigl(H(rz) - H(z)\bigr)\,dm\right].
\end{multline*}
As $r \uparrow 1$, the kernel $r|w - rz|^{-2} = (|w - z|^2 + (1-r)^2/r)^{-1}$ increases monotonically to $|w - z|^{-2}$. On the right-hand side, since $\Lambda f \in L^{s'}(\T)$ and $r\partial_r P_r f = P_r \Lambda f$, we have $(P_r f - f)/(1 - r) \to -\Lambda f$ in $L^{s'}(\T)$, while $\rho \in L^s(\T)$.

Off the prevertices, $(H(rz) - H(z))/(1 - r) \to -zH'(z)$. To pass to the limit under the integral by dominated convergence, set $t = |z - \xi_j|$ near a prevertex $\xi_j$ of weight $b_j \neq 0$. From the product formula and $|1 - r\bar{\xi}_j z|^2 \ge r t^2$, we have $|H(z)| \lesssim t^{-b_j}$ and $\sup_{1/2 < r < 1} |H'(rz)| \lesssim 1 + t^{-b_j - 1}$. When $b_j > 0$, the vanishing $f(\xi_j) = 0$ and H\"older's inequality yield $|f(z)| \lesssim t^{1/s}$ with $1/s > 2b_j$; when $b_j < 0$, boundedness of $f$ suffices. In both cases the integrand admits an integrable majorant. Because $2\Re(\bar{H} z H') = -\beta \rho$ on $\T$, taking $r \uparrow 1$ yields the first equality in \eqref{eq:phase-product}. Decomposing a signed $f$ into $f_+ - f_-$ establishes absolute convergence in general.

For the second equality, fix $z \in \T$ away from the prevertices. Then $D_z H \in H^2(\D)$, and writing $\Delta H(w) = H(w) - H(z)$ gives
\[
  \int_\T \frac{\Delta H(w)^2}{|w - z|^2}\,dm(w)
  = -z\int_\T w\bigl(D_z H(w)\bigr)^2 dm(w) = 0,
\]
since $w(D_z H(w))^2 \in H^1(\D)$ vanishes at $0$. Consequently, for every unimodular constant $\alpha \in \T$,
\begin{equation}\label{eq:hardy-projection}
  \int_\T \frac{|\Delta H(w) - \alpha\overline{\Delta H(w)}|^2}{|w - z|^2}\,dm(w) = 2\|D_z H\|_2^2.
\end{equation}
Choosing $\alpha = u(z) = H(z)/\overline{H(z)}$ makes the numerator equal to $|H(w) - u(z)\overline{H(w)}|^2 = |F(w)|\,|u(w) - u(z)|^2$, which completes the proof.
\end{proof}

In particular, Lemma~\ref{lem:boundary-phase-product} applies whenever $f$ and $\Lambda f$ are continuous on $\T$ and $f$ vanishes at every strictly convex prevertex, because boundedness of the circular Hilbert transform then places $f'$ in $L^r(\T)$ for all $1 < r < \infty$.

Using \eqref{eq:phase-product} together with \eqref{eq:duality}, we now optimize a gap-supported trial function to bound $d(P)$ solely in terms of the balanced gap lengths $s_1, \dots, s_p$.

\begin{lem}[Gap profile optimization]\label{lem:gap}
Recall $D(s) = c + bs - es^3$ and $\mathcal{H}(s) = s^4/D(s)$ from \eqref{eq:gap-constants}, and define
\[
  \mathcal{H}_{\mathrm{opt}}(s) \coloneqq 2s^4 \sum_{j=1}^\infty \frac{1}{j^3 + 2s j^2 - 2s^3}.
\]
In arc-length balanced coordinates, put $y = 1/p$ and $Q = \sum_{i=1}^p (s_i/y - 1)^2$, and normalize $H_F^2 = F$ by $|F| = \rho$ and $H_F(0) > 0$. Then
\[
  d(P) \ge B_\mu \ge \Bigl(1 + \frac{1 - H_F(0)^2}{2}\Bigr)^2 \sum_{i=1}^p \mathcal{H}_{\mathrm{opt}}(s_i),
  \qquad \mathcal{H}_{\mathrm{opt}}(s) \ge \mathcal{H}(s).
\]
Moreover, setting $m = \sum_{i=1}^p s_i^2 = y(1 + Q/p)$, both $K = \mathcal{H}$ and $K = \mathcal{H}_{\mathrm{opt}}$ satisfy
\begin{equation}\label{eq:gap}
  \sum_{i=1}^p K(s_i) \ge \frac{K(m)}{m}
  \ge p\Bigl(1 + \frac{Q}{p}\Bigr)^2 K(y)
  \ge (p + 2Q) K(y).
\end{equation}
In particular, $B_\mu \ge (c\nu^3 + b\nu^2 - e)^{-1}$ where $\nu = (\sum_{i=1}^p s_i^2)^{-1} \le p$.
\end{lem}

\begin{proof}
Consider first finitely many nonnegative coefficients $A_{ij}$. Parametrizing each gap $I_i$ by its unit coordinate $v \in [0,1]$, set $f(v) = \sum_j A_{ij}(1 - \cos(2\pi jv))$ on $I_i$. Then $f \ge 0$ is piecewise smooth with Lipschitz derivative, and both $f$ and $f'$ vanish at every convex prevertex; thus $\Lambda f$ is continuous. Applying \eqref{eq:phase-product} with $|F| = \rho$ and $\beta = 1$ to $G = (I + \Lambda)f$ yields
\[
  -\mu G = 2\int_\T f(z)\|D_z H_F\|_2^2\,dm(z).
\]
Since $|w - z| \le 2$ on $\T$, writing $|H_F(w) - H_F(z)|^2 = |H_F(w) - H_F(0)|^2 + |H_F(z) - H_F(0)|^2 - 2\Re\bigl((H_F(w) - H_F(0))\overline{(H_F(z) - H_F(0))}\bigr)$ and integrating in $w$ (using $\int_\T H_F\,dm = H_F(0)$ and $\|H_F\|_2 = 1$) gives
\[
  2\|D_z H_F\|_2^2 \ge \frac{1}{2}\bigl(\|H_F - H_F(0)\|_2^2 + |H_F(z) - H_F(0)|^2\bigr)
  \ge \frac{1 - H_F(0)^2}{2}.
\] Hence $-\mu G \ge \frac{1}{2}(1 - H_F(0)^2)\bar{f}$. Setting $L = 1 + (1 - H_F(0)^2)/2$ and invoking \eqref{eq:duality}, we obtain
\[
  L^2\Bigl(\sum_{i,j} s_i A_{ij}\Bigr)^2 \le (\bar{G} - \mu G)^2 \le B_\mu\,\mathcal{N}(G).
\]

On the Fourier side, angular differentiation $f'$ satisfies
\[
  \mathcal{N}(G) = \E(f') + 2\|f'\|_2^2 - 2\Var_m(f) - \|f\|_{\dot{H}^{-1/2}}^2.
\]
To bound $\E(f')$ from above, view $\D$ conformally as a half-cylinder of angular circumference $2\pi$ and vertical coordinate $t > 0$. On the half-strip above gap $I_i$, extend the $j$-th sine mode $(j A_{ij}/s_i)\sin(2\pi j v)$ of $f'$ by the factor $e^{-jt/s_i}$. Because these extensions vanish on every vertical ray above a convex prevertex, they glue continuously to a finite-energy competitor on the half-cylinder, whose normalized Dirichlet energy bounds $\E(f')$ from above. The $j$-th mode on $I_i$ contributes $\frac{s_i}{2}(j/s_i)^3 A_{ij}^2$ to this energy. Combining this upper bound with the exact value of $\|f'\|_2^2$ and the within-gap variances (dropping the between-gap variance and $\|f\|_{\dot{H}^{-1/2}}^2$) yields
\[
  \mathcal{N}(G) \le \sum_{i,j} \frac{j^3 + 2s_i j^2 - 2s_i^3}{2s_i^2} A_{ij}^2.
\]
Every denominator coefficient is positive since $0 < s_i < 1$. Choosing $A_{ij} = s_i^3/(j^3 + 2s_i j^2 - 2s_i^3)$ maximizes the resulting quotient by Cauchy--Schwarz, and letting the frequency cutoff tend to infinity gives $B_\mu \ge L^2\sum_i \mathcal{H}_{\mathrm{opt}}(s_i)$.

For later sections, it is convenient to isolate the single profile minimizing the leading strip cost:
\begin{equation}\label{eq:phi-profile}
  \Phi(v) \coloneqq \frac{1}{Z}\sum_{j=1}^\infty \frac{1 - \cos(2\pi jv)}{j^3}, \qquad v \in [0,1].
\end{equation}
By Cauchy--Schwarz, $(\sum_j A_j)^2 \le Z\sum_j j^3 A_j^2$ with equality when $A_j \propto j^{-3}$. Splitting the series at $j \approx 1/v$ shows $\Phi(v) = O(v^2(1 + |\log v|))$ near $v = 0$, and the corresponding split for $\Phi'$ shows that $\Phi'$ is H\"older continuous of every exponent below $1$. Hence the preceding estimates extend by truncation to $f|_{I_i} = a_i \Phi$:
\begin{equation}\label{eq:phi-energy-bounds}
  \|f\|_{\dot{H}^{3/2}}^2 \le c\sum_{i=1}^p \frac{a_i^2}{s_i^2}, \qquad
  2\|f'\|_2^2 = b\sum_{i=1}^p \frac{a_i^2}{s_i}, \qquad
  2\Var_m(f) \ge e\sum_{i=1}^p a_i^2 s_i.
\end{equation}
We will also require a lower bound for $\E(f)$: restricting the harmonic extension of $f$ to each half-strip and relaxing the continuity constraint across vertical rays produces the Neumann extension, which damps the $j$-th cosine mode by $e^{-jt/s_i}$ at cost $jA_{ij}^2/2$. Summing over $j$ gives
\begin{equation}\label{eq:phi-neumann-lower}
  \E(f) \ge 2c^2\zeta(5)\sum_{i=1}^p a_i^2.
\end{equation}

Testing the quotient for $\mathcal{H}_{\mathrm{opt}}(s)$ with $A_j \propto j^{-3}$ immediately gives $\mathcal{H}_{\mathrm{opt}}(s) \ge \mathcal{H}(s)$. More generally, if $\gamma_i \ge 0$ satisfy $\bar{G} - \mu G \ge \sum_{i=1}^p a_i s_i \gamma_i$ for all $a_i \ge 0$, then setting $a_i = s_i^3 \gamma_i / D(s_i)$ in \eqref{eq:duality} yields
\begin{equation}\label{eq:gap-opt}
  d(P) \ge B_\mu \ge \sum_{i=1}^p \mathcal{H}(s_i)\gamma_i^2.
\end{equation}

Finally, assign gap $i$ probability $s_i$, so that the random gap length $S$ has expectation $\mathbb{E} S = m = \sum_i s_i^2$. Each summand $K_0(s) = s^4/D_0(s)$ appearing in $\mathcal{H}$ or $\mathcal{H}_{\mathrm{opt}}$ has a positive concave denominator $D_0$. Applying Cauchy--Schwarz and Jensen's inequality (using convexity of $s^{3/2}$ and concavity of $D_0$) gives
\[
  \sum_{i=1}^p K_0(s_i) = \mathbb{E}\!\left[\frac{S^3}{D_0(S)}\right]
  \ge \frac{(\mathbb{E} S^{3/2})^2}{\mathbb{E} D_0(S)}
  \ge \frac{m^3}{D_0(m)} = \frac{K_0(m)}{m},
\]
with strict inequality unless all gaps are equal. Summing over the terms of $K$ establishes the first inequality in \eqref{eq:gap}. Since $D_0(s) - s D_0'(s) > 0$, the function $s/D_0(s)$ is increasing, whence $K(m)/m \ge (m/y)^2 K(y)/y$, proving the rest of \eqref{eq:gap}. Substituting $m = 1/\nu$ in $\mathcal{H}(m)/m$ yields the final estimate.
\end{proof}

When a gap contains at most one reflex prevertex, the phase quotient in \eqref{eq:phase-product} stays close to $1$ throughout the gap and yields an explicit lower bound on the within-gap phase contribution.

\begin{lem}[Within-gap source average with at most one reflex prevertex]\label{lem:single-gap-loss}
For $0 < s \le U < 1/2$, set $v = (1 - 1/(4U))_+$ and define
\[
  \kappa(U) \coloneqq \left[1 + 2\pi U v^2 - \frac{U}{\pi Z}\bigl(1 - \cos(2\pi v)\bigr)\right]^{-1}.
\]
On a gap of length $s$ containing at most one reflex prevertex, for every fixed target point $w$ in that gap, the average of the phase quotient $|u(w) - u(z)|^2/|w - z|^2$ over the source $z$ with respect to the probability density $\Phi$ in gap coordinates is at least $\kappa(U)$. In particular, $\kappa(1/4) = 1$, $\kappa(1/3) > 0.95$, and $\kappa(5/12) > 9/11$.
\end{lem}

\begin{proof}
At angular separation $t \in (0,1/2)$ turns, the phase quotient equals $1$ whenever $z$ and $w$ lie on the same straight edge. Across a single reflex corner of turning magnitude $\gamma \in (0,\pi)$, the phase changes by $2\pi t + \gamma$ in turns-normalized form, giving
\[
  \left|\frac{u(w) - u(z)}{w - z}\right|^2
  = \frac{\sin^2(\pi t + \gamma/2)}{\sin^2(\pi t)} \ge \min(1, \cot^2 \pi t).
\]
Thus only pairs separated by $t > 1/4$ can experience a loss below $1$. Measuring $w' \in [0,1]$ from the endpoint farther from the target and setting $v_s = (1 - 1/(4s))_+$, we have $2\pi(t - 1/4)_+ \le X \coloneqq 2\pi s(v_s - w')_+$. Since $1 - \cot^2(\pi/4 + \theta) = 2\sin(2\theta)/(1 + \sin(2\theta))$, the pointwise loss is at most $2X/(1 + X)$. Setting $L_s = \int_0^{v_s} \Phi(w')\,dw'$ and $A_s = 2\int_0^{v_s} X\Phi(w')\,dw'$, Jensen's inequality on $[0,v_s]$ bounds the total lost fraction above by
\[
  e(s) \coloneqq \frac{2L_s A_s}{2L_s + A_s} \le \frac{A_s}{1 + A_s},
\]
with $e(s) = 0$ for $s \le 1/4$. Here $L_s \le 1/2$ by symmetry of $\Phi$, and both $L_s$ and $A_s$ increase with $s$. Termwise integration of \eqref{eq:phi-profile} gives
\[
  A_U = 2\pi U v^2 - \frac{U}{\pi Z}\sum_{n=1}^\infty \frac{1 - \cos(2\pi n v)}{n^5}.
\]
Retaining only the $n = 1$ term bounds $A_s \le A_U \le \kappa(U)^{-1} - 1$, yielding the stated lower bound $\kappa(U)$. At $U = 1/3$ and $U = 5/12$, this upper bound on $A_U$ equals $\pi/24 - 1/(3\pi Z) < 0.05$ and $2\pi/15 - 5(1 + \cos(\pi/5))/(12\pi Z) < 2/9$, respectively.

Retaining each within-gap contribution $(1 - e(s_i))\mu_i$ (where $\mu_i = \mu(I_i)$) in the phase identity and invoking \eqref{eq:gap-opt} yields
\begin{equation}\label{eq:four-mass-hand}
  d(P) \ge \sum_{i=1}^p \mathcal{H}(s_i)\bigl[1 + (1 - e(s_i))\mu_i\bigr]^2.
\end{equation}
At $s = 3/10$, keeping the first two terms of the series gives $A_s \le \pi/60 - 21/(128\pi Z) < 0.01$; at $s = 5/12$, direct integration gives $L_s \le 2/5 - [\sin(\pi/5) - (W - 1)]/(2\pi Z) < 1/3$. Consequently $e(3/10) < 0.01$, $e(1/3) < 0.05$, and $e(5/12) < 1/6$.
\end{proof}

\section{Spectral analysis of the regular polygons}\label{sec:regular}

Write $d_n = d(R_n) = 2\pi/J(R_n) - 1$, $Z = \zeta(3)$, and $W = \zeta(4)$. In this section we establish the strict monotonicity of $J(R_n)$, its fifth-order asymptotic expansion
\begin{equation}\label{eq:regular-two}
  d_n = 2Z n^{-3} + 8W n^{-4} + 26\zeta(5)n^{-5} + O(n^{-6}),
\end{equation}
and explicit rational upper bounds for $d_n$ across all $n \ge 4$, building on the regular-polygon Fourier matrix introduced in \cite{regular}. As an immediate application of Lemma~\ref{lem:gap} and these spectral bounds, we then prove the comparison with two additional sides, settling Theorem~\ref{thm:main} for all polygons with at least two reflex vertices.

Throughout our numerical bounds we use the following elementary zeta estimates: for any real $s \ge 3$, comparing the tail $\sum_{k=9}^\infty k^{-s}$ against $\int k^{-s}\,dk$ by the trapezoid and midpoint rules places the tail between $9^{1-s}/(s-1) + 9^{-s}/2$ and $(17/2)^{1-s}/(s-1)$. In particular,
\[
  1.202 < Z < 1.2022, \qquad 1.082 < W < 1.0825, \qquad 1.036 < \zeta(5) < 1.0371.
\]

\subsection{Fourier matrix, first eigenmodes, and strict monotonicity}
For $b < 1/2$, let $w_b$ denote the probability density on $\T$ proportional to $|1 - z|^{-2b}$. Its Fourier coefficients are given by
\[
  r_0(b) = 1, \qquad r_j(b) = \frac{(b)_j}{(1 - b)_j} \quad (j \ge 1),
  \qquad (b)_j \coloneqq b(b+1)\cdots(b+j-1),
\]
where the recurrence $(j + 1 - b)r_{j+1}(b) = (j + b)r_j(b)$ follows from integrating the derivative of $(\sin(\theta/2))^{1-2b} e^{i(j+1/2)\theta}$ over $[0,2\pi]$.

In centered disk coordinates on $R_n$, the boundary density is $w_x(z^n)$ with $x = 1/n$. On the closed rotational subspace of $H^{1/2}(\T)/\C$ spanned by the energy-orthonormal basis $\{z^{1+ni}/\sqrt{\ell_i} : i \in \Z\}$, where $\ell_i = |1 + ni|$, the boundary variance operator $T_\mu$ is represented by the infinite symmetric matrix
\[
  M_{ik} = \frac{r_{|i-k|}(x)}{\sqrt{\ell_i \ell_k}}, \qquad i, k \in \Z.
\]
To identify $\|M\|$ with $1 + d_n$, we first record a general nondegeneracy property of first Steklov eigenfunctions on arbitrary polygons.

\begin{lem}[Nonvanishing first Fourier mode]\label{lem:first-coordinate}
In any disk coordinate representation of $P$ with boundary measure $\mu = \rho\,dm$, every nonconstant real first Steklov eigenfunction $u$ normalized by $\mu u = 0$ satisfies $\E(z, u) \neq 0$.
\end{lem}

\begin{proof}
Write $\lambda = 1 + d(P)$ and suppose for contradiction that $\E(z, u) = 0$. The eigenvalue equation $\lambda \E(z, u) = \mu(z u)$ then yields $\mu(z u) = 0$. Put $a = \int_\T u\,dm$ and $g = u - a$. Pairing opposite Fourier modes of the real mean-zero function $g$ gives $\E(zg) = \E(g) = \E(u)$, whereas
\[
  \Var_\mu(zg) = \Var_\mu(u) + a^2(1 - |\mu z|^2).
\]
Since $\Var_\mu(u) = \lambda\E(u)$ and $|\mu z| < 1$, the variational maximality of $\lambda$ forces $a = 0$. Thus $zu$ is also a maximizer in the complex energy space with $\mu(zu) = 0$. Comparing the two eigenvalue equations $\lambda \Lambda(zu) = \rho z u$ and $\lambda \Lambda u = \rho u$ yields $[\Lambda, z]u = 0$. However, the commutator acts on monomials by $[\Lambda, z]z^k = (|k+1| - |k|)z^{k+1}$, whose multiplier is $\pm 1$ for every $k \in \Z$, contradicting $u \neq 0$.
\end{proof}

Applying a disk automorphism to move an arbitrary interior point to the origin shows that the gradient map $u \mapsto \nabla u(0)$ is injective on the real first eigenspace in any disk coordinate. In centered coordinates on $R_n$, this map is equivariant with respect to the dihedral group $D_n$, so its nonzero image in $\R^2$ is invariant under rotation by $2\pi/n$ and therefore equals all of $\R^2$. Consequently, the real first eigenspace of $R_n$ is two-dimensional and carries the standard representation of $D_n$; moreover, the equivariant first eigenfields form a one-dimensional real space, because any $2 \times 2$ real matrix commuting with both a rotation of order $n \ge 3$ and a reflection is scalar. By Lemma~\ref{lem:first-coordinate}, the orthogonal projection of $z$ onto the first eigenspace is nonzero and lies in the rotational subspace spanned by $\{z^{1+ni}\}$. Therefore $\|M\| = 1 + d_n$.

Every entry $M_{ik}$ decreases with $n$, strictly so whenever $(i,k) \neq (0,0)$. For $n + 1$, replacing the entries of a unit maximizing eigenvector of $M^{(n+1)}$ by their absolute values preserves or increases the Rayleigh quotient, and the eigenvalue equation then forces every entry to be strictly positive. Testing $M^{(n)}$ against this positive unit vector yields $1 + d_n > 1 + d_{n+1}$, establishing the strict monotonicity
\begin{equation}\label{eq:reg-monotonicity}
  J(R_n) < J(R_{n+1}) \qquad (n \ge 3).
\end{equation}

\begin{rem}[Recovering the eigenvalue from iterates of $z$]
For an arbitrary polygon $P$, define recursively $f_0(z) = z$ and $f_{k+1} = T_\mu f_k$, represented in mean-zero form by
\[
  (T_\mu f)(z) = -2\int_\T \log|z - w|\bigl(f(w) - \mu f\bigr)\,d\mu(w),
\]
since the Fourier coefficients of $-2\log|z - w|$ are $1/|j|$ for $j \neq 0$ and $0$ at $j = 0$. Let $\nu$ be the spectral measure of $z$ for $T_\mu$ in the energy inner product, with moments $m_j = \int t^j\,d\nu(t)$. Then $\Var_\mu(f_k)/\E(f_k) = m_{2k+1}/m_{2k}$, which increases with $k$ by Cauchy--Schwarz. Since Lemma~\ref{lem:first-coordinate} guarantees that $\nu$ has a positive atom at $1 + d(P)$, dominated convergence after dividing by $(1 + d(P))^{2k}$ implies $m_{2k+1}/m_{2k} \uparrow 1 + d(P)$. In particular, every strict lower bound for $2\pi/J(P)$ is witnessed by finitely many iterations of $T_\mu$ on $z$.
\end{rem}

\subsection{Reciprocal-density compression, chord integral, and asymptotics}
Let $M_*$ denote the submatrix of $M$ obtained by deleting the zeroth row and column, and let $v$ be the coupling column with entries $v_i = r_{|i|}(x)/\sqrt{\ell_i}$ for $i \neq 0$.

\paragraph{A quadratic form from the reciprocal density.}
Fix $x = 1/n$ and a parameter $\lambda \in [x, 1/2)$. Set $r_m = r_m(x)$, $c_m = -r_m(-\lambda) > 0$, and $b_m = -r_m(x - \lambda) \ge 0$ for $m \ge 1$; the recurrence for $r_m$ remains valid at negative parameters. The pointwise product $w_x w_{-\lambda} = D w_{x-\lambda}$ has constant mean $D = 1 - 2\sum_{m=1}^\infty r_m c_m > 0$. Comparing nonconstant Fourier coefficients (justified by absolute convergence since $c_m = O(m^{-1-2\lambda})$) shows that the positive vector $h_i = \sqrt{\ell_i}\,c_{|i|}$ ($i \neq 0$) satisfies
\begin{equation}\label{eq:schur-row-ratio}
  \frac{(M_* h)_i}{h_i} = \frac{r_{|i|} + D b_{|i|}}{\ell_i c_{|i|}}, \qquad i \neq 0.
\end{equation}
For any finitely supported vector $f$ indexed by $\Z \setminus \{0\}$, expanding the difference yields the positive semidefinite identity
\[
  \sum_{i \neq 0} \frac{(M_* h)_i}{h_i} |f_i|^2 - \langle M_* f, f \rangle
  = \frac{1}{2}\sum_{i,k \neq 0} M_{ik} h_i h_k \left|\frac{f_i}{h_i} - \frac{f_k}{h_k}\right|^2 \ge 0.
\]
Specializing to $\lambda = x$ gives $b_m = 0$ and $M_* h = v$. For $i = \pm m$, the larger row ratio is $Q_m = v_{-m}/h_{-m}$, which satisfies
\[
  Q_1 = \frac{x(1+x)}{(1-x)^2}, \qquad
  \frac{Q_{m+1}}{Q_m} = \frac{(m+x)(m+1+x)}{(m+1-x)^2} < 1
\]
for $n \ge 4$, because $(m+1-x)^2 - (m+x)(m+1+x) = m(1 - 4x) + 1 - 3x > 0$. By the quadratic form bound above, we obtain
\begin{equation}\label{eq:compression}
  \|M_*\| \le Q_1 = \frac{n + 1}{(n - 1)^2} \qquad (n \ge 4).
\end{equation}
By codimension-one min--max, the second variance eigenvalue of $T_\mu$ in the rotational class is likewise bounded above by $(n+1)/(n-1)^2$.

\paragraph{A chord integral and the scalar coupling series.}
For scalar estimates we extend $M_*$ and $v$ to continuous $x \in (0, 1/6]$ by the same formulas, with $\ell_i = |1 + i/x|$. Eliminating the nonzero components in the top eigenvector equation at $x = 1/n$ yields the Schur complement equation
\[
  d_n = \bigl\langle v, ((1 + d_n)I - M_*)^{-1} v \bigr\rangle.
\]
We compare $d_n$ with the zero-shift resolvent $U(x) \coloneqq \langle v, (I - M_*)^{-1} v \rangle$, whose leading term is governed by the scalar series
\begin{equation}\label{eq:scalar-positive}
  S(x) \coloneqq \sum_{j=1}^\infty \frac{r_j(x)^2}{j}
  = \sum_{k=2}^\infty (2^{k-1} - 1)\zeta(k + 1)x^k, \qquad 0 \le x < \frac{1}{2}.
\end{equation}
To derive \eqref{eq:scalar-positive}, let $a, b < 1/2$, $s > \max(-1/2, a+b-1)$, and write $(c)_s = \Gamma(c+s)/\Gamma(c)$. We first evaluate the circular chord moment
\begin{equation}\label{eq:circle-chord-moment}
  I_{ab}(s) \coloneqq \iint_{\T \times \T} |z - w|^{2s} w_a(z) w_b(w)\,dm(z)\,dm(w)
  = 4^s \frac{(1/2)_s (1 - a - b)_s}{(1 - a)_s (1 - b)_s}.
\end{equation}
Equip two independent copies of the upper half-plane with the normalized densities proportional to $v^{-2a} e^{-u^2 - v^2}\,du\,dv$ and $(v')^{-2b} e^{-(u')^2 - (v')^2}\,du'\,dv'$, and compute the expectation of $|\det(X, Y)|^{2s}$ in two ways. In polar coordinates, $z = e^{2i\arg X}$ has density $w_a$ and $|X|^{2s}$ has mean $(1-a)_s$ (and similarly for $Y, w$); since $2|\det(X, Y)| = |X|\,|Y|\,|z - w|$, the expectation equals $4^{-s}(1-a)_s(1-b)_s I_{ab}(s)$. In Cartesian coordinates, conditioning on the heights $v, v'$ and rotating the horizontal Gaussian pair $(u, u')$ reduces $uv' - u'v$ to a single standard Gaussian multiplied by $\sqrt{v^2 + (v')^2}$, contributing $(1/2)_s(v^2 + (v')^2)^s$. Setting $v^2 = r t$ and $(v')^2 = r(1 - t)$ shows that $r = v^2 + (v')^2$ has normalized density $r^{-a-b} e^{-r}\,dr/\Gamma(1 - a - b)$, whose $s$-th moment is $(1 - a - b)_s$. Equating the two evaluations proves \eqref{eq:circle-chord-moment}.

By the Fourier series of $-\log|z - w|$ and Abel limits (justified since small positive and negative chord moments dominate $|\log|z-w||$), differentiating \eqref{eq:circle-chord-moment} at $s = 0$ via Euler's gamma product yields
\[
  2\sum_{j=1}^\infty \frac{r_j(a) r_j(b)}{j}
  = -I'_{ab}(0)
  = \sum_{j=1}^\infty \left(\frac{1}{j - a - b} - \frac{1}{j - a} - \frac{1}{j - b} + \frac{1}{j}\right).
\]
Setting $a = b = x$ gives $S(x) = \sum_{j=1}^\infty x^2/[j(j-x)(j-2x)]$, and expanding in geometric series proves \eqref{eq:scalar-positive}. Integrating the partial fractions produces the integral representations
\begin{equation}\label{eq:regular-log-energy}
  2S(x) = \int_0^1 \frac{(t^{-x} - 1)^2}{1 - t}\,dt, \qquad
  2S(1/p) = p\int_0^1 \frac{v^{p-3}(1 - v)}{1 + v + \cdots + v^{p-1}}\,dv.
\end{equation}

\paragraph{Neumann expansion and the fifth-order coefficient.}
In the Neumann series for $U(x) = \langle v, (I - M_*)^{-1} v \rangle$, grouping consecutive runs of any repeated nonzero index $i$ replaces its diagonal weight by the geometric sum $D_i = \sum_{r=0}^\infty \ell_i^{-1-r} = (\ell_i - 1)^{-1}$, so that $D_j = x/j$ and $D_{-j} = x/(j - 2x)$ for $j \ge 1$. This gives $U(x) = \sum_{k=0}^\infty U_k(x)$, where
\[
  U_k(x) = \sum_{\substack{i_0, \dots, i_k \neq 0 \\ i_a \neq i_{a-1}}}
  r_{|i_0|}(x) r_{|i_k|}(x) \prod_{a=1}^k r_{|i_a - i_{a-1}|}(x) \prod_{a=0}^k D_{i_a}.
\]
By \eqref{eq:compression}, this series converges at $x = 1/6$ and hence throughout $(0, 1/6]$ by entrywise monotonicity. Every summand in $U_k$ has a Taylor series in $x$ with nonnegative coefficients starting in degree $3 + 2k$, so coefficientwise summation is justified on $(0, 1/6]$.

Pairing opposite indices $i_0 = \pm j$ in $Y(x) \coloneqq U_0(x)$ yields
\[
  Y(x) = 2x S(x) + 2x^2 \sum_{j=1}^\infty \frac{r_j(x)^2}{j(j - 2x)} = 2Z x^3 + 8W x^4 + O(x^5),
\]
since $S(x) = Z x^2 + 3W x^3 + O(x^4)$ and the second sum begins with $W x^2$. In particular, $U(x) = Y(x) + O(x^5)$. At $x = 1/n$, comparing the spectral integrals of $U(1/n)$ and $d_n$ against the spectral measure of $M_*$ (supported on $[0, \|M_*\|]$) inserts the factor $d_n/(1 - t)$ into the integrand for $d_n$, yielding
\begin{equation}\label{eq:regular-resolvent-error}
  d_n^2 \le U(1/n) - d_n \le \frac{d_n^2}{1 - \|M_*\|} = O(n^{-6}).
\end{equation}

To extract the fifth-order coefficient, set $V = \zeta(5)$, $P_4(x) = 2Z x^3 + 8W x^4$, and define the nonnegative-coefficient power series
\begin{equation}\label{eq:regular-path-five}
  \mathcal{U}(x) \coloneqq \frac{U(x) - P_4(x)}{x^5}, \qquad 0 < x \le \frac{1}{6}.
\end{equation}
Only $Y = U_0$ and $U_1$ contribute to $\mathcal{U}(0)$. Write $H_m = \sum_{k=1}^m 1/k$ and $T = \sum_{j=1}^\infty H_{j-1}/j^4$. Substituting $r_j(x) = (x/j)[1 + (2H_{j-1} + 1/j)x + O(x^2)]$ into the paired formula for $Y$ gives $[x^5]Y = 22V + 8T$. On the left-hand side of Euler's beta identity
\begin{equation}\label{eq:euler-beta-moments}
  a\int_0^1 t^{a-1}(1 - t)^b\,dt = \frac{\Gamma(1+a)\Gamma(1+b)}{\Gamma(1+a+b)},
\end{equation}
the coefficient of $b^2$ is $a\sum_{j=1}^\infty H_{j-1}/[j(j+a)]$ by the power series of $\frac{1}{2}\log^2(1-t)$. Extracting the coefficient of $a^3 b^2$ on both sides via Euler's gamma product gives $T = 2V - \zeta(2)Z$.

For $U_1$, group the index pairs $(j, j-k)$ by their difference $k \ge 1$. Telescoping $k/[j(j-k)] = 1/(j-k) - 1/j$ gives $\sum_{j \in \Z \setminus \{0,k\}} 1/[j(j-k)] = 2/k^2$, and squaring that identity produces
\[
  \sum_{j \in \Z \setminus \{0,k\}} \frac{1}{j^2(j - k)^2} = \frac{4\zeta(2)}{k^2} - \frac{6}{k^4}.
\]
Multiplying by $2/k$ and summing over $k \ge 1$ yields $[x^5]U_1 = 8\zeta(2)Z - 12V$. Consequently $\mathcal{U}(0) = [x^5](Y + U_1) = 26V$, and \eqref{eq:regular-resolvent-error} transfers this coefficient directly to $d_n$, proving \eqref{eq:regular-two}.

\subsection{Explicit bounds for the reference spectrum}
We now record the quantitative bounds on $J(R_n)$ and $d_n$ needed in the finite-$p$ comparisons, beginning with the explicit first eigenfunctions of the square (compare \cite{Tan2017, GLPS2019} for the Steklov spectrum of rectangles and cuboids).

\begin{lem}[The reference square]\label{lem:reference-square}
Let $\tau \in (0,\pi/2)$ be the unique root of $\tan\tau\tanh\tau = 1$, and set $t = \tanh\tau = \cot\tau$. On the square $S = (-1,1)^2$, the functions $\sin(\tau x)\cosh(\tau y)$ and $\sin(\tau y)\cosh(\tau x)$ are first Steklov eigenfunctions, and
\[
  \frac{11}{15} < t < 0.735, \qquad 0.93 < \tau < 0.94, \qquad 5.5 < J(R_4) = 8\tau t < 5.52.
\]
\end{lem}

\begin{proof}
By the gradient isomorphism established after Lemma~\ref{lem:first-coordinate}, the first eigenfunction on $S$ with gradient $(1,0)$ at the origin is odd in $x$. Replacing it by its absolute value on $\{x > 0\}$ and extending oddly to $S$ preserves the Rayleigh quotient, yielding a first eigenfunction whose boundary inner product with $u(x,y) = \sin(\tau x)\cosh(\tau y)$ is strictly positive. Because $u$ is harmonic and satisfies $\partial_\nu u = \tau t\,u$ on all four sides, orthogonality of distinct eigenspaces identifies $\sigma_1(S) = \tau t$, and swapping $x$ and $y$ produces the second eigenfunction.

Rewriting $\arctan t + \operatorname{artanh} t = \pi/2$ as $\pi/4 = t S_0(t^4)$, where $S_0(u) = \sum_{k=0}^\infty u^k/(4k+1)$, we bound the tail past $u^2/9$ at $u = 0.29$ between $u^3/13$ and $u^3/[13(1-u)]$ to obtain $1.0692 < S_0(0.29) < 1.0700$. Since $(11/15)^4 < 0.29 < 0.735^4$, monotonicity gives $11/15 < t < 0.735$. The function $z \mapsto z\arccot z$ increases on $(0,1)$; expressing $\arccot z = \pi/4 + \arctan w$ with $w = (1-z)/(1+z)$ and applying $w - w^3/3 < \arctan w < w$ at the endpoints $w = 2/13$ and $w = 53/347$ yields $11/16 < t\tau < 0.69$ and $0.93 < \tau < 0.94$.
\end{proof}

\begin{lem}[Reciprocal-density upper bounds for $d_n$]\label{lem:regular-six-short}
For $d_n = d(R_n)$ and $U(x) = \langle v, (I - M_*)^{-1} v \rangle$, one has
\[
  d_5 < 0.0518, \qquad d_6 < 0.02420, \qquad d_7 < 0.013273, \qquad
  d_8 < 0.00812, \qquad d_9 < 0.00531,
\]
and $U(1/10) < 0.003668$. Moreover, the fifth-order remainder $\mathcal{U}(x)$ in \eqref{eq:regular-path-five} satisfies
\begin{equation}\label{eq:regular-chord}
  \mathcal{U}(x) < 27 + 130x \quad (0 < x \le 1/10), \qquad
  \mathcal{U}(1/8) < 44.97, \qquad \mathcal{U}(1/9) < 42.21.
\end{equation}
\end{lem}

\begin{proof}
For $n \in \{5, 6, 7, 10\}$, set $x = 1/n$, $\lambda = \min\{1/3, 1/2 - x\}$, $\delta = \lambda - x$, and $t_m = r_m c_m$. Putting $z = 2m + 1$, $\alpha = 1 - 2x$, $\beta = 1 + 2\lambda$, and $u_m = m + 1/2 + \delta$, the two recurrences give
\[
  \frac{t_{m+1}}{t_m}
  = \frac{z^2 + \alpha\beta - z(\alpha + \beta)}{z^2 + \alpha\beta + z(\alpha + \beta)}
  \ge \frac{z - \alpha - \beta}{z + \alpha + \beta}
  = 1 - \frac{2 + 2\delta}{u_{m+1}}.
\]
Rearranging yields $u_m t_m - u_{m+1} t_{m+1} \le (1 + 2\delta)t_m$. Since $u_m t_m \to 0$, summing from $m = 3$ gives
\[
  D \le D_0 \coloneqq 1 - 2t_1 - 2t_2 - \frac{7 + 2\delta}{1 + 2\delta} t_3.
\]
By \eqref{eq:schur-row-ratio}, the reciprocal-density quadratic form yields the operator bound $M_* \le \operatorname{diag}(T_{|i|}/\ell_i)$, where $T_m = (r_m + D_0 b_m)/c_m$.

The sequence $T_m$ is concave: for either ratio $v_m = r_m/c_m$ (with $\theta = x$) or $v_m = b_m/c_m$ (with $\theta = x - \lambda$), the forward differences are positive and satisfy
\[
  \frac{v_{m+2} - v_{m+1}}{v_{m+1} - v_m}
  = \frac{(m + \theta)(m + 1 + \lambda)(2m + 3)}{(m + 2 - \theta)(m + 1 - \lambda)(2m + 1)} \le 1,
\]
because the denominator exceeds the numerator by $2[(1 - 2\theta - 2\lambda)(m+1)^2 + \lambda(1 - \theta)] \ge 0$.

For each of the four parameter choices tabulated below, direct evaluation checks that $2T_3 \le 3T_2$, $T_1 < (n-1)a$, and $T_2 < (2n-1)a$. Concavity then implies that $T_m/m$ decreases for $m \ge 2$, and hence so does
\[
  C_m \coloneqq \frac{1}{na - T_m/m + a/m} + \frac{1}{na - T_m/m - a/m},
\]
as seen by expanding in an even geometric series in $(a/m)/(na - T_m/m)$. Inverting the diagonal upper bound on $M_*$ and pairing opposite indices yields
\[
  \bigl\langle v, (aI - M_*)^{-1} v \bigr\rangle
  \le C_3 S(x) + \sum_{m=1}^2 (C_m - C_3)\frac{r_m^2}{m}.
\]
To bound $S(x)$ rationally, put $\gamma = 15/64 - x^2/2$ and $g(u) = [2(u^2 + \gamma)]^{-1}$. For $u = j - x$ with $j \ge 3$, one has $(\gamma + 1/4)^2 < u^2/32$, which ensures $g(u - 1/2) - g(u + 1/2) \ge [u(u^2 - x^2)]^{-1}$. Telescoping \eqref{eq:scalar-positive} therefore gives
\[
  S(x) \le x^2\left[\sum_{j=1}^2 \frac{1}{j(j - x)(j - 2x)} + g(5/2 - x)\right].
\]

Evaluating these rational upper bounds for $n = 5, 6, 7, 10$ and weighting $C_1, C_2, C_3$ by the nonnegative quantities $r_1^2$, $r_2^2/2$, and $S - r_1^2 - r_2^2/2$, respectively, yields the following bounds:
\begin{itemize}[leftmargin=*]
\item For $n = 5$ and $a = 1.0518$: $(C_1, C_2, C_3) < (0.5621, 0.531, 0.522)$, $S < 0.09392$, and $\langle v, (aI - M_*)^{-1} v \rangle < 0.05166$.
\item For $n = 6$ and $a = 1.02420$: $(C_1, C_2, C_3) < (0.4316, 0.4110, 0.4037)$, $S < 0.05700$, and $\langle v, (aI - M_*)^{-1} v \rangle < 0.02419$.
\item For $n = 7$ and $a = 1.013273$: $(C_1, C_2, C_3) < (0.351929, 0.33531, 0.329034)$, $S < 0.038304$, and $\langle v, (aI - M_*)^{-1} v \rangle < 0.0132724$.
\item For $n = 10$ and $a = 1$: $(C_1, C_2, C_3) < (0.229, 0.220, 0.217)$, $S < 0.01619$, and $\langle v, (aI - M_*)^{-1} v \rangle < 0.003668$.
\end{itemize}
For example, at $n = 7$ the tail estimate gives $D_0 < 0.8875$ with $r_1^2 = 1/36$ and $r_2^2/2 = 8/1521$; at $n = 10$ one has $(c_1, c_2, c_3) = (1/4, 1/14, 1/28)$ and $0.926 < D_0 < 0.927$, so the resolvent is at most $0.217(0.01619) + 0.012/81 + 0.003(11/171)^2/2 < 0.003668$. For $n = 5, 6, 7$, the upper bound lies strictly below $a - 1$, so the Schur complement equation gives $d_n < a - 1$; at $n = 10$ it yields $U(1/10) < 0.003668$ directly.

By \eqref{eq:regular-resolvent-error} and \eqref{eq:compression}, we also have $U(1/7) < 0.013273 + (9/7)(0.013273)^2 < 0.0135$. Subtracting $P_4(x)$ using $Z > 1.202$ and $W > 1.082$ gives $\mathcal{U}(1/10) < 40$ and $\mathcal{U}(1/7) < 48.51$. Because $\mathcal{U}(x)$ has nonnegative Taylor coefficients, it is convex on $[0, 1/7]$ with $\mathcal{U}(0) = 26\zeta(5) < 27$. Its chords on $[0, 1/10]$ and $[1/10, 1/7]$ (where $1/8$ and $1/9$ divide $[1/10, 1/7]$ in fractions $7/12$ and $7/27$) immediately give \eqref{eq:regular-chord}. Substituting into $U = P_4 + x^5\mathcal{U}$ and using $d_n(1 + d_n) \le U(1/n)$ yields $d_8 < 0.00812$ and $d_9 < 0.00531$.
\end{proof}

\begin{cor}\label{cor:regular-uniform-nine}
Set $R'_p = p^4 d(R_{p+1}) - 2Zp$ and $A = 8W - 6Z$. Then
\[
  R'_6 < 2.78, \qquad R'_p < A + \frac{8.9}{p} \quad (p \ge 4).
\]
In particular, $R'_7 < 2.741$, $R'_8, R'_9 < 2.60$, and $R'_p < 2.34$ for all $p \ge 10$.
\end{cor}

\begin{proof}
The bound $d_7 < 0.013273$ from Lemma~\ref{lem:regular-six-short} directly gives $R'_6 < 2.78$. For $4 \le p \le 8$, each bound for $d_{p+1}$ in Lemma~\ref{lem:regular-six-short} is strictly less than $2.404/p^3 + 1.442/p^4 + 8.9/p^5$.

For $p \ge 9$, put $y = 1/p$, $x = y/(1+y)$, and $B = 2Z y^3 + A y^4 + 8.9 y^5$. By \eqref{eq:regular-resolvent-error}, it suffices to verify $B(1 + B) > U(x)$. Expanding $x^{-5}[B(1+B) - P_4(x)]$ as a polynomial in $y$, its first three coefficients are
\[
  32W - 12Z + 8.9, \qquad 20Z + 10A + 44.5 + 4Z^2, \qquad 10Z + 10A + 89 + 20Z^2 + 4ZA,
\]
which exceed $29$, $88$, and $150$, respectively, while all higher coefficients are nonnegative. On the other hand, \eqref{eq:regular-chord} gives $(1 + y)\mathcal{U}(x) < 27 + 157y$. Multiplying the difference by $1 + y$ yields at least
\[
  2 - 40y + 238y^2 = 2(1 - 10y)^2 + 38y^2 > 0.
\]
The specific numerical bounds for $p \ge 7$ follow from $A < 1.448$.
\end{proof}

\subsection{Comparison with two additional sides}
Combining Lemma~\ref{lem:gap} with Corollary~\ref{cor:regular-uniform-nine} resolves the comparison against $R_{p+2}$ for all polygons, with an explicit quantitative margin in terms of the balanced gap lengths.

\begin{prop}[Comparison with two additional sides]\label{prop:pure-two-finite}
Let $s_1, \dots, s_p$ be the balanced convex-vertex gaps of a simple polygon $P$, and set $m = \sum_{i=1}^p s_i^2$ and $f(s) = s^3/D(s) = \mathcal{H}(s)/s$. For every integer $k \ge 3$,
\[
  d(P) - d(R_{k+2}) > \frac{1}{(k+1)^4} + f(m) - f(1/k).
\]
In particular, $m \ge 1/k$ implies $J(P) < J(R_{k+2})$, independently of the number of reflex vertices. Setting $k = p$ and $x = 1/p$ yields
\[
  d(P) - d(R_{p+2}) > \frac{1}{(p+1)^4} + \frac{x^2(3c + 2bx)}{D(x)^2}\sum_{i=1}^p (s_i - x)^2.
\]
\end{prop}

\begin{proof}
By \eqref{eq:gap}, we have $d(P) \ge f(m)$, so the first assertion reduces to proving $d_* - d(R_{k+2}) > (k+1)^{-4}$ for the polygon-independent quantity $d_* = f(1/k)$. Set $y = 1/(k+1)$, $u = 2 - 2W/Z$, $h = 1 + u$, $\ell = 1 - 2u$, and $A = 8W - 6Z$. Writing $1/k = y/(1-y)$ gives $d_* = 2Z y^3/L$, where
\[
  L = 1 - hy - \ell y^2 - \Bigl(\frac{2\zeta(6)}{Z} + u - 1\Bigr)y^3 < 1 - hy - \ell y^2.
\]
By Corollary~\ref{cor:regular-uniform-nine} (applied at $p = k + 1 \ge 4$), it suffices to verify $[2Z + (A + 1)y + 8.9 y^2]L < 2Z$. Subtracting the left-hand side from $2Z$, dividing by $y$, applying $L < 1 - hy - \ell y^2$, and dropping the resulting positive cubic term leaves a lower bound of
\[
  12(Z - W) - 1 + \bigl[2Z\ell + (A + 1)h - 8.9\bigr]y + \bigl[(A + 1)\ell + 8.9h\bigr]y^2
  > 0.43 - 4.54y + 12.1y^2 > 0,
\]
where the numerical bounds follow from $1.202 < Z < 1.2022$ and $1.082 < W < 1.0825$ (which give $0.198 < u < 0.200$ and $A > 1.442$), and the final quadratic has negative discriminant. This establishes $d_* - d(R_{k+2}) > y^4$.

Finally, $f(s) = s^3/D(s)$ is increasing and convex on $(0,1)$ because $D$ is positive and concave (Lemma~\ref{lem:gap}). Bounding $f(m) - f(x)$ from below by its tangent at $x = 1/p$, whose slope is
\[
  f'(x) = \frac{x^2\bigl(3D(x) - xD'(x)\bigr)}{D(x)^2} = \frac{x^2(3c + 2bx)}{D(x)^2},
\]
and noting that $m - x = \sum_{i=1}^p (s_i - x)^2$ gives the final inequality.
\end{proof}

As an illustration, for any $p \ge 4$, the single dispersion condition $\sum_{i=1}^p (s_i - 1/p)^2 \ge 1/[p(p-1)]$ already forces $m \ge 1/(p-1)$ and hence $d(P) > d(R_{p+1}) + p^{-4}$ by taking $k = p - 1$ in Proposition~\ref{prop:pure-two-finite}, without using any information about the reflex angles.

\section{Convex polygons: global comparison and quantitative stability}\label{sec:convex}

In this section we prove the sharp comparison $J(Q) \le J(R_p)$ with uniqueness of equality for every convex $p$-gon ($p \ge 3$), together with the quantitative stability theorem (Theorem~\ref{thm:convex}). The argument proceeds in physical space rather than conformal coordinates. First, a cubic support-function identity and a tangential half-strip correction yield an a priori geometric budget controlling area, support heights, and side lengths. Next, a sharp equivariant trace inequality on $R_p$ governs the energy-variance balance of the linearly transplanted regular eigenfield on a concrete box around $R_p$. Combining these tools settles $p \ge 10$ directly, $p = 3$ via affine pullback, and $4 \le p \le 9$ via explicit low-mode vertex-fixing harmonic fields, after which quantitative stability follows uniformly across all $p \ge 3$.

\subsection{Boundary inertia and the tangential half-strip trial}\label{sec:convex-strip}
Normalize the convex $p$-gon $Q$ to unit perimeter and place its boundary centroid at the origin. Set $x = 1/p$, $H = 2|Q|$, $H_* = (x/2)\cot(\pi x)$, and $I = \int_{\partial Q} |X|^2\,ds$. In the outward unit normal and positive unit tangent frame $(\nu_i, \tau_i)$ of the $i$-th edge, we parametrize the edge by
\[
  X_i(t) = h_i \nu_i + (m_i + y_i t)\tau_i, \qquad -\frac{1}{2} \le t \le \frac{1}{2},
\]
where $y_i$ is the edge length, $h_i > 0$ is its support height from the origin, and $m_i$ is the signed displacement of its midpoint from the foot of the perpendicular. Write $\mathbb{E} f = \sum_{i=1}^p y_i f_i$. The divergence theorem gives $\mathbb{E} h_i = H$, and integration along each edge gives $I = \mathbb{E}(h_i^2 + m_i^2 + y_i^2/12)$. Testing \eqref{eq:variational-physical} with the coordinate field $X$ (which has Dirichlet energy $H$ and boundary variance $I$) yields $j I \le H$ for every $0 < j \le J(Q)$. The polygonal isoperimetric inequality gives $H \le H_*$.

Let $u(\phi) = \max_{Z \in \overline{Q}} Z \cdot (\cos\phi, \sin\phi)$ be the support function of $Q$. On each open interval of normal directions corresponding to a vertex, $u'' = -u$, whereas across an edge normal $u'$ jumps from $m_i - y_i/2$ to $m_i + y_i/2$, so $\int_0^{2\pi} u(\phi)\,d\phi = \sum_i y_i = 1$. Similarly, $u^2 u' + \frac{2}{3}(u')^3$ has derivative $-u^3$ between edge normals and jumps by $\int_{m_i - y_i/2}^{m_i + y_i/2} (h_i^2 + 2t^2)\,dt$ across the $i$-th edge normal. Summing over all edges yields the cubic inertia identity
\[
  2I - \mathbb{E} h_i^2 = \int_0^{2\pi} u(\phi)^3\,d\phi.
\]
Between the $i$-th and $(i+1)$-th edge normals, denote the exterior turning angle by $2\theta_i \in (0,\pi)$ (so $\sum_{i=1}^p \theta_i = \pi$) and write $u = c_0\cos t + d_0\sin t$ for $|t| \le \theta_i$. From the endpoint values $h_i$ and $h_{i+1}$, we have $w_i \coloneqq \int_{-\theta_i}^{\theta_i} u\,dt = 2c_0\sin\theta_i = (h_i + h_{i+1})\tan\theta_i$ and $2d_0\sin\theta_i = h_{i+1} - h_i$. Define
\[
  k(t) \coloneqq \frac{2\sin t}{\sqrt{1 - \frac{1}{3}\sin^2 t}}, \qquad k_i = k(\theta_i), \qquad
  K = \sum_{i=1}^p k_i, \qquad K_* = p\,k(\pi/p), \qquad z_i = \frac{w_i}{k_i}.
\]
Pairing $t$ and $-t$ gives $\int_{-\theta_i}^{\theta_i} u^3\,dt = w_i^3/k_i^2 + w_i(h_{i+1} - h_i)^2/4$. Because $\sum_i k_i z_i = 1$ and $K_*^{-2} = H_*^2 + x^2/6$, expanding about $1/K$ yields the nonnegative support-function inertia remainder
\begin{equation}\label{eq:support-remainder}
\begin{aligned}
  \mathcal{I} &\coloneqq 2I - \mathbb{E} h_i^2 - K_*^{-2} \\
  &= \frac{1}{4}\sum_{i=1}^p w_i(h_{i+1} - h_i)^2
    + \sum_{i=1}^p k_i\Bigl(z_i - \frac{1}{K}\Bigr)^2\Bigl(z_i + \frac{2}{K}\Bigr)
    + K^{-2} - K_*^{-2}.
\end{aligned}
\end{equation}
Notice that a straight angle ($\theta_i = 0$) contributes zero to both sums by continuity.

To see how $\mathcal{I}$ controls the exterior angles, observe that on $[0, \pi/2]$ we have
\[
  -k''(t) = \frac{4\sin^3 t}{3(1 - \frac{1}{3}\sin^2 t)^{5/2}} \ge \frac{32 t^3}{3\pi^3}.
\]
In the integral Taylor remainder of $k$ about $b = \pi/p$, using $b + s(t - b) \ge (1 - s)b$ for $s \in [0,1]$ and summing over $i$ (with $\sum_i \theta_i = \pi$) gives $K_* - K \ge \frac{32}{15 p^3}\sum_{i=1}^p (\theta_i - \pi/p)^2$. Since $k$ is concave with $k'(0) = 2$, we have $K \le K_* \le 2\pi$, and therefore
\begin{equation}\label{eq:angle-deficit}
  \mathcal{I} \ge K^{-2} - K_*^{-2} \ge \frac{2(K_* - K)}{K_*^3}
  \ge \frac{8}{15\pi^3 p^3}\sum_{i=1}^p \Bigl(\theta_i - \frac{\pi}{p}\Bigr)^2.
\end{equation}
Rewriting \eqref{eq:support-remainder} in terms of $V_h \coloneqq \mathbb{E}(h_i - H)^2$ gives
\begin{equation}\label{eq:pure-strip-inertia}
  I = \frac{H^2 + V_h + H_*^2 + \mathcal{I}}{2} + \frac{x^2}{12},
\end{equation}
or equivalently, relative to the regular apothem $H_*$,
\begin{equation}\label{eq:area-deviation}
  2H_*(H_* - H) + \mathcal{I} = \mathbb{E}\Bigl[(h_i - H_*)^2 + 2m_i^2 + \frac{y_i^2 - x^2}{6}\Bigr].
\end{equation}

We also record two elementary bounds on individual support heights and edge midpoints. For any single edge, rotate and translate so that it lies on $\{Y = 0\}$ with $Q \subset \{Y \ge 0\}$, and set $T = \max_Q Y$ and $h = \mathbb{E} Y$ (the support height of that edge from the centroid). Both boundary arcs joining a point of minimum height to a point of maximum height project onto $[0,T]$ with $ds \ge |dY|$, so $\mathbb{E}[Y(T - Y)] \ge 2\int_0^T u(T - u)\,du = T^3/3$. Maximizing $hT - T^3/3$ over $T > 0$ gives the upper bound below, while combining the two coordinate trials with \eqref{eq:pure-strip-inertia} yields the lower bounds: for every $0 < \sigma \le J(Q)$,
\begin{equation}\label{eq:edge-altitude}
  \frac{2}{3}h^{3/2} - h^2 \ge \Var(Y)
  \ge \frac{H_*^2 + x^2/6}{2} + \frac{H^2}{2} - \frac{H}{2\sigma}
  \ge \frac{H_*^2 + x^2/6}{2} - \frac{1}{8\sigma^2}.
\end{equation}
By the same projection argument onto the base of edge $i$, one has $y_i^3/6 \le \Var(X) \le H/(2\sigma)$. Next, parametrize $\partial Q$ at unit speed starting from the midpoint $P$ of an edge of length $y$, and choose centroid-centered coordinates $(X, Y)$ with the $X$-axis directed along $P$. In $X(0) = \int_0^1 (t - 1/2)X'(t)\,dt$, the two half-edges adjacent to $P$ cancel because $X'$ is constant there, so Cauchy--Schwarz on $[y/2, 1 - y/2]$ gives $12|P|^2 \le (1 - y)^3 \int_0^1 (X')^2\,dt$. Green's formula and the coordinate trial yield $|Q|^2 = (\int_0^1 X Y'\,dt)^2 \le (|Q|/\sigma)\int_0^1 (Y')^2\,dt$. Combining this with $(X')^2 + (Y')^2 = 1$ and $H = 2|Q|$ gives
\begin{equation}\label{eq:pure-strip-heights}
  24|P|^2 \le (1 - y)^3(2 - \sigma H).
\end{equation}

To sharpen the coordinate trial $X$ along each edge tangentially, introduce the odd harmonic half-strip profile
\[
  w(s, t) \coloneqq \sum_{k=1}^\infty \frac{(-1)^{k+1}}{k^2} e^{-2\pi k s}\sin(2\pi k t),
  \qquad s > 0, \quad |t| < \frac{1}{2},
\]
which vanishes on the lateral rays $|t| = 1/2$. All double integrals below extend over $(0,\infty) \times (-1/2, 1/2)$. Orthogonality gives
\begin{gather*}
  E_0 \coloneqq 2\iint w_s^2 = 2\iint w_t^2 = \pi Z, \\
  G_1 \coloneqq \int_{-1/2}^{1/2} t w(0,t)\,dt = \frac{Z}{2\pi}, \qquad
  G_2 \coloneqq \int_{-1/2}^{1/2} w(0,t)^2\,dt = \frac{W}{2}.
\end{gather*}
To evaluate the cross moment $2\iint t w_s w_t$, set $A = w_t^2 - w_s^2$ and $B = 2w_s w_t$. Harmonicity gives $B_s = -A_t$, and orthogonality gives $\int_{-1/2}^{1/2} A(s,t)\,dt = 0$ at each $s > 0$. Since $w_s = 0$ at $|t| = 1/2$ and $w_t$ is even in $t$, integrating by parts yields $\frac{d}{ds}\int t B\,dt = -w_t(s, 1/2)^2$. Integrating from $\infty$ and then over $s \in (0,\infty)$, with $w_t(s, 1/2) = 2\pi\log(1 - q)$ where $q = e^{-2\pi s}$, we obtain
\[
  2\iint t w_s w_t = \int_0^\infty s\,w_t(s, 1/2)^2\,ds
  = \int_0^1 \frac{-\log q}{q}\log^2(1 - q)\,dq
  = 3W - \zeta(2)^2 = G_2.
\]
Here the integral is $-2$ times the $a^2 b^2$-coefficient on the left-hand side of Euler's beta identity \eqref{eq:euler-beta-moments}; writing $R$ for its right-hand side, every term in $\log R$ is divisible by $ab$ with $[ab]\log R = -\zeta(2)$ and $[a^2 b^2]\log R = -3W/2$ by Euler's gamma product, and $\zeta(2)^2 = 5W/2$.

On the centroid sector $X = r X_i(t)$ ($0 < r \le 1$, $|t| \le 1/2$), introduce the logarithmic depth $u = -(h_i/y_i)\log r$ and define $W_i(X) = y_i \tau_i w(\gamma_i u, t)$ for a scaling parameter $\gamma_i > 0$, extending $W_i$ by zero to the other sectors. Because $W_i$ vanishes on the lateral rays and has finite energy, $W_i \in H^1(Q; \R^2)$. Denote by $E$ and $V$ the Dirichlet energy and arc-length boundary variance (together with their bilinear polarizations). Disjoint supports and zero boundary mean make the fields $W_i$ mutually orthogonal for both $E$ and $V$. Green's identity and tangentiality give $E(X, W_i) = 0$, while $V(X, W_i) = G_1 y_i^3$ and $V(W_i, W_i) = G_2 y_i^3$.

Setting $f(u,t) = w(\gamma_i u, t)$, the sector Jacobian is $|\det\partial_{(r,t)} X| = r h_i y_i$ and in the local $(\nu_i, \tau_i)$ frame
\[
  \nabla t = \Bigl(-\frac{m_i + y_i t}{r h_i y_i}, \frac{1}{r y_i}\Bigr), \qquad
  \nabla u = \Bigl(-\frac{1}{r y_i}, 0\Bigr),
\]
so that
\[
  E(W_i) = y_i^2 \iint \left[f_t^2 + \Bigl(f_u + \frac{m_i + y_i t}{h_i} f_t\Bigr)^2\right] du\,dt.
\]
The two terms odd in $t$ integrate to zero, while the cross moment contributes $G_2 y_i/h_i$ to $E(W_i)/y_i^2$. For the quadratic moment in $t$, the elementary bound $t^2 \le (1 - \cos(2\pi t))/8$ couples only adjacent cosine modes of $w_t$, giving
\[
  \frac{\iint t^2 w_t^2}{\iint w_t^2}
  \le \frac{1}{8}\left(1 + \frac{2}{Z}\sum_{j=1}^\infty \frac{1}{j(j+1)(2j+1)}\right) < 0.173,
\]
where the sum equals $3 - 4\log 2$ by partial fractions $1/j + 1/(j+1) - 4/(2j+1)$. Consequently, the terms involving $w_s^2$ and $w_t^2$ are bounded above by $\frac{E_0}{2}[\gamma_i + (1 + (m_i^2 + 0.173 y_i^2)/h_i^2)/\gamma_i]$. Optimizing over $\gamma_i > 0$ yields
\[
  \frac{E(W_i) - j V(W_i)}{y_i^2} \le K_i \coloneqq E_0\sqrt{1 + \frac{m_i^2 + 0.173 y_i^2}{h_i^2}} + G_2 y_i(h_i^{-1} - j).
\]
For any $0 < j \le J(Q)$, testing $X + \sum_{i=1}^p c_i W_i$ gives $0 \le H - j I + \sum_i (y_i^2 K_i c_i^2 - 2j G_1 y_i^3 c_i)$ for all real $c_i$, which forces every $K_i > 0$ and, upon minimizing over $c_i$, yields $H - j I \ge j^2 G_1^2 \sum_{i=1}^p y_i^4/K_i$.

Setting $s = H_* - H$, $b_i = y_i/x$, $k_i = K_i/E_0$, $\bar{k} = \mathbb{E} k_i$, $C_0 = (Z/\pi)(j/(2\pi))^2 x^3$, and $D = H_* - j(H_*^2 + x^2/12)$, and substituting \eqref{eq:pure-strip-inertia} together with $x\sum_i b_i = 1$, we arrive at the global strip inequality
\begin{equation}\label{eq:pure-strip}
  (1 - j H_*)s + \frac{j}{2}(s^2 + V_h + \mathcal{I}) + C_0 x\sum_{i=1}^p \Phi_{k_i}(b_i)
  \le D - C_0 + C_0(\bar{k} - 1),
\end{equation}
where $\Phi_k(b) \coloneqq b^4/k + bk - 5b + 3$. With $\psi(b) \coloneqq b^{5/2} - 1 - \frac{5}{2}(b - 1)$, AM--GM gives $\Phi_k(b) \ge 2\psi(b)$, and the factorization $\psi(b) - \frac{3}{2}(b - 1)^2 = b(\sqrt{b} - 1)^2(\sqrt{b} + 1/2)$ yields $\Phi_k(b) \ge 2\psi(b) \ge 3(b - 1)^2$.

\subsection{An equivariant trace inequality and sector transplantation}
Fix $p \ge 3$, normalize $R_p$ to unit perimeter with center at the origin, and orient it so that one outward edge normal points along the positive real axis. Write $x = 1/p$, $H_* = (x/2)\cot(\pi x)$, and $j = J(R_p)$. A vector field $Y \in H^1(R_p; \R^2)$ (identified with a complex-valued function) is \emph{equivariant} if $Y(\omega z) = \omega Y(z)$ and $Y(\bar{z}) = \overline{Y(z)}$, where $\omega = e^{2\pi i/p}$.

\begin{lem}[Equivariant trace inequality on the energy complement of $X$]\label{lem:sector-trace}
On a regular $p$-gon $R$ of apothem $h$, let $Y \in H^1(R; \R^2)$ be equivariant, nonzero, and energy-orthogonal to the coordinate field $X(z) = z$. Put $a = \tan(\pi/p)$. Then
\[
  \int_{\partial R} |Y|^2\,ds < h\,\frac{1 + a^2}{2}\,E(Y).
\]
Furthermore,
\[
  \int_{\partial R} |Y|^2\,ds \le h\,\frac{a}{\pi}\mathcal{B}_*(a)\,E(Y) \le h\,\frac{a}{\pi}\mathcal{B}(a)\,E(Y),
\]
where $\mathcal{B}(a) = 1 + \frac{3}{5}a(1 + a)$ and $\mathcal{B}_*(a)$ is defined in \eqref{eq:B-star} below.
\end{lem}

\begin{proof}
Rescale so that $h = 1$, denote the uniform average on the right edge $\{1 + it : |t| \le a\}$ by brackets $[\cdot]$, and write $d(U) = E(U)/|\partial R|$. Since $E(Y, X) = 0$ and $Y$ is equivariant, Green's formula and reflection symmetry give $[Y] = 0$. Consequently, the edge primitive $P(t) \coloneqq \int_{-a}^t Y(1 + is)\,ds$ vanishes at $t = \pm a$. Extend $P$ to the right sector by $U(r(1 + it)) = r P(t)$ and to $R$ by rotation. For any $b \in \C$, Stokes' theorem and $[Y] = 0$ yield
\[
  [|Y|^2] = |\partial R|^{-1}\int_R dY \wedge d\overline{(U - bX)}
  \le \sqrt{d(Y)\,d(U - bX)},
\]
extending from smooth fields to $H^1(R; \R^2)$ by density. Radial integration and an integration by parts on $[-a, a]$ give
\[
  2d(U) = \bigl[|P - tY|^2 + |Y|^2\bigr] = \bigl[(1 + t^2)|Y|^2\bigr] + 2\bigl[|P|^2\bigr].
\]
By Green's formula, the complex energy pairing of $U$ with $X$, divided by $|\partial R|$, equals $[P]$, while $d(X) = 1$; minimizing over $b \in \C$ therefore replaces $[|P|^2]$ by $\Var P$. Applying the Legendre inequality $2\Var P \le [(a^2 - t^2)|P'|^2]$ gives $2\min_b d(U - bX) \le (1 + a^2)[|Y|^2]$. (To verify the Legendre inequality, write $4a^2\Var P = \int_{-a \le t < s \le a} |P(s) - P(t)|^2\,dt\,ds$, apply $|P(s) - P(t)|^2 \le (s - t)\int_t^s |P'(u)|^2\,du$, and note that $\int_{-a}^u \int_u^a (s - t)\,ds\,dt = a(a^2 - u^2)$.) Equality forces $P'$ to be constant, whence $P \equiv 0$ by the endpoint conditions; since any nonzero equivariant field with vanishing boundary trace has strictly positive energy, the first inequality is strict.

For the sharper second bound, write $Y = f + ig$ on the right edge and put $c = [tg]$, so that $[P] = -ic$ by parity. Apply the same Stokes estimate with $U_\gamma(r(1 + it)) = r^\gamma K(t)$, where $K(t) = P(t) + ic(1 + it)$ and $\gamma = \pi/a > 1$. Radial integration gives $2\gamma\,d(U_\gamma) = [|\gamma K - tK'|^2 + |K'|^2]$. Integrating $t(|K|^2)'$ by parts and using $[t\Re P] = -[t^2 f]/2$, we obtain
\[
  2\gamma\,d(U_\gamma) - 2d(U_1)
  = (\gamma - 1)\left[(\gamma + 2)\bigl(\Var P + c[t^2 f]\bigr)
    + \Bigl(\frac{(\gamma - 1)a^2}{3} - 1\Bigr)c^2\right].
\]
On $[0, a]$, $\Re P$ vanishes at both endpoints and $\Im P + c$ is even with mean zero, so Wirtinger's inequality gives $\Var P \le (a/\pi)^2 [|Y|^2]$. Putting $u = [f^2]^{1/2}$, $v = \sqrt{3}\,|c|/a$, and $q = a^3/(3\sqrt{15})$, Cauchy--Schwarz gives $u^2 + v^2 \le [|Y|^2]$ and $|c[t^2 f]| \le 2quv$, because $[f] = 0$ and $[(t^2 - a^2/3)^2] = 4a^4/45$. Combining this with the negative $c^2$ term yields
\[
  q(2uv - h_0 v^2) \le q R(u^2 + v^2), \qquad
  h_0 = \frac{\sqrt{15}(1 - \pi a/3 + a^2/3)}{\pi + 2a}, \qquad
  R = \frac{\sqrt{h_0^2 + 4} - h_0}{2},
\]
since $R(R + h_0) = 1$. Substituting into the Stokes inequality gives the second bound with
\begin{equation}\label{eq:B-star}
  \mathcal{B}_*(a) \coloneqq 1 + \frac{a}{2\pi} + \frac{1 - 2/\pi^2}{2}a^2
  + R\,\frac{a(\pi - a)(\pi + 2a)}{6\sqrt{15}}.
\end{equation}
Because $\pi a - a^2 \le \pi^2/4 < 3$, we have $h_0 > 0$ and hence $R < 1$; replacing $R$ by $1$ yields coefficients of $a$ and $a^2$ strictly below $3/5$ and a negative cubic term, so $\mathcal{B}_*(a) \le \mathcal{B}(a)$.
\end{proof}

For the coordinate field $X$ on the unit-perimeter $R_p$, we have $E(X) = H_*$ and $V(X) = H_*^2 + x^2/12$. Expanding the cotangent and tangent series gives
\begin{equation}\label{eq:coordinate-ratio}
  \mathcal{L} \coloneqq \frac{2\pi V(X)}{E(X)}
  = \pi x\Bigl(\cot(\pi x) + \frac{1}{3}\tan(\pi x)\Bigr)
  = 1 + \sum_{k=2}^\infty \frac{2(4^k - 4)}{3}\zeta(2k)x^{2k} > 1 + 8x^4.
\end{equation}
On the energy complement of $X$, Lemma~\ref{lem:sector-trace} gives $2\pi V(Y)/E(Y) \le x\mathcal{B}_*(a)$, since $H_* = x/(2a)$. Splitting any field in the complex rotational class into real and imaginary parts under the reflection $Y(z) \mapsto \overline{Y(\bar{z})}$ preserves orthogonality to $X$ and adds both quadratic forms; codimension-one min--max therefore bounds the second variance eigenvalue $\lambda_2$ in the rotational class by $x\mathcal{B}_*(a)$. Since $(\pi - a)(\pi + 2a) \le 9\pi^2/8$, we have $\mathcal{B}_*(a) < 1 + 0.16a + 0.4a^2 + 0.48a R$; using $h_0 > 0.3$ (hence $R < 0.9$) for $a \le 3/4$, $h_0 > 0.5$ (hence $R < 0.8$) for $a \le 1/2$, and $a < 0.73, 0.484, 0.33$ at $p = 5, 7, 10$, we obtain
\begin{equation}\label{eq:sector-compression}
  \lambda_2 < 0.334 \quad (p \ge 5), \qquad
  \lambda_2 < 0.194 \quad (p \ge 7), \qquad
  \lambda_2 < 0.124 \quad (p \ge 10),
\end{equation}
since both $x$ and $a = \tan(\pi x)$ decrease with $p$. With $j = 2\pi/(1 + d_p)$, every equivariant $Y$ energy-orthogonal to $X$ satisfies
\begin{equation}\label{eq:coordinate-complement}
  E(Y) - j V(Y) \ge \eta E(Y), \qquad \eta \coloneqq 1 - \frac{x\mathcal{B}_*(a)}{1 + d_p},
\end{equation}
with $\eta > 0.64$ for $p \ge 5$, $\eta > 0.79$ for $p \ge 7$, and $\eta > 0.872$ for $p \ge 10$.

Taking the trace on the energy-orthonormal pair built from $X$ and $Y$ also yields
\[
  \frac{2\pi V(Y)}{E(Y)} \le 1 + d_p + \lambda_2 - \mathcal{L} < 1 + d_p,
\]
so the equivariant first eigenline pairs nontrivially with $X$. Let $U$ be the unique equivariant first eigenfield normalized by $E(U, X) = H_*$. By Green's formula, $\int_{\partial R_p} U \cdot \nu\,ds = H_*$. Writing $U = (u, w)$ in the $(\nu, \tau)$ frame of each centroid sector, the components $u$ and $w$ are respectively even and odd under reflection in the edge bisector, with edge means $H_*$ and $0$.

\paragraph{Sector transplantation and calibration.}
On a unit-perimeter convex $p$-gon $Q$ with boundary centroid at the origin, put $a_i = h_i/H_*$, $b_i = y_i/x$, and $v_i = m_i/H_*$. In the oriented $(\nu_i, \tau_i)$ frames of the reference and target sectors, define the affine sector map $L_i$ and the transplanted field $F$ by
\[
  L_i = \begin{pmatrix} a_i & 0 \\ v_i & b_i \end{pmatrix}, \qquad F(L_i Z) = L_i U(Z).
\]
Along every ray from the origin to a vertex of $R_p$, reflection symmetry forces $U$ to be parallel to that ray; since adjacent maps $L_i$ agree on the shared ray, the traces match and $F \in H^1(Q; \R^2)$. Because $U$ has mean $(H_*, 0)$ on each edge of $R_p$, the mean of $F$ on the $i$-th edge of $Q$ equals its midpoint $h_i\nu_i + m_i\tau_i$, so $F$ has zero global boundary mean on $\partial Q$.

Set $G = U - X$ and $D = E(X) - j V(X) = H_* - j(H_*^2 + x^2/12) > 0$ (equality would force $\nu = jX$ along each straight edge). Define
\[
  e = \frac{E(G)}{D}, \qquad r = \frac{j x^2}{12 D}, \qquad \tau = \frac{x^2}{12 H_*^2}.
\]
Because $E(G, X) = 0$ and the eigenvalue equation gives $E(G) - j V(G) = D$, inequality \eqref{eq:coordinate-complement} yields $e \le 1/\eta$. Let $X_0$ and $Y_0$ be the squared arc-length $L^2$-norms on $\partial R_p$ of the normal and tangential components of $U = (u, w)$. The normalization $E(U, X) = H_*$ gives $j\int_{\partial R_p} G \cdot X\,ds = D$, whose normal part vanishes because $G$ has zero normal mean on each edge. Expanding $U = X + G$ and using $E(U) = H_* + E(G)$ yields
\begin{equation}\label{eq:pure-calibration}
  \frac{j X_0}{D} = \frac{r}{\tau} + b_n, \qquad
  \frac{j Y_0}{D} = r + 2 + b_t, \qquad
  b_n + b_t = e - 1,
\end{equation}
where $b_n, b_t \ge 0$ are $j/D$ times the squared boundary norms of the normal and tangential components of $G$. Define the sector-integrated derivative moments
\[
  n = \frac{1}{D}\int_{R_p} u_t^2, \qquad
  k = \frac{1}{D}\int_{R_p} w_n^2, \qquad
  m = \frac{1}{D}\int_{R_p} \bigl[(u_n - w_t)^2 - 2u_t w_n\bigr].
\]
From the entries of $DG$, we have $n + k \le e$ and $m \le 2(e - n - k) + 2\sqrt{nk}$. Put $C = j(X_0 - Y_0) + 2D(k - n)$ and $T = j Y_0 - \tau C$. Equation~\eqref{eq:pure-calibration} gives
\[
  \frac{C}{D} \ge r(\tau^{-1} - 1) - 1 - 3e, \qquad
  \frac{T}{D} \ge 2 + \tau(r + 3 - e - 2k + 2n).
\]

On each edge, write $s = a - 1$, $z = b - 1$, and
\[
  A = \frac{k}{a} + n\frac{a + 2b}{b^2}, \qquad
  B = \frac{m}{a} + n\frac{2a + v^2/a}{b^2}.
\]
Changing variables on each sector yields
\begin{align*}
  V(F) &= X_0\,\mathbb{E}(a^2 + v^2) + Y_0\,\mathbb{E} b^2, \\
  E(F) &= \mathbb{E}\bigl\{a E(U) + 2D(k - n)(b - a) + D\bigl[A(s - z)^2 + B v^2\bigr]\bigr\}.
\end{align*}
To verify the energy formula, apply $|M|^2 = (\tr M)^2 - 2\det M + (M_{21} - M_{12})^2$ to $M = L_i (DU) L_i^{-1}$, whose skew difference is $[b^2 w_n - (a^2 + v^2)u_t + bv(u_n - w_t)]/(ab)$; multiplying by the Jacobian $ab$, using parity to drop terms odd in $t$, and substituting $b^2/a - a = 2(b - a) + (b - a)^2/a$ and $a^3/b^2 - a = 2(a - b) + (a + 2b)(a - b)^2/b^2$ gives the identity.

Put $g_0 = \mathbb{E} z = \mathbb{E}(z^2/b)$ and $\Delta = -\mathbb{E} s = 1 - H/H_*$. The linear terms in $j V(F) - E(F)$ equal $-C\Delta + 2(j Y_0 - Dk + Dn)g_0$. Rewriting \eqref{eq:area-deviation} as $\Delta = \mathbb{E}(s^2/2 + v^2) + \tau(2g_0 + \mathbb{E} z^2) - \mathcal{I}/(2H_*^2)$ and using $E(U) = H_* + eD$ with $e \ge n + k$, we obtain for $C \ge 0$:
\begin{equation}\label{eq:pure-sector}
  \frac{j V(F) - E(F)}{D}
  \ge \frac{C\mathcal{I}}{2D H_*^2}
    + \mathbb{E}\bigl\{U_1 s^2 + V_1 z^2 - A(s - z)^2 + (r + 2 - 2k + 2n - B)v^2\bigr\},
\end{equation}
where $U_1 = H_*/(2D) + (3n - k)/2$ and $V_1 = T/D + 2(T/D - k + n)/b$.

\begin{lem}[Sector stability box]\label{lem:sector-box}
Suppose that the calibration quantities of $R_p$ satisfy $C \ge 0$, $e \le E$, $n \le N$, $H_*/(2D) \ge S$, $r + 2 \ge R$, and $T/D \ge t \coloneqq 2 + \theta n$ for the constants in the appropriate row of the table below, and that every edge of $Q$ satisfies $a_i \in [a_-, a_+]$, $b_i \in [b_-, b_+]$, and $v_i^2 \le u(a_i)$:
\[
\begin{array}{@{}lcccccccc@{}}
  \toprule
  p & S & E & N & R & \theta & [a_-, a_+] & [b_-, b_+] & u(a) \\
  \midrule
  5, 6 & 14.5 & 1.57 & 1.57 & 6.2 & 3/5 & [0.84, 1.16] & [11/15, 19/15] & 1/16 \\
  7, 8, 9 & 49 & 1.266 & 1.266 & 9 & 0 & [0.84, 1.16] & [19/36, 53/36] & 1/16 \\
  \ge 10 & 173 & 1.15 & 0.185 & 14 & 0 & [0.69, 1.35] & [1/4, 7/4] & 1.82 - a^2 \\
  \bottomrule
\end{array}
\]
Then the transplanted regular eigenfield $F$ satisfies
\begin{equation}\label{eq:common-sector-margin}
  j V(F) - E(F)
  \ge \frac{C\mathcal{I}}{2H_*^2}
    + \frac{1}{4000}\,\mathbb{E}\bigl[H_*(a_i - 1)^2 + j x^2 v_i^2 + D(b_i - 1)^2\bigr].
\end{equation}
\end{lem}

\begin{proof}
Write $s = a - 1$, $z = b - 1$, and $V = t + 2(t - k + n)/b$. In \eqref{eq:pure-sector}, the form $U_1 s^2 + V z^2 - A(s - z)^2$ decreases with $k$, so we replace $k$ by its upper bound $E - n$. With $y = 1/b$, we have $A \le f \coloneqq k/a_- + n(2y + a_+ y^2)$. For fixed $(s, z)$, the expression $0.96(U_1 s^2 + V z^2) - f(s - z)^2$ is affine in $n$ and concave in $y$ separately, so it is nonnegative on $[0, N] \times [1/b_+, 1/b_-]$ provided it is nonnegative at all four corners. By weighted Cauchy--Schwarz, the corner condition is $f(1/U_1 + 1/V) \le 0.96$. Taking the larger value over $b \in \{b_-, b_+\}$ at $n = 0$ and $n = N$ yields upper bounds $(0.84, 0.96)$, $(0.54, 0.91)$, and $(0.58, 0.63)$ in the three rows. For instance, at the tightest corner ($n = N$ and $b = b_-$ in the first row), $U_1 > 16.85$, $f < 7.68$, and $V > 15.24$. Thus at least $0.04(U_1 s^2 + V z^2)$ remains in every row.

For the shear coefficient of $v^2$, put $M = \max(2, 2/a)$ and $X = (2a + u(a)/a)/b_-^2 - 2$. Using $m \le 2(E - n) + 2\sqrt{n(E - n)}$ and Young's inequality $2\sqrt{n(E - n)}/a \le \frac{7}{5}(E - n) + \frac{5}{7a^2}n$, the reduction from $r + 2$ is at most
\[
  M(E - n) + \frac{2\sqrt{n(E - n)}}{a} + Xn
  \le \Bigl(M + \frac{7}{5}\Bigr)(E - n) + \Bigl(X + \frac{5}{7a^2}\Bigr)n.
\]
In the first row, $M < 2.39$, $X < 2.46$, and $5/(7a^2) < 1.02$, bounding the loss by $3.79 E$. In the second row, $M + 7/5 < 3.79$, and the convex function $X + 5/(7a^2)$ takes endpoint values below $5.32$ and $7.06$, giving a loss below $7.06 E$. In the third row, $M < 2.90$, $X < 51.3$, and $5/(7a^2) < 1.51$, so the loss is at most $4.30(0.965) + 52.81(0.185) < 13.93$. In all three rows the loss stays below $0.996 R$, leaving a shear coefficient of at least $0.004(r + 2) \ge 0.004 r$.

Finally, $U_1 \ge H_*/(2D) - E/2$ and $V > 2$; multiplying by $D$ yields \eqref{eq:common-sector-margin}.
\end{proof}

\subsection{Convex polygons with at least ten sides}
Combining the strip budget \eqref{eq:pure-strip} with Lemma~\ref{lem:sector-box} settles all convex polygons with ten or more sides.

\begin{thm}\label{thm:convex-ten}
Every convex $p$-gon $Q$ with $p \ge 10$ satisfies $J(Q) \le J(R_p)$, with equality if and only if $Q$ is similar to $R_p$.
\end{thm}

\begin{proof}
Set $x = 1/p \le 1/10$, $d = d(R_p)$, $j = 2\pi/(1 + d)$, $H_* = (x/2)\cot(\pi x)$, $D = H_* - j(H_*^2 + x^2/12)$, $\gamma_0 = W/(2\pi Z) < 0.144$, and $C_\infty = Z/\pi \in (0.382, 0.383)$. Because $U(x)/x^3$ has nonnegative Taylor coefficients and Lemma~\ref{lem:regular-six-short} gives $U(1/10) < 0.003668$, we have $d \le U(x) < 3.668 x^3 < 3.684 x^3$. With $h = 2\pi H_*$, the cotangent expansion gives $h \le 1 - 2\zeta(2)x^2$ and $h^2 - h + \pi^2 x^2/3 \ge 8W x^4 + 12\zeta(6)x^6$. Inserting \eqref{eq:regular-chord} into $2\pi(1 + d)D = hd - (h^2 - h + \pi^2 x^2/3)$ gives
\[
  2\pi D < 2Z x^3 + x^5\bigl[27 - 4\zeta(2)Z + \bigl(130 - 16\zeta(2)W - 12\zeta(6)\bigr)x\bigr].
\]
The bracket is at most $19.1 + 90x \le 28.1$, whence $D < x^3(C_\infty + 4.48 x^2) < 0.428 x^3$.

\paragraph{Global a priori box from the strip trial.}
Assume $J(Q) \ge j$, and write $H = 2|Q|$, $s = H_* - H$, $z = jH$, $z_* = j H_*$, $\kappa = 1 - z_*$, $L = H/j - H^2$, $V_h = \mathbb{E}(h_i - H)^2$, $M_i = m_i^2 + 0.173 y_i^2$, and $C_0 = C_\infty (j/(2\pi))^2 x^3$. Since $h^2 + 2\pi^2 x^2/3 \ge 1$ by the cotangent series, \eqref{eq:pure-strip-inertia} and $I \le H/j$ imply $2H/j \ge H^2 + (2\pi)^{-2}$, so
\[
  z \ge 1 - \sqrt{1 - (j/(2\pi))^2} \ge 1 - \sqrt{2d} > 0.913, \qquad H > 0.14.
\]
Each support height satisfies $h_i \ge H/2 > 0.07$: in the inward coordinate $Y \ge 0$ measured from edge $i$, Green's formula gives $|Q| = \int_{\partial Q} Y\nu_Y\,ds \le \int_{\partial Q} Y\,ds = h_i$. We also have $L = H(1 - z)/j < 0.0024$.

The coordinate bound $V_h + \mathbb{E}(m_i^2 + y_i^2/12) \le L$ implies $\mathbb{E} M_i \le 2.076 L$ and $\mathbb{E} y_i^2 \le 12L$. Replacing every $h_i$ by $H$ in $k_i = K_i/E_0$ and using $\sqrt{1 + t} \le 1 + t/2$ bounds the resulting average strip cost above by
\[
  1 + (1 - z)\left[\frac{2.076}{2z} + \gamma_0\sqrt{\frac{12(1 - z)}{z}}\right] < 1 + 1.30(1 - z).
\]
Because $t \mapsto \sqrt{1 + t^2}$ is $1$-Lipschitz, restoring the individual heights $h_i \ge 0.07$ adds at most
\[
  \frac{(\sqrt{2.076} + \gamma_0\sqrt{12})\sqrt{L V_h}}{(0.14)(0.07)} < 9.8\sqrt{V_h}.
\]
From the cotangent series and $d < 3.684 x^3$, we have $3.28 x^2 < \kappa < 3.71 x^2$ and $H_* > 0.1538$. Substituting $1 - z = \kappa + j s$ and $9.8 C_0\sqrt{V_h} \le j V_h/4 + 9.8^2 C_0^2/j$ into \eqref{eq:pure-strip} yields the unified budget
\begin{equation}\label{eq:strip-boot}
  2.84 x^2 s + \frac{j}{4} V_h + C_0 x\sum_{i=1}^p \Phi_{k_i}(b_i)
  \le x^5(4.48 + 1.85 + 5.10x) < 6.85 x^5,
\end{equation}
using $C_0 > 0.3796 x^3$, $D - C_0 < 4.48 x^5 + 2.83 x^6$ (from $C_0 = C_\infty x^3/(1+d)^2$ and $C_\infty > 0.3825$), $\kappa - 1.30 j C_0 > 2.84 x^2$, $1.30 C_0 \kappa < 1.85 x^5$, and $9.8^2 C_0^2/j < 2.27 x^6$.

From \eqref{eq:strip-boot} we immediately read off $s < 2.64 x^3 \le 0.264 x^2$ and $V_h < 4.79 x^5 \le 0.479 x^4$. If any edge had $|b_i - 1| \ge 3/4$, then $\psi(b_i) \ge 29/32$, so that single edge would contribute more than $0.7592(29/32)x^4 > 6.85 x^5$ to \eqref{eq:strip-boot}, a contradiction. Hence $y_i > x/4$ for every $i$, and $|h_i - H| \le \sqrt{V_h/y_i} < 0.439 x$. For every edge midpoint $P_i$, applying \eqref{eq:pure-strip-heights} with $\sigma = j$ gives $24|P_i|^2 < 1$, since $2 - jH = 1 + \kappa + js < 1 + 5.38 x^2 < 1 + 3x/4 \le (1 - x/4)^{-3}$. Since $a_i^2 + v_i^2 = |P_i|^2/H_*^2$ and $H_* > 0.1538$, every edge lies in the box
\begin{equation}\label{eq:convex-box}
  \frac{1}{4} < b_i < \frac{7}{4}, \qquad 0.69 < a_i < 1.35, \qquad a_i^2 + v_i^2 < 1.82.
\end{equation}

\paragraph{Application of the exact eigenfield.}
By \eqref{eq:coordinate-complement}, we have $e < 1/0.872 < 1.15$, while $D < 0.428 x^3$ yields $r > 12$ and $H_*/(2D) > 173$. To bound the derivative moment $n = D^{-1}\int_{R_p} u_t^2$, let $U_s$ be the coordinate field $X$ on $R_p$ plus the optimal odd strip correction from Subsection~\ref{sec:convex-strip}. Since this correction is purely tangential, $E(U_s - U, X) = 0$ and the normal component of $U_s$ has zero tangential derivative in each sector. The regular strip cost satisfies $K_t = \sqrt{1 + 0.173 x^2/H_*^2} + \gamma_0 x\kappa/H_* < 1 + 5.64 x^2 < 1.057$. Applying the eigenvalue equation and \eqref{eq:coordinate-complement} to $U_s - U$ gives
\[
  \eta n D \le \eta E(U_s - U) \le E(U_s) - j V(U_s) \le D - \frac{C_0}{K_t} < (1 - 0.839)D,
\]
because $0.839(0.428)(1.057) < 0.3796$; thus $n < 0.161/0.872 < 0.185$.

The calibration bounds give $C > 0$ and $T/D \ge 2$, so \eqref{eq:convex-box} and the spectral constants satisfy the $p \ge 10$ row of Lemma~\ref{lem:sector-box}. Since $J(Q) \ge j$ implies $j V(F) - E(F) \le 0$, inequality \eqref{eq:common-sector-margin} forces $(a_i, b_i, v_i) = (1, 1, 0)$ for every edge. Hence $H = H_*$, and equality in the polygonal isoperimetric inequality forces $Q \sim R_p$.
\end{proof}

\subsection{Triangles}\label{sec:triangles}
Before turning to $4 \le p \le 9$, we treat triangles directly using Lemma~\ref{lem:sector-trace} and affine pullback from an equilateral triangle.

\begin{thm}\label{thm:pure-triangles}
Every triangle $P$ satisfies $J(P) \le J(R_3)$, with equality if and only if $P$ is equilateral. More precisely, if $P = L(T)$ for an equilateral triangle $T$ and an affine map $L$ whose linear part has singular-value ratio $r \in (0,1]$, then $J(R_3) - J(P) \ge c(1 - r)$ for an absolute constant $c > 0$.
\end{thm}

\begin{proof}
Take the equilateral reference triangle $T$ with vertices $-2, -2\omega, -2\omega^2$, where $\omega = e^{2\pi i/3}$. Its right edge lies on $\{x = 1\}$, its area is $|T| = 3\sqrt{3}$, and its perimeter is $|\partial T| = 6\sqrt{3} = 2|T|$. Write $D = E/|\partial T|$ for the normalized Dirichlet inner product, $B$ for the arc-length boundary inner product, and $[\cdot]$ for the uniform average on the right edge $\{1 + iy : |y| \le \sqrt{3}\}$. Set $Z(X) = X$, $K(Z) = -\bar{Z}^2/2$, and $\lambda = 1/\sigma$ with $\sigma = \sigma_1(T)$.

On the right edge, $K = (y^2 - 1)/2 + iy$ and $\partial_\nu K = Z - 2$. Since $[y^2] = 1$, Green's formula makes $Z$ and $K$ orthonormal for $D$. Testing $Z + sK$ for $s \in \R$ and using $[y^4] = 9/5$ gives
\[
  \lambda(1 + s^2) \ge B(Z + sK, Z + sK) = 1 + (1 + s)^2 + \frac{s^2}{5}.
\]
Dropping $s^2/5$ and taking the discriminant yields $(\lambda - 2)(\lambda - 1) \ge 1$, so $\lambda > 2$ and $(\lambda - 1)^2 \ge \lambda$. Let $V$ be an equivariant first eigenfield of $T$. By Lemma~\ref{lem:sector-trace}, $V$ pairs nontrivially with $Z$; normalize $V$ by $D(V, Z) = [V] = 1$. Writing $W = V - Z \neq 0$ (since $Z$ is not a first eigenfield), the eigenvalue equation and Lemma~\ref{lem:sector-trace} (with $(1 + \tan^2(\pi/3))/2 = 2$) yield
\[
  (\lambda - 2)\bigl(D(V, V) - 2\bigr) = B(W, W) - 2D(W, W) < 0,
\]
and therefore $D(V, V) < 2$.

Up to similarity, write $P = L(T)$ where the linear part of $L$ has singular values $1$ and $r \in (0,1]$ with major axis $e^{i\gamma}$. Put $Y = (1 - r)Z + r V$. From $B(V, Z) = \lambda$ and $[V] = 1$, we find $[y\Im V] = \lambda - 1$, whence Cauchy--Schwarz gives $2[(\Im V)^2] \ge 2(\lambda - 1)^2 \ge 2\lambda > \lambda D(V, V) = [|V|^2]$. By reflection symmetry, $[V^2] = [(\Re V)^2 - (\Im V)^2] < 0$, while $[ZV] = 2 - \lambda < 0$ and $[Z^2] = 0$; expanding the square therefore gives $[Y^2] \le 0$.

Let $s_1, s_2, s_3$ be the lengths of the images under $L$ of the unit tangents to the three edges of $T$, set $w_i = s_i/\sum_k s_k$, and put $\xi = \sum_{i=1}^3 w_i \nu_i$ for the outward unit normals $\nu_i$ of $T$. Because $c_i \coloneqq \cos(2\arg\nu_i - 2\gamma)$ satisfies $\sum_i c_i = 0$ and $s_i^2 = (1 + r^2 - (1 - r^2)c_i)/2$ decreases with $c_i$, we have $\sum_i w_i c_i \le 0$. Transplanting the scalar function $\Re(e^{-i\gamma} Y)$ via $L^{-1}$ therefore changes its boundary mean square from $B_r \coloneqq [|Y|^2]/2$ by $\frac{1}{2}[Y^2]\sum_i w_i c_i \ge 0$. Since $[Y] = 1$, its boundary mean on $\partial P$ is $\Re(e^{-i\gamma}\xi)$, so its variance on $\partial P$ is at least $B_r - |\xi|^2$. Subtracting the median of $\{s_1, s_2, s_3\}$ from each $s_i$ leaves $\sum_i s_i\nu_i$ unchanged and bounds its modulus by $\max s_i - \min s_i \le 1 - r$; combined with $\sum_i s_i \ge \sum_i s_i^2 = \frac{3}{2}(1 + r^2)$, this gives $|\xi| \le \frac{2}{3}(1 - r)$.

Set $x_\gamma = \Re(e^{-i\gamma} Z)$ and $g = \Re(e^{-i\gamma}(V - Z))$. Because $g$ has zero mean on every edge of $T$, we have $\int_T \nabla g = 0$, while rotational symmetry and $|\partial T| = 2|T|$ give $|T|^{-1}\int_T |\nabla g|^2 = D(V, V) - 1$. Pulling $x_\gamma + rg$ forward to $P$ by $L^{-1}$ yields a trial function whose Dirichlet energy on $P$ is at most $r|T| D(V, V)$ (the cross term integrates to zero and $\|r L^{-1}\| = 1$). On the other hand, $B(V, Z) = \lambda$ gives $2\sigma B_r = 2\sigma(1 - r)^2 + 2r(1 - r) + r^2 D(V, V)$. Subtracting the energy bound divided by $|T|$ from $2\sigma$ times the variance bound leaves at least
\[
  2\sigma\bigl[(1 - r)^2 - |\xi|^2\bigr] + r(1 - r)\bigl(2 - D(V, V)\bigr)
  \ge (1 - r)\min\Bigl\{\frac{10\sigma}{9},\; 2 - D(V, V)\Bigr\}.
\]
In particular, the trial function has strictly positive variance, and because its variance is uniformly bounded above and $2\sigma|T| = J(T)$, the variational principle yields $J(T) - J(P) \ge c(1 - r)$, with equality precisely when $r = 1$.
\end{proof}

\subsection{Convex polygons with four through nine sides}
For $4 \le p \le 9$, we construct explicit smooth equivariant fields on $R_p$ that fix the vertices, allowing a simple sector transplantation valid over the entire class of convex $p$-gons.

\paragraph{A vertex-fixing transplantation framework.}
Scale the regular $p$-gon $R_p$ to unit edge length and let $s = \frac{1}{2}\cot(\pi/p)$ be its apothem. Let $U$ be a smooth equivariant vector field on $R_p$ fixing all vertices, so that along every vertex ray $U(rV) = r c(r) V$ for a common real profile $c(r)$ with $c(1) = 1$; set $a(r) = (r c(r))'$. Writing $Z = s + it$ ($|t| \le 1/2$) on the right edge, denote the uniform $t$-average by $[\cdot]$ and the sector radial average by $\av{f} = \int_0^1 [f]\,r\,dr$. Define the moments
\begin{gather*}
  \alpha = \av{|U_y - i a|^2}, \qquad \beta = s^2 \av{|U_x - a|^2}, \\
  n = s - [\Re U(Z)], \qquad k_0 = [|U(Z) - s|^2], \qquad M = [|U(Z)|^2],
\end{gather*}
where the derivatives are evaluated at $rZ$. Since $E(U)$ is scale-invariant and $\Var(U) = M$, the perimeter-normalized Rayleigh quotient of $U$ is $J_0 = E(U)/M$.

On a unit-perimeter convex $p$-gon $Q$ with boundary centroid at the origin, put $x = 1/p$, $H_* = sx$, $w = h - H_*$, $d = y - x$, and $\Theta = 1 + \alpha - \beta/s^2$. On the centroid sector $X = r\nu[h + i(m + yt)]$, define the transplanted field
\[
  F(X) = \nu\bigl\{y U(rZ) + r c(r)(h - sy + im)\bigr\}.
\]
On every vertex ray, $F(X)$ reduces to $c(r)X$, so adjacent sectors glue continuously. In the local $(\nu, \tau)$ frame of the edge, the tangential and normal derivatives of $F$ are $U_y$ and $a + [sy(U_x - a) + m(ia - U_y)]/h$. Reflection parity eliminates the terms linear in $m$, while $[U_y] = i c(r)$ gives $\av{a^2 + |U_y|^2} - \alpha = 2\av{ac} = 1$. Integrating against the area element $r h y\,dr\,dt$ shows that $E(F)$ equals $(1 + \alpha)H + \mathbb{E}[(\alpha m^2 + \beta y^2)/h]$ plus a constant multiple of $\mathbb{E} y$. Rewriting $y^2/h = 2y/s - h/s^2 + (d - w/s)^2/h$ and calibrating the coefficient of $\mathbb{E} y$ on $Q = R_p$ (where $F(X) = x U(X/x)$ has energy $x^2 E(U) = J_0 x^2 M$ while $H = sx$ and $\mathbb{E} y = x$) determines that coefficient as $J_0 x M - s\Theta$. On the boundary, the edge mean of $F$ is $\nu(h + im - ny)$, and $\mathbb{E}\nu(h + im) = \mathbb{E}\nu = 0$. Therefore
\begin{align*}
  E(F) &= \Theta H + (J_0 x M - s\Theta)\mathbb{E} y
    + \mathbb{E}\!\left[\frac{\alpha m^2 + \beta(d - w/s)^2}{h}\right], \\
  V(F) &= \mathbb{E}\bigl[h^2 + m^2 - 2nyh + k_0 y^2\bigr] - n^2|\mathbb{E}(d\nu)|^2.
\end{align*}

When $U = P + \overline{Q}$ is harmonic with holomorphic $P, Q$ vanishing at $0$, write $P(rv) = v f(r)$ and $Q(rv) = \bar{v} g(r)$ at the vertex $v = s + i/2$. Equivariance forces $f$ and $g$ to be real with $f(1) + g(1) = 1$, and the sector moments simplify to
\begin{equation}\label{eq:harmonic-ray-moments}
\begin{aligned}
  \alpha &= \frac{E(U)}{2ps} - 1 + \int_0^1 r(f'^2 + g'^2)\,dr + 2\int_0^1 f' g'\,dr, \\
  \frac{\beta}{s^2} - \alpha &= 2g(1)^2 + 4\int_0^1 (1 - r)f'' g\,dr.
\end{aligned}
\end{equation}
To verify \eqref{eq:harmonic-ray-moments}, note that integrating along the edge segment gives $[\Re U_x] = (f - g)/r$ and $[\Im U_y] = (f + g)/r$, while taking a primitive of $P' Q'$ along the same segment gives $r[\Re(P' Q')] = 2\Im\int_0^r v P'(uv)Q'(uv)\,du = -\int_0^r f' g'\,du$. Thus $\av{|U_x|^2 - |U_y|^2} = -4\int_0^1 (1 - r)f' g'\,dr$, whereas $\av{|U_x|^2 + |U_y|^2} = E(U)/(ps)$; expanding the defining squares with $a = f' + g'$ and integrating by parts yields \eqref{eq:harmonic-ray-moments}.

Now put $c_* = 2J_0 x(s - n) - \Theta$. Suppose every edge satisfies $h_i \ge h_0 > 0$, and choose $B \ge 0$ such that $c_* \le 2H_* B$. Using $M = s^2 - 2sn + k_0$, the centered linear terms of $J_0 V(F) - E(F)$ are $-c_*(H_* - H)$ and $(J_0 x M - s c_*)\mathbb{E} d$. Applying $|\mathbb{E}(d\nu)|^2 \le \mathbb{E} d^2$, $\mathbb{E} d = \sum_i d_i^2 = \mathbb{E}(d^2/y)$, and \eqref{eq:area-deviation} produces
\begin{equation}\label{eq:quadform}
  J_0 V(F) - E(F) \ge \mathbb{E}\Bigl[A w^2 + 2L wd + \Bigl(C + \frac{D_0}{y}\Bigr)d^2 + G m^2\Bigr],
\end{equation}
where
\begin{align*}
  A &= J_0 - B - \frac{\beta}{s^2 h_0}, &
  L &= -J_0 n + \frac{\beta}{s h_0}, &
  G &= J_0 - 2B - \frac{\alpha}{h_0}, \\
  C &= J_0(k_0 - n^2) - \frac{B}{6} - \frac{\beta}{h_0}, &
  D_0 &= J_0 x M - s c_* - \frac{B x}{3}.
\end{align*}
Because $y_i \le 1/2$ on every convex polygon of unit perimeter, completing the square shows that whenever $A > t$ and $D_0 \ge Y \ge 0$,
\begin{equation}\label{eq:shared-small-form}
  (A - t)\bigl[C + 2(D_0 - Y)\bigr] \ge L^2
  \implies
  J_0 V(F) - E(F) \ge t\mathbb{E} w^2 + G\mathbb{E} m^2 + Y\sum_{i=1}^p d_i^2.
\end{equation}

\paragraph{Convex quadrilaterals.}
Applying \eqref{eq:quadform} directly to the exact first eigenfield of the square yields both rigidity and quadratic stability for convex quadrilaterals.

\begin{prop}[Convex quadrilaterals]\label{prop:convex-quad}
Every convex quadrilateral $Q$ satisfies $J(Q) \le J(R_4)$, with equality if and only if $Q$ is a square. Moreover, if $Q$ has unit perimeter and $J(Q) \ge 5.5$, the transplanted square eigenfield satisfies
\begin{equation}\label{eq:quad-stability-form}
  J(R_4) V(F) - E(F) \ge \frac{1}{20}\,\mathbb{E}\bigl(w^2 + d^2 + m^2\bigr).
\end{equation}
\end{prop}

\begin{proof}
Normalize $Q$ to unit perimeter, and on the unit-side reference square $[-1/2, 1/2]^2$ ($s = 1/2$) consider the field
\[
  U(z) = \kappa\bigl(\Re\sin(2\tau z) + i\Im\sinh(2\tau z)\bigr),
  \qquad \tan\tau\tanh\tau = 1, \quad \kappa = (2\sin\tau\cosh\tau)^{-1}.
\]
By Lemma~\ref{lem:reference-square}, both components of $U$ are first eigenfunctions of $R_4$. Since $U$ is equivariant and fixes all four vertices, \eqref{eq:quadform} applies with $J_0 = j = J(R_4)$.

Put $t = \tanh\tau = \cot\tau$, so $j = 8\tau t$ and $\kappa^2 = (1 - t^4)/4$. On the right edge $z = (1 + iv)/2$ ($|v| \le 1$),
\[
  U(z) = \frac{1}{2}\left(\frac{\cosh(\tau v)}{\cosh\tau} + i\frac{\sin(\tau v)}{\sin\tau}\right),
\]
which gives $n = (1 - t/\tau)/2$ and $4[|U|^2] = \frac{1}{2}(\operatorname{sech}^2\tau + \csc^2\tau) + (\tanh\tau - \cot\tau)/(2\tau) = 1$, hence $M = 1/4$ and $k_0 = n$. In \eqref{eq:harmonic-ray-moments}, the holomorphic parts are $2P = \kappa(\sin + \sinh)(2\tau z)$ and $2Q = \kappa(\sin - \sinh)(2\tau z)$, so $f'' = -2\tau^2 g$, $g'' = 2\tau^2 f$, and $f'g' = \tau^2(f^2 - g^2)$. The integrals in \eqref{eq:harmonic-ray-moments} admit the exact primitives
\[
  (f^2 + g^2)'' = 2(f'^2 + g'^2), \qquad
  (fg)'' = 4f' g', \qquad
  \left(\frac{f'^2 + g'^2}{2\tau^2} - fg\right)'' = 8\tau^2 g^2.
\]
Using $(f', g') = (2\kappa\tau, 0)$ at $r = 0$, together with $(f, g) = ((1 + t^2)/2, (1 - t^2)/2)$ and $f' = g' = \tau t$ at $r = 1$, we obtain
\[
  \alpha = \frac{j - 5 - t^4}{4}, \qquad 4\beta = \frac{j}{4} - 2t^2 - \frac{t^4}{2}, \qquad c_* = \kappa^2.
\]
By Lemma~\ref{lem:reference-square}, these satisfy $5.5 < j < 5.52$, $0.10 < n < 0.12$, $\kappa^2 < 0.18$, $\alpha < 0.06$, and $4\beta = 2t(\tau - t) - t^4/2 < 0.16$.

Assume $J(Q) \ge 5.5$. Inequality \eqref{eq:edge-altitude} gives $(2/3)h_i^{3/2} - h_i^2 \ge 5/384 - 1/242 > 0.00888$, whence $h_i > 0.08$ on every edge. Taking $h_0 = 0.08$ and $B = 3/4$ in \eqref{eq:quadform} (valid since $c_* = \kappa^2 < 2H_* B = 3/16$), we find $\beta/h_0 < 1/2$, $j n(1 - n) > 0.495$, and $D_0 = j/16 - \kappa^2/2 - B/12 > 0.191$, yielding
\[
  A > \frac{11}{4}, \qquad |L| < \frac{2}{3}, \qquad C + 2D_0 > \frac{1}{4}, \qquad G > 3.
\]
Applying \eqref{eq:shared-small-form} with $t = 1/20$ and $Y = 1/40$ (since $(11/4 - 1/20)(1/4 - 1/20) > 4/9$) and using $\sum_i d_i^2 \ge 2\mathbb{E} d_i^2$ establishes \eqref{eq:quad-stability-form}. In particular, $J(Q) \ge J(R_4)$ forces $w = d = m = 0$ on every edge, so $H = H_*$ and $Q$ is a square.
\end{proof}

\paragraph{Convex polygons with five through nine sides.}
For $5 \le m \le 9$, we first build an explicit three-mode harmonic field on $R_m$ whose Rayleigh quotient approximates $J(R_m)$ accurately enough to confine any competitor to the stability box of Lemma~\ref{lem:sector-box}.

\begin{lem}[Three-mode harmonic reference fields]\label{lem:three-mode-reference}
Write $j_m = J(R_m)$. For each $5 \le m \le 9$, there exists a smooth equivariant harmonic field $W$ on $R_m$ fixing the vertices whose normalized Rayleigh quotient $\widehat{j}$ satisfies
\[
  0 \le \widehat{j} - j_m <
  \begin{cases}
    0.0062 & (m = 5), \\
    0.004 & (m = 6, 7), \\
    0.0027 & (m = 8), \\
    0.0023 & (m = 9),
  \end{cases}
  \qquad j_5 > 5.9725.
\]
\end{lem}

\begin{proof}
Normalize $R_m$ to unit edge length and put $u = \sin(\pi/m)$, $c = \cos(\pi/m)$, $R = 1/(2u)$, $s = Rc$, $h = m - 1$, and $l = m + 1$. With $z = X/R$ and $g = 1 + (P - 1)t$, define
\[
  W(X) = \frac{R}{g}\bigl(z - P t \bar{z}^h + t z^l\bigr),
\]
where the parameters $(P, t)$ are chosen as $(4, 3/80)$ for $m = 5$, $(14/5, 1/32)$ for $m = 6$, $(3, 0.019)$ for $m = 7$, $(5/2, 0.016)$ for $m = 8$, and $(2, 0.014)$ for $m = 9$.
At every vertex $V$, each nonlinear monomial divided by $z$ equals $-1$, so the radial trace is shared across all vertex rays and $W(V) = V$. On the edge $z = c + iuv$ ($|v| \le 1$), denoting the uniform $v$-average by $[\cdot]$, integrating $z^k$ and differentiating $v(c^2 + u^2 v^2)^k$ give
\[
  [z^k] = \frac{\sin((k + 1)\pi/m)}{(k + 1)u}, \qquad
  L_k \coloneqq [|z|^{2k}] = \frac{1 + 2k c^2 L_{k-1}}{2k + 1}, \qquad L_0 = 1.
\]

Setting $f = g W/R$, $d_0 = g\partial_\nu W$, $S = [|f|^2]$, and $cT = [\Re(d_0 \bar{f})]$, the edge averages and Green's identity give (using $\bar{z} = 2c - z$ to evaluate $[|z|^2 z^m] = -1/(m+3) - 4c^2/[(m+2)(m+3)]$)
\begin{align*}
  S &= L_1 + t\left[\frac{2P}{m + 1} - \frac{2}{m + 3} - \frac{8c^2}{(m + 2)(m + 3)}\right]
    + t^2\Bigl(P^2 L_h + L_l - \frac{2P}{2m + 1}\Bigr), \\
  T &= 1 - \frac{4t}{m + 2} + t^2\bigl(P^2 h L_{m-2} + l L_m\bigr), \\
  M &= \frac{R^2 S}{g^2}, \qquad \widehat{j} = \frac{2muc T}{S}, \qquad
  \mu = s\left[1 - \frac{1 - 2t/(m + 2)}{g}\right], \qquad v_* = M - (s - \mu)^2.
\end{align*}
In \eqref{eq:harmonic-ray-moments}, the ray traces are $f(r) = (r - t r^l)/g$ and $g(r) = P t r^h/g$, with $E(W)/(2ms) = T/(2g^2)$. Direct polynomial integration yields
\begin{align*}
  2g^2 \alpha &= T - 1 + \frac{4t}{m + 2}
    + t^2\left[P^2(m - 3) + m - 1 - \frac{4Pm(m - 2)}{2m - 1}\right], \\
  g^2\Bigl(\frac{\beta}{s^2} - \alpha\Bigr) &= 2t^2 P\Bigl(P - \frac{l}{2m - 1}\Bigr).
\end{align*}

Aside from $c$, only $L_{m-2}$ is needed, as $L_{m-1}, L_m, L_{m+1}$ follow from the recurrence. With $\varphi = (1 + \sqrt{5})/2$, we have $2c_5 = \varphi$, $2c_6 = \sqrt{3}$, and $c_8^2 = (2 + \sqrt{2})/4$, while alternating Taylor bounds at $3.1415926/m < \pi/m < 3.1415927/m$ give $0.81174 < c_7^2 < 0.81175$ and $0.883021 < c_9^2 < 0.883024$. Because the recurrence increases in both $c^2$ and $L_{k-1}$, endpoint evaluation gives
\[
  L_3 = \frac{7 + 4\varphi}{28}, \quad L_4 = \frac{319}{630}, \quad
  L_5 \in (0.532, 0.536), \quad L_6 \in (0.558, 0.562), \quad L_7 \in (0.582, 0.586),
\]
for $m = 5, 6, 7, 8, 9$, respectively, and hence
\[
  \widehat{j} < 5.9787, \quad 6.138, \quad 6.204, \quad 6.235, \quad 6.252,
\]
respectively. For example, at $m = 5$ the bounds $1.6180339 < \varphi < 1.6180340$ and $1.175570 < R^{-1} < 1.175571$ yield $80^2 S > 5215.67$ and $80^2 cT < 5305.11$; at $m = 8$, $S > 0.90820$ and $T < 1.00099$ satisfy $5.65686(1.00099) < 6.235(0.90820)$. On the other hand, Lemma~\ref{lem:regular-six-short} gives $j_m = 2\pi/(1 + d_m) > 5.9725, 6.134, 6.20, 6.2323, 6.2497$ for $m = 5, \dots, 9$. Subtracting yields all five stated bounds.
\end{proof}

\begin{thm}\label{thm:convex-five-nine}
For $5 \le m \le 9$, every convex $m$-gon $Q$ satisfies $J(Q) \le J(R_m)$, with equality if and only if $Q$ is similar to $R_m$.
\end{thm}

\begin{proof}
Set $x = 1/m$, $j = J(R_m)$, and $H_* = sx$. On a unit-perimeter convex $m$-gon $Q$, let $0 < \sigma \le \min\{J(Q), j\}$ be a test parameter (for the global comparison we take $\sigma = j$). Apply \eqref{eq:quadform} to the harmonic field $W$ of Lemma~\ref{lem:three-mode-reference} with $J_0 = \widehat{j}$, $n = \mu$, $k_0 = k_* \coloneqq v_* + \mu^2$, and $B = c_*/(2H_*) > 0$, noting that $D_0 = s\Theta + \widehat{j}x(k_* - s^2) - Bx/3$. Using $I \le H/\sigma$, $\mu > 0$, and $y_i \le 1/2$, the boundary variance of $F$ satisfies
\[
  V(F) \le \frac{H_*}{\sigma} + \frac{1}{2}\Bigl[\frac{k_* - 1/12}{2} - 2\mu h_0\Bigr]_+ \eqqcolon \overline{V}(\sigma).
\]
Whenever \eqref{eq:shared-small-form} holds with $t = 2$, the variational inequality $\sigma V(F) \le E(F)$ yields
\begin{equation}\label{eq:small-mode-budget}
  2\mathbb{E} w^2 + G\mathbb{E} m_i^2 + Y\sum_{i=1}^m d_i^2
  \le (\widehat{j} - \sigma)\overline{V}(\sigma) \eqqcolon \varepsilon_m(\sigma).
\end{equation}

Assume $J(Q) \ge j$ and set $\sigma = j$. By the cotangent series, $H_*^2 + x^2/6 \ge (2\pi)^{-2}$. Since $j > 5.9$ and $(2/3)(0.08)^{3/2} - (0.08)^2 < 1/(8\pi^2) - 1/[8(5.9)^2]$, inequality \eqref{eq:edge-altitude} implies $h_i > 0.08$ on every edge. Thus $h_0 = 0.078$ is admissible for all $5 \le m \le 9$, and continues to hold for all $\sigma$ in a left neighborhood of $j$.

Choose $Y = 1/20$ for $m = 5, 6$ and $Y = 1/50$ for $m = 7, 8, 9$. Direct substitution of the moment formulas from Lemma~\ref{lem:three-mode-reference} gives $B > 0$, $D_0 > Y$, $A > 3.30$, and
\[
  |L| < 0.114, \qquad C + 2(D_0 - Y) > 0.011, \qquad
  G > (1.9,\; 1.2,\; 0.87,\; 0.60,\; 0.42) \quad (m = 5, \dots, 9).
\]
Because $1.30(0.011) > 0.114^2$, condition \eqref{eq:shared-small-form} holds for all five values of $m$. (For instance, at $m = 5$ and $m = 6$ one has $0.0256 < \beta < 0.02575$ and $0.0165 < \beta < 0.0168$, giving $0.020 < L < 0.028$, $-0.082 < L < -0.068$, and $C + 2(D_0 - Y) > 0.08, 0.04$, respectively.) The same moments give $\overline{V}(j) < 0.0265$ in all five cases.

Write $\varepsilon = \varepsilon_m(j)$ and $q_m^2 = m(m - 1)\varepsilon/Y$ (the case $\varepsilon = 0$ is immediate). Since $\sum_k d_k = 0$, Cauchy--Schwarz yields $\sum_k d_k^2 \ge m d_i^2/(m - 1)$, so \eqref{eq:small-mode-budget} gives
\[
  \frac{(b_i - 1)^2}{q_m^2} + \frac{x H_*^2 b_i}{\varepsilon}\bigl[2(a_i - 1)^2 + G v_i^2\bigr] \le 1.
\]
Hence $|b_i - 1| \le q_m$, with $q_m < 4/15$ for $m = 5, 6$ and $q_m < 17/36$ for $m = 7, 8, 9$. Completing the square gives $1 - (b - 1)^2/q_m^2 \le 2b/(1 + \sqrt{1 - q_m^2})$, whence
\[
  2(a_i - 1)^2 + G v_i^2 \le \frac{2m\varepsilon}{H_*^2(1 + \sqrt{1 - q_m^2})} < \frac{16m\varepsilon}{15 H_*^2}.
\]
Using $2\pi H_* > 1 - \frac{7}{2}x^2$ from the cotangent series, the preceding bounds yield $m\varepsilon/H_*^2 < 0.047$ and $m\varepsilon/(G H_*^2) < 0.057$, so $|a_i - 1| < 0.16$ and $|v_i| < 1/4$ for every $5 \le m \le 9$.

To verify the spectral hypotheses of Lemma~\ref{lem:sector-box}, retain the first two positive terms of \eqref{eq:coordinate-ratio}: set $\gamma = 8.64 x^4 + 40 x^6$ and $\rho_m = (d_m - \gamma)/x^2$. Because $\mathcal{L} > 1 + \gamma$ and $D > 0$, we have $\rho_m > 0$. Using $D/H_* = 1 - \mathcal{L}/(1 + d_m)$ and $2\pi H_* < 1 - 3.28 x^2$, we obtain
\[
  \frac{H_*}{2D} > \frac{1}{2} + \frac{1 + \gamma}{2\rho_m x^2}, \qquad
  r > \frac{3.28}{(1 - 3.28 x^2)\rho_m}.
\]
Lemma~\ref{lem:regular-six-short} shows that $\rho_m = m^2 d_m - 8.64/m^2 - 40/m^4$ is strictly below $0.886, 0.601, 0.46, 0.46, 0.46$ for $m = 5, \dots, 9$. Thus $H_*/(2D) > 14.5$ and $r > 4.2$ for all $m \ge 5$, while $H_*/(2D) > 49$ and $r > 7$ for $m \ge 7$. By \eqref{eq:coordinate-complement}, $e < 1/0.64$ for $m \ge 5$ and $e < 1/0.79$ for $m \ge 7$. Since $\tau < 1/5$, calibration yields $C > 0$; for $m = 5, 6$, $\tau \ge 1/9$ and $r + 3 - 3e > 7n/5$ give $T/D \ge 2 + 3n/5$, while for $m \ge 7$ we have $T/D \ge 2$.

Thus $Q$ lies strictly inside the corresponding box of Lemma~\ref{lem:sector-box}. Since $J(Q) \ge j$ gives $j V(F) - E(F) \le 0$, inequality \eqref{eq:common-sector-margin} forces all edge deviations to vanish, and the polygonal isoperimetric equality yields $Q \sim R_m$.
\end{proof}

\subsection{Quantitative stability and affine images}
Having established convex rigidity for every $p \ge 3$, we now prove Theorem~\ref{thm:convex} and its corollaries. We begin with the sharp behavior for triangles and quadrilaterals.

\begin{cor}[Stability for triangles and convex quadrilaterals]\label{cor:small-stability}
Let $Q$ be a unit-perimeter convex $p$-gon with $p \in \{3, 4\}$, and put $\mathcal{D} = \mathbb{E}[(a_i - 1)^2 + (b_i - 1)^2 + v_i^2]$. There exist absolute constants $c, C > 0$ such that
\[
  c\sqrt{\mathcal{D}} \le J(R_3) - J(Q) \le C\sqrt{\mathcal{D}} \quad (p = 3), \qquad
  J(R_4) - J(Q) \ge c\mathcal{D} \quad (p = 4).
\]
For triangles, one has $\mathcal{D} \asymp 1 - |Q|/|R_3| \asymp (1 - r)^2$, where $R_3$ has the same perimeter as $Q$ and $r$ is the singular-value ratio of an affine map sending $R_3$ onto $Q$.
\end{cor}

\begin{proof}
For triangles ($p = 3$), Theorem~\ref{thm:pure-triangles} already supplies the lower bound $J(R_3) - J(Q) \ge c(1 - r)$. Conversely, affine pullback by $L^{-1}$ multiplies Dirichlet energy by at least $r$ and arc-length boundary variance by at most $1/\bar{s}$ (where $\bar{s} = (s_1 + s_2 + s_3)/3$, since $s_i \le 1$ and variance minimizes over additive constants). Hence $J(Q) \ge r\bar{s}\,J(R_3) \ge r^2 J(R_3)$, giving $J(R_3) - J(Q) \le 2J(R_3)(1 - r)$.

Write $\rho = |Q|/|R_3|$ at equal perimeter and $W_\ell = \mathbb{E}(b_i - 1)^2$. The edge midpoints $M_1, M_2, M_3$ have weighted average equal to the boundary centroid and pairwise distances $|M_i - M_j| = y_k/2$ ($k \neq i, j$), so their weighted variance equals $\sum_{i < j} y_i y_j |M_i - M_j|^2 = y_1 y_2 y_3/4$. Dividing by $H_*^2 = 1/108$ gives $\mathbb{E}(a_i^2 + v_i^2) = 27 y_1 y_2 y_3$; since $\mathbb{E} a_i = \rho$, it follows that $\mathcal{D} = 1 + 27 y_1 y_2 y_3 - 2\rho + W_\ell$. In terms of $e_i = 3y_i - 1$ (which satisfy $\sum_i e_i = 0$), we have $W_\ell = \frac{1}{3}\sum_i e_i^2 + e_1 e_2 e_3$ and Heron's formula $\rho^2 = \prod_{i=1}^3 (1 - 2e_i)$. Direct algebraic expansion yields the identity
\[
  \mathcal{D} + 8W_\ell = \frac{1}{4}(1 - \rho)(13 + 5\rho).
\]
Since $0 \le W_\ell \le \mathcal{D}$ and $0 < \rho \le 1$, this gives $\frac{13}{36}(1 - \rho) \le \mathcal{D} \le \frac{9}{2}(1 - \rho)$. Moreover, $\rho = r/\bar{s}^2$; taking $L$ symmetric positive definite and averaging $\tau_i \cdot L\tau_i \le s_i$ and $s_i^2 = \tau_i \cdot L^2\tau_i$ via $\sum_i \tau_i \otimes \tau_i = \frac{3}{2}I$ yields $(1 + r)/2 \le \bar{s} \le \sqrt{(1 + r^2)/2}$, and thus
\[
  \frac{(1 - r)^2}{(1 + r)^2} \le 1 - \rho \le \frac{(1 - r)^2}{1 + r^2}.
\]
Consequently $J(R_3) - J(Q) \asymp 1 - r \asymp \sqrt{1 - \rho} \asymp \sqrt{\mathcal{D}}$.

For convex quadrilaterals ($p = 4$), $\mathcal{D}$ is uniformly bounded above because $\operatorname{diam}(Q) \le 1/2$, $H_* = 1/8$, and $b_i \le 2$. If $J(Q) < 5.5$, the gap $J(R_4) - J(Q)$ is bounded below by $J(R_4) - 5.5 > 0$. Otherwise, Proposition~\ref{prop:convex-quad} yields \eqref{eq:quad-stability-form}: $J(R_4)V(F) - E(F) \ge \frac{1}{20}\mathbb{E}(w^2 + d^2 + m^2)$. Since $a_i - 1 = 8w_i$, $b_i - 1 = 4d_i$, $v_i = 8m_i$, and $V(F)$ is uniformly bounded above (as $|h_i + im_i| \le 1/2$, $y_i \le 1/2$, and $U$ is fixed), the variational inequality $J(Q)V(F) \le E(F)$ proves $J(R_4) - J(Q) \ge c\mathcal{D}$.
\end{proof}

\begin{proof}[Proof of Theorem~\ref{thm:convex}]
Put $x = 1/p$, $j = J(R_p)$, $\sigma = J(Q)$, and $g = j - \sigma$. By Theorems~\ref{thm:convex-ten}, \ref{thm:pure-triangles}, \ref{thm:convex-five-nine}, and Proposition~\ref{prop:convex-quad}, we have $g \ge 0$ for every convex $p$-gon. Since $\mathcal{I} \le 2I \le 1/2$, $h_i^2 + m_i^2 \le 1/4$, $H_* \ge 1/(6\sqrt{3})$, and $b_i \le 1/(2x)$, the geometric quantity
\[
  \mathcal{G}(Q) \coloneqq \mathcal{I} + \mathbb{E}\bigl[(a_i - 1)^2 + x^2 v_i^2 + x^3(b_i - 1)^2\bigr]
\]
is uniformly bounded above across all $p \ge 3$ and all unit-perimeter convex $p$-gons $Q$. It therefore suffices to assume $g$ is sufficiently small.

Whenever $Q$ lies in the corresponding box of Lemma~\ref{lem:sector-box}, the transplanted regular eigenfield satisfies $V(F) \le 4(X_0 + Y_0) = 4(H_* + eD)/j$, while calibration gives $C/H_* \ge j H_*(1 - \tau) - (1 + 3e)D/H_* > 1/3$ (using $j > 5.9$, $H_* \ge 1/8$, $\tau < 1/5$, $e \le 1.57$, and $D/H_* \le 1/29$). Since $D \asymp x^3$ by \eqref{eq:regular-two} and $D > 0$ at every fixed $p$, inequality \eqref{eq:common-sector-margin} combined with $j V(F) - E(F) \le g V(F)$ yields $\mathcal{G}(Q) \lesssim g$.

For $p \in \{3, 4\}$, the bound follows directly from Corollary~\ref{cor:small-stability} and \eqref{eq:area-deviation}: with $d_i = y_i - x$, one has $\mathbb{E}(y_i^2 - x^2) = \mathbb{E} d_i^2 + 2x\sum_i d_i^2$ and $\sum_i d_i^2 \le 2(\sum_{d_i > 0} d_i)^2 \le 2p^2 \mathbb{E} d_i^2$. For each fixed $p \ge 5$, the a priori estimates in Theorems~\ref{thm:convex-five-nine} and~\ref{thm:convex-ten} place every edge strictly inside the box of Lemma~\ref{lem:sector-box}, and all intermediate scalar bounds depend continuously on the parameter $\sigma$ at $\sigma = j$ with strictly positive denominators; hence for each fixed $p \ge 5$, there is $\delta_p > 0$ such that $g < \delta_p$ places $Q$ inside that box. This handles any finite range of $p$, so it remains to consider $x = 1/p$ and $g$ both sufficiently small.

Write $s = H_* - H$, $V_h = \mathbb{E}(h_i - H)^2$, and $L = H/\sigma - H^2$, so the coordinate bound $I \le H/\sigma$ reads $V_h + \mathbb{E}(m_i^2 + y_i^2/12) \le L$. For $x$ and $g$ small, $\sigma > 5.9$, so \eqref{eq:edge-altitude} forces $h_i > 0.08$ on every edge. Thus the normalized strip costs $k_i = K_i/E_0$ evaluated at $\sigma$ satisfy $0 < k_i \le K$ and $\bar{k} = \mathbb{E} k_i \le 1 + C L$ for absolute constants $K, C > 0$: the radical in $k_i$ is at most $1 + C(m_i^2 + y_i^2)$, while $\mathbb{E}[y_i(h_i^{-1} - \sigma)] \le C\mathbb{E}[y_i|h_i - H|] + C(1 - \sigma H) \le C' L$, and $L = O(1 - \sigma H) = O(x^2 + s + g)$.

Apply \eqref{eq:pure-strip} with $j$ replaced by $\sigma$, writing $D_\sigma = H_* - \sigma(H_*^2 + x^2/12)$. On $0 < k \le K$, the function $\Phi_k(b)$ satisfies $\Phi_k(b) \ge c_K(1 + b)(b - 1)^2$ (since $2\psi(b)$ is strictly positive away from its double zero at $b = 1$, while $b^4/K$ dominates as $b \to \infty$). Setting $D_y = x\sum_i (b_i - 1)^2$, the length term therefore controls $C_0(D_y + \mathbb{E}(b_i - 1)^2)$ up to an absolute factor. Since the regular expansion \eqref{eq:regular-two} gives $D_\sigma - C_0 \le C(g + x^5)$ and $1 - \sigma H_* \ge c x^2$, the error $C_0(\bar{k} - 1) \le C x^3(x^2 + s + g)$ is absorbed into the left-hand side of \eqref{eq:pure-strip} for small $x$, yielding
\[
  s^2 + V_h + \mathcal{I} + x^2 s + x^3\bigl(D_y + \mathbb{E}(b_i - 1)^2\bigr) \le C(g + x^5).
\]
Notice that $s^2 + V_h = H_*^2 \mathbb{E}(a_i - 1)^2$, and since $\mathbb{E} y_i^2 \ge x^2$, \eqref{eq:area-deviation} gives $2\mathbb{E} m_i^2 \le \mathcal{I} + 2H_* s$. Consequently, the augmented quantity $\mathcal{A} \coloneqq \mathcal{G}(Q) + x^3 D_y$ satisfies $\mathcal{A} \le C(g + x^5)$.

If $Q$ lies outside the box \eqref{eq:convex-box}, then $\mathcal{A} \ge c_0 x^4$: either some edge has $|b_i - 1| \ge 3/4$, giving $D_y \ge \frac{9}{16}x$, or else every $y_i \ge x/4$ while some edge has $|a_i - 1| \ge 1/10$ or $|v_i| \ge 1/10$, contributing at least $c_0 x^4$ to $\mathbb{E}[(a_i - 1)^2 + x^2 v_i^2]$. For small $x$, the error $C x^5$ is then absorbed by $\mathcal{A}$, giving $\mathcal{G}(Q) \le \mathcal{A} \lesssim g$. Inside \eqref{eq:convex-box}, Lemma~\ref{lem:sector-box} already gives $\mathcal{G}(Q) \lesssim g$. Finally, the exterior-angle estimate follows immediately from \eqref{eq:angle-deficit}, and the triangle characterization is Corollary~\ref{cor:small-stability}.
\end{proof}

\begin{rem}[Optimality of the side-length exponent]\label{rem:length-weight-sharp}
The weight $p^{-3}$ in front of $\mathbb{E}(b_i - 1)^2$ in Theorem~\ref{thm:convex} cannot be improved. Starting from the unit-perimeter regular $(p-1)$-gon $R_{p-1}$, bisect one edge and push its midpoint outward along the edge normal by a small distance $t > 0$, then rescale back to unit perimeter. For each fixed $p \ge 4$ and sufficiently small $t > 0$, this produces a strictly convex $p$-gon $Q_{p,t}$. As $t \downarrow 0$, the piecewise affine maps on the radial fan from the center through the vertices and the added midpoint have derivatives converging uniformly to the identity, so both Dirichlet energies and centered boundary variances converge multiplicatively for all trial functions, giving $J(Q_{p,t}) \to J(R_{p-1})$. In this limit, the two half-edges have combined weight $1/(p-1)$ with $b_i = p/[2(p-1)]$, while the remaining $p - 2$ edges have $b_i = p/(p-1)$, so $\mathbb{E}(b_i - 1)^2 \sim 1/(4p)$. On the other hand, \eqref{eq:regular-two} implies $J(R_p) - J(R_{p-1}) \sim 12\pi Z p^{-4}$. Choosing $t = t_p$ sufficiently small yields a sequence of convex $p$-gons with $J(R_p) - J(Q_{p,t_p}) \asymp p^{-3}\mathbb{E}(b_i - 1)^2$.
\end{rem}

\begin{cor}[Affine images of regular polygons]\label{prop:affine-hand}
For every $p \ge 3$ and every nonsingular affine transformation $T$ of $\R^2$, one has $J(T R_p) \le J(R_p)$, with equality if and only if $T$ is a similarity.
\end{cor}

\begin{proof}
Because $T R_p$ is a convex $p$-gon, the convex rigidity theorems above give $J(T R_p) \le J(R_p)$, with equality precisely when $T R_p$ is regular. The covariance matrix of the equally weighted vertices of $R_p$ is a positive scalar multiple of the identity; if $T R_p$ is also regular, its vertex covariance matrix has the same property, forcing the singular values of the linear part of $T$ to coincide.
\end{proof}

\section{Fourth-order asymptotics and Blaschke-transported rigidity}\label{sec:asymptotics-blaschke}

Every estimate in this section holds for simple planar polygons with $p$ strictly convex vertices and an \emph{arbitrary} number of reflex vertices. In Subsection~\ref{sec:fourth-order} we prove Theorem~\ref{thm:unrestricted-reflex-asymptotic} by averaging boundary phase differences across cyclic shifts of the convex gaps. In Subsections~\ref{sec:blaschke-moments}--\ref{sec:hausdorff}, we introduce a degree-$p$ Blaschke factorization adapted to the convex prevertices, establish $H^{1/p}$ Sobolev coercivity and arc-mass rigidity for the transported boundary measure, and prove Theorem~\ref{thm:intro-blaschke} together with its sharp refinements.

\subsection{The sharp fourth-order bound at fixed convex-vertex count}\label{sec:fourth-order}

\begin{proof}[Proof of Theorem~\ref{thm:unrestricted-reflex-asymptotic}]
Work in arc-length balanced coordinates and write $y = 1/p$, $u = \max_i s_i$, $M = \sum_{i=1}^p s_i^4$, $T_5 = \sum_{i=1}^p s_i^5$, and $d = d(P)$. Assume first that $d \le 2Z y^3 + 22W y^4$; all implicit constants below are absolute and uniform in $P$.

\paragraph{Gap concentration.}
Put $r_i = s_i/y$. Since $p\mathcal{H}(y) = 2Z y^3 + O(y^4)$, inequality \eqref{eq:gap} implies $Q \coloneqq \sum_{i=1}^p (r_i - 1)^2 = O(1)$, and in particular $u = O(y)$. Expanding $r_i^4$ and $r_i^5$ about $1$ and using $\sum_i (r_i - 1) = 0$ together with $r_i = O(1)$ gives $M = y^3 + O(y^4)$ and $T_5 = y^4 + O(y^5)$.

\paragraph{Cyclic averaging of the phase differences.}
Define the gap-supported profile $f|_{I_i} = s_i^3\Phi$ with $\Phi$ as in \eqref{eq:phi-profile}. Normalize the conformal derivative $F$ so that $|F| = \rho$, and choose its analytic square root $H = F^{1/2}$ with $H(0) > 0$. By Lemma~\ref{lem:boundary-phase-product} and Lemma~\ref{lem:gap}, the common value $\mathcal{P}$ of the two phase integrals in \eqref{eq:phase-product} satisfies
\[
  0 \le \mathcal{P} \le \sqrt{d(c M + b T_5)} - M = O(y^4).
\]
Set $\eta^2 = \|H - 1\|_2^2 = 2(1 - H(0))$, choose the branch of $\log H$ that is real at $0$, and define the good set
\[
  G \coloneqq \Bigl\{z \in \T : \sup_{0 < r < 1} |H(rz) - 1| \le \frac{1}{2}\Bigr\}.
\]
Since $|H(rz) - 1| \le P_r|H - 1|(z)$, the Hardy--Littlewood maximal inequality gives $m(G^c) \lesssim \eta^2$. For every $z \in G$, the radial path $r \mapsto H(rz)$ remains in the disk $\{|w - 1| \le 1/2\}$, so $\log H(z)$ agrees with the principal logarithm. Applying \eqref{eq:hardy-projection} with $\alpha = 1$ expresses the phase integral in terms of $h \coloneqq \Im H$ as
\[
  \mathcal{P} = 4\iint_{\T \times \T} f(z)\frac{(h(z) - h(w))^2}{|z - w|^2}\,dm(z)\,dm(w).
\]

Set $\epsilon = \sqrt{y Q}$. For any signed arc length $s \in (-1, 1)$ in turns, pairing opposite terms in the partial fraction expansion yields
\[
  |1 - e^{2\pi i s}|^{-2} = \frac{1}{4\pi^2}\sum_{k \in \Z} (s + k)^{-2} \ge \frac{1}{(2\pi s)^2} + \frac{1}{12}.
\]
The constant term $1/12$ contributes at least $\frac{1}{3}M\Var_m(h) = \frac{1}{6}M(1 - H(0)^2) \ge \frac{1}{12}M\eta^2$ to $\mathcal{P}$.

Retain only the gaps satisfying $r_i \ge 1/2$ and relabel them in cyclic order. Since each discarded gap contributes at least $1/4$ to $Q$, the number of retained gaps is $\ell = p - O(Q)$. Let $\mathbb{E}_i$ denote the uniform average over these $\ell$ indices, write $U_i(v) = H(z_i(v))$ for the parametrization $z_i : [0,1] \to I_i$, and put $t = w - v$, $d\omega = \Phi(v)\,dv\,dw$, $A_n = (U_i(v) - U_{i+n}(w))/y$, and $X_n = \Im A_n/\pi$. All sums over shifts $n$ below are restricted to $|n| \le \ell/3$.

Let $a$ be the shortest unsigned arc length between $z_i(v)$ and $z_{i+n}(w)$ measured in units of $y$, and define $K_n = r_i^{-4} r_{i+n}^{-1} a^2$. Deleting the intervening discarded short gaps replaces $a$ by $a_0$, where $a_0/|n + t|$ is a convex combination of the retained gap ratios $r_j$. Because the retained ratios are bounded above and below and $\mathbb{E}_i |r_i - 1| \le \sqrt{Q/\ell} = O(\epsilon)$, we obtain $\mathbb{E}_i(r_i^{-4} r_{i+n}^{-1} a_0^2) \le (1 + C\epsilon)(n + t)^2$. The lengths of all discarded gaps sum to $O(Q)$ in units of $y$, and each discarded gap is traversed at most $|n|$ times as $i$ varies, so $0 \le a - a_0 \le C Q$ and $\mathbb{E}_i(a - a_0) \le C|n| y Q$. Since $Q = O(1)$, the difference $(a - a_0)(a + a_0)$ contributes $O(yQ)$ when $|n| = 1$ and $O((n+t)^2 y Q)$ when $|n| \ge 2$ (where $|n+t| \ge |n|/2$). Hence there is an absolute constant $C_0 > 0$ such that $\mathbb{E}_i K_n \le (1 + C\epsilon)D_n$, where $D_n = (n + t)^2 + C_0\epsilon^2\mathbf{1}_{\{|n| = 1\}}$.

Define the cyclic mean profile $h_{\mathrm{av}}(v) = \mathbb{E}_i \Im U_i(v)/(\pi y)$, so that $\mathbb{E}_i X_n = h_{\mathrm{av}}(v) - h_{\mathrm{av}}(w)$ for every shift $n$. Put $S = (h_{\mathrm{av}}(v) - h_{\mathrm{av}}(w))/t$ and $d\lambda = t^2\sum_n D_n^{-1}\,d\omega$; the diagonal $\{t = 0\}$ has measure zero, and the $n = 0$ summand in the density equals $1$. Applying Cauchy--Schwarz in the cyclic average for each $n$ gives
\[
  \frac{\mathcal{P} - M\eta^2/12}{y^4}
  \ge \frac{\ell}{p}\int \sum_n \mathbb{E}_i\frac{X_n^2}{K_n}\,d\omega
  \ge \frac{\ell}{p(1 + C\epsilon)}\int S^2\,d\lambda.
\]
For $n \neq 0$, the only singular corners occur when $|n + t| = x_0 + b_0$, where $x_0, b_0$ are the distances of $v, w$ to the adjacent gap endpoints; since $\Phi(v) \le C x_0^{3/2}$, both $\sum_{n \neq 0} D_n^{-1}\,d\omega$ and $\lambda$ have uniformly bounded marginal densities. Replacing $D_n$ by $(n+t)^2$ changes the total mass of $\lambda$ by $O(\epsilon)$ because the reciprocal difference is $O(\epsilon/(x_0 + b_0)^3)$ and $x_0^{3/2}/(x_0 + b_0)^3$ is integrable. Furthermore, $\int \mathbb{E}_i |A_0/t|^2\,d\lambda = O(1)$: the $n = 0$ term is controlled directly by the same-gap Hardy integral from $\mathcal{P} = O(y^4)$ and $r_i \ge 1/2$, while the remaining terms follow by inserting an independent intermediate point of density $\Phi$ in the same gap and using the bounded marginals.

For $n = 0$, write $z = z_i(v)$ and $\xi = z_i(w)$. On the open gap $I_i$, $\Im\log H$ has angular slope $-1/2$ along straight edges and jumps strictly downward across any interior reflex prevertex; consequently $\Im(\log H(z) - \log H(\xi))/(\pi y t) \ge r_i$ for either sign of $t = w - v$. When $(z, \xi) \in G \times G$, integrating $1/\zeta - 1$ along the straight segment from $H(\xi)$ to $H(z)$ gives
\[
  \frac{X_0}{t} \ge r_i - \frac{2}{\pi}\bigl(|H(z) - 1| + |H(\xi) - 1|\bigr)\left|\frac{A_0}{t}\right|.
\]
Under $\mathbb{E}_i\,d\lambda$, both marginals have bounded densities with respect to $m$. Cauchy--Schwarz therefore bounds the integral of the error term on $G \times G$ by $O(\eta)$, while on the exceptional set $(G \times G)^c$ (which has measure $O(\eta^2)$) the same $L^2$-bound controls the integral of $1 - X_0/t$ by $O(\eta)$. Since $S = \mathbb{E}_i(X_0/t)$ and $\mathbb{E}_i |r_i - 1| = O(\epsilon)$, we obtain $\int (1 - S)\,d\lambda \le C(\epsilon + \eta)$, and hence $\int S^2\,d\lambda \ge \int d\lambda - O(\epsilon + \eta)$ via $S^2 \ge 2S - 1$. Using $\ell/p = 1 - O(yQ)$, noting that the omitted tail $|n| > \ell/3$ contributes $O(1/\ell) = O(y)$, and recalling $M \ge y^3$, we arrive at
\[
  \frac{\mathcal{P}}{y^4}
  \ge \frac{\eta^2}{12y}
    + \int_0^1 \Phi(v)\int_0^1 \sum_{n \in \Z} \frac{(v - w)^2}{(n + w - v)^2}\,dw\,dv
    - O(\eta + \epsilon + y).
\]

To evaluate this double integral, put $L(t) = -\log(2|\sin\pi t|)$. By partial fractions, the inner integral equals $\int_{v-1}^v t^2\pi^2\csc^2(\pi t)\,dt$. Since $L''(t) = \pi^2\csc^2(\pi t)$ off $\Z$ and $L$ is $1$-periodic with mean zero, integration by parts gives $(2v - 1)L'(v) - 2L(v)$ (the boundary terms at $t = 0$ vanish since $t^2 L'(t)$ and $t L(t)$ tend to $0$). Using $\Phi'' = 4\pi^2 L/Z$ on $(0,1)$ with $\Phi(0) = \Phi'(0) = \Phi(1) = \Phi'(1) = 0$, two integrations by parts and Parseval's identity yield
\begin{align*}
  \int_0^1 \Phi(v)\bigl[(2v - 1)L'(v) - 2L(v)\bigr]\,dv
  &= -\int_0^1 \bigl[(2v - 1)\Phi'(v) + 4\Phi(v)\bigr]L(v)\,dv \\
  &= \frac{5Z}{4\pi^2}\int_0^1 (\Phi')^2\,dv
   = \frac{5W}{2Z}.
\end{align*}

\paragraph{Variational comparison.}
Set $h = M - (1 + \Lambda)f$. Then $\bar{h} = 0$, $\mu h \ge M + \mathcal{P}$, and Fourier expansion combined with \eqref{eq:phi-energy-bounds} gives $\E(h) - \|h\|_2^2 \le c M + b T_5/2$. Because $f' = 0$ at every gap endpoint and $\Phi' \in C^{1/2}$, rescaling yields $[f']_{C^{1/2}} = O(u^{3/2})$ globally on $\T$. The symmetric second difference in the singular integral for $\Lambda f$ is therefore $O(\min\{u^3, u^{3/2}|t|^{3/2}\})$; integrating against $O(t^{-2})\,dt$ on both sides of $|t| = u$ gives $\|h\|_\infty = O(u^2) = O(y^2)$ and hence $\E(h) = O(y^3)$. Since $\|\rho - 1\|_1 \le \||H| + 1\|_2\|H - 1\|_2 \le 2\eta$, replacing $m$ by $\mu$ alters $\Var(zh)$ by $O(\eta\|h\|_\infty^2) = O(\eta y^4)$. Moreover, $|\int_\T z h\,dm| = 2|\widehat{f}(1)| \le 2M$. Testing \eqref{eq:variational-disk} with $z(1 + h/c)$ and using $\mu z = \bar{h} = 0$ gives
\[
  d \ge \frac{2\mu h}{c} - \frac{(1 + d)\E(h) - \Var_\mu(zh)}{c^2}
  \ge \frac{M + 2\mathcal{P}}{c} - \frac{b T_5}{2c^2} - O(\eta y^4 + y^6).
\]
Substituting $T_5 = y^4 + O(y^5)$, the lower bound for $\mathcal{P}$, $c = 1/(2Z)$, and $b = W/Z^2$, and noting that $4Z \cdot 5W/(2Z) - 2W = 8W$, we obtain
\[
  d \ge 2Z M + 8W y^4 + \frac{Z}{3}y^3 \eta^2 - C y^4(\eta + \sqrt{yQ} + y).
\]
Using $(M - y^3)/y^4 = \sum_{i=1}^p (r_i - 1)^2(r_i^2 + 2r_i + 3) \ge 3Q$ and completing the square in $\eta$ and $\sqrt{Q}$ yields
\[
  d \ge 2Z y^3 + 8W y^4 + 2y^4 Q + \frac{1}{6}y^3 \eta^2 - O(y^5).
\]

When $d > 2Z y^3 + 22W y^4$, average this threshold bound with Lemma~\ref{lem:gap}: the prefactor $[1 + (1 - H(0)^2)/2]^2$ is at least $1 + \eta^2/2$, while $(p + 2Q)\mathcal{H}(y) \ge 2Z y^3 - 5W y^4 + 4y^4 Q$ and $p\mathcal{H}(y) \ge 2y^3$ for small $y$, recovering the exact same lower bound. Combining with \eqref{eq:regular-two} establishes the stability inequality with $d(R_p)$ and the $L^1$-bound $\|\rho - 1\|_1 \le 2\eta$. Regular polygons provide the matching upper bound for the infimum, and replacing $p$ by $p + 1$ in \eqref{eq:regular-two} reduces the fourth-order coefficient by $6Z$, ensuring $J(P) < J(R_{p+1})$ for all sufficiently large $p$.
\end{proof}

\subsection{Blaschke factorization, transported moments, and \texorpdfstring{$H^{1/p}$}{H\string^(1/p)} coercivity}\label{sec:blaschke-moments}
Let $\xi_1, \dots, \xi_p \in \T$ be the preimages of the $p$ strictly convex vertices of $P$, and define the prevertex polynomial $S(z) = \prod_{i=1}^p (1 - \bar{\xi}_i z)$. Any factorization of the form
\begin{equation}\label{eq:blaschke-factor}
  S(z) = A(z)\bigl(1 - B(z)\bigr),
\end{equation}
where $A$ is holomorphic and zero-free on $\overline{\D}$ and $B$ is a finite Blaschke product of degree $p$ (so that $B$ maps $\T$ onto $\T$ as a $p$-fold orientation-preserving covering sending each $\xi_i$ to $1$), transports the boundary measure $\mu$ to a circle measure whose Fourier moments satisfy exact Fej\'er positivity identities.

\begin{lem}[Fej\'er identity for transported moments]\label{lem:hardy-marker-moments}
Given a factorization \eqref{eq:blaschke-factor}, normalize the conformal derivative $F$ so that $|F| = \rho$, and define
\begin{align*}
  H &= (F S^{2/p})^{1/2}, \qquad u = \frac{|B'|}{p}, \qquad \nu = B_*(u|S|^{2/p}\mu), \\
  \sigma &= B_*\bigl(|1 - B(z)|^2 \|D_z H\|_2^2\,dm(z)\bigr), \qquad
  D_z H(w) = \frac{H(w) - H(z)}{w - z}.
\end{align*}
Then for every integer $j \ge 1$ and the Fej\'er kernel $K_j(\zeta) = j^{-1}|1 + \zeta + \cdots + \zeta^{j-1}|^2$,
\[
  p\int_\T \Re\zeta^j\,d\nu(\zeta) = \int_\T K_j(\zeta)\,d\sigma(\zeta).
\]
Furthermore, let $w = w_{1/p}$ be the probability density proportional to $|1 - \zeta|^{-2/p}$, write $r_j = r_j(1/p)$, and set $a_{jk} = r_{|j-k|} + r_{j+k} - 2r_j r_k$. Then
\begin{equation}\label{eq:positive-tilt}
  \int_\T w(\zeta)(\Re\zeta^j - r_j)\,d\nu(\zeta)
  = \frac{1}{p}\sum_{k=1}^\infty a_{jk}\int_\T K_k(\zeta)\,d\sigma(\zeta) \ge 0.
\end{equation}
Every tilted moment on the left-hand side is strictly positive if $H$ is nonconstant, and vanishes if $H$ is constant.
\end{lem}

\begin{proof}
Write $x = 1/p$ and $w(\zeta) = C_p |1 - \zeta|^{-2x}$, so that $w\,d\nu = C_p B_*(u|A|^{2x}\mu)$. Choose $s' \in (1,\infty)$ satisfying $2x < 1/s' < 1 - \max_i (2b_i - 2x)_+$ (which is nonempty since $2b_i < 1$). Then $w \in L^{s'}(\T)$ while $|H|^2$ lies in the conjugate space $L^s(\T)$, because its only positive singularity exponents are among $2b_i - 2x$.

For any real $g \in L^{s'}(\T)$, put $U = \Lambda^{-1}(g - \bar{g}) \in W^{1,s'}(\T)$. Checking on monomials and extending by density gives the conformal covariance $\Lambda(U \circ B) = p u(g - \bar{g}) \circ B$. Since the turning weights of $H^2 = F S^{2x}$ sum to $1 - p x = 0$ and the function $f(z) = U(1) - U(B(z))$ vanishes at every strictly convex prevertex, Lemma~\ref{lem:boundary-phase-product} (with $\beta = 0$) yields
\begin{equation}\label{eq:anchored-transport-potential}
  p\int_\T (g - \bar{g})\,d\nu
  = 2\int_\T \bigl[U(1) - U(B(z))\bigr]\|D_z H\|_2^2\,dm(z).
\end{equation}
Taking $g(\zeta) = \Re\zeta^j$ gives $U(1) - U(\zeta) = (1 - \Re\zeta^j)/j = \frac{1}{2}|1 - \zeta|^2 K_j(\zeta)$, proving the Fej\'er identity. At $j = 1$, this shows that $\sigma$ is a finite nonnegative measure with $\sigma(\T) > 0$ if and only if $H$ is nonconstant.

Let $X$ be a $\mathrm{Beta}(x, 1 - 2x)$ random variable, so that $\mathbb{E} X^k = r_k$. Positive covariance of increasing functions gives $r_{j+k} \ge r_j r_k$, and since $r_{|j-k|} > r_{j+k}$ we have $a_{jk} > 0$ for all $j, k \ge 1$. For $g(\zeta) = w(\zeta)(\Re\zeta^j - r_j)$, we have $\bar{g} = 0$ and $\widehat{g}(k) = a_{jk}/2$ for $k \ge 1$. By Abel means and monotone convergence, the continuous potential satisfies $U(1) - U(\zeta) = \sum_{k=1}^\infty a_{jk}(1 - \Re\zeta^k)/k$. Substituting into \eqref{eq:anchored-transport-potential} and applying Tonelli's theorem gives \eqref{eq:positive-tilt}, with strict positivity whenever $\sigma(\T) > 0$ since $K_1 \equiv 1$.
\end{proof}

At $j = 1$, the identity \eqref{eq:positive-tilt} evaluates to an exact radial-difference integral that controls the $H^{1/p}$ Sobolev seminorm of $H$ with optimal exponent and optimal constant.

\begin{lem}[First transported moment and $H^{1/p}$ coercivity]\label{lem:radial-marker}
In the notation of Lemma~\ref{lem:hardy-marker-moments}, set $x = 1/p$, $h = H/\sqrt{C_p}$, $V = u|A|^{2x}$, and $\delta = \mu[V(\Re B - r_1)]$. With $X \sim \mathrm{Beta}(x, 1 - 2x)$, one has
\begin{align}
  p(1 - x)\delta
  &= \mathbb{E}\int_\T \left|\frac{(1 + X)(1 - B(z))}{1 - X B(z)}\right|^2 \|D_z h\|_2^2\,dm(z) \label{eq:radial-moment} \\
  &\ge \sum_{n=1}^\infty |\widehat{h}(n)|^2\,\mathbb{E}\!\left[\frac{1 + X}{1 - X}(1 - X^n)^2\right]
   \ge (1 - r_2)\sum_{n=1}^\infty n^{2/p}|\widehat{h}(n)|^2. \label{eq:radial-coercivity}
\end{align}
Identity~\eqref{eq:radial-moment} holds for every factorization \eqref{eq:blaschke-factor}, while \eqref{eq:radial-coercivity} holds whenever $B(0) = 0$. The factor $1 - r_2$ is sharp for the radial multipliers in \eqref{eq:radial-coercivity}, and the Sobolev exponent $1/p$ is sharp even within the class of convex $p$-gons.
\end{lem}

\begin{proof}
Set $c_1 = 1 - r_1 = (1 - 2x)/(1 - x)$. Using $r_{n+1}/r_n = (n + x)/(n + 1 - x)$, we rewrite $a_{1n}$ for $n \ge 1$ as
\begin{equation}\label{eq:first-covariance}
  a_{1n} = \frac{n(r_{n-1} - r_{n+1})}{1 - x} = 2c_1 r_n \frac{n^2}{n^2 - (1 - x)^2}.
\end{equation}
Substituting $|1 - \zeta|^2 K_n(\zeta) = 2(1 - \Re\zeta^n)/n$ into \eqref{eq:positive-tilt} at $j = 1$ and evaluating the nonnegative sum via
\[
  (1 - v^2)\sum_{n=1}^\infty v^{n-1}(1 - \cos nt) = \frac{(1 + v)^2(1 - \cos t)}{|1 - v e^{it}|^2}
\]
at $v = X$ yields \eqref{eq:radial-moment}.

Assume now that $B(0) = 0$. Fix $v \in (0,1)$ and a point $w \in \T$ away from the prevertices. As a function of $z$, $G(z) = (1 + v)(1 - B(z))D_w h(z)/(1 - v B(z))$ lies in $H^2(\D)$. By Schwarz's lemma, $|B(vw)| \le v$; for any complex number $a$ with $|a| \le v$, the triangle inequality gives $|1 - va| \le v|1 - a| + (1 - v) \le (1 + v)|1 - a|$. The Hardy evaluation bound $\|G\|_2^2 \ge (1 - v^2)|G(vw)|^2$ therefore implies
\[
  \|G\|_2^2 \ge \frac{1 + v}{1 - v}|h(w) - h(vw)|^2.
\]
Integrating in $w$, applying Parseval's identity, and taking expectation at $v = X$ establishes the first inequality in \eqref{eq:radial-coercivity}.

Denote the $n$-th radial multiplier by $M_n = \mathbb{E}[(1 + X)(1 - X^n)^2/(1 - X)]$. Expanding $(1 - X^n)/(1 - X)$ and telescoping the blocks $0 \le k < n$ and $n \le k < 2n$ via $x(r_k + r_{k+1}) = (k + 1)r_{k+1} - k r_k$ yields $M_n = (2n/x)(r_n - r_{2n})$. Substituting $t = -n\log X$ into this beta-moment difference gives
\[
  \frac{M_n}{n^{2x}} = \frac{2C}{x}\int_0^\infty \frac{e^{-t} - e^{-2t}}{[2n\sinh(t/(2n))]^{2x}}\,dt,
\]
where $C$ is the beta normalizing constant. Because $u \mapsto \sinh(u)/u$ increases strictly on $(0,\infty)$, $M_n/n^{2x}$ increases strictly with $n$, attaining its minimum at $n = 1$, where $M_1 = 1 - r_2 = 2(1 - 2x)/[(1 - x)(2 - x)]$.

To verify sharpness of the exponent $1/p$, let $P$ be a convex $p$-gon with one arbitrarily small interior angle $\alpha > 0$ (obtained by truncating other corners of a sharp triangle). Near the corresponding prevertex $\xi$, $h(e^{it}\xi)$ behaves like a nonvanishing smooth factor times $(1 - e^{it})^{-b_0}$, where $b_0 = 1/2 - x - \alpha/(2\pi) > 0$. For $s \in (0,1)$, restricting the $H^s$ double integral to $0 < t < \varepsilon$ and $2t < u < 3t$ bounds $\|h\|_{\dot{H}^s}^2$ below by a positive multiple of $\int_0^\varepsilon t^{-2b_0 - 2s}\,dt$, which diverges whenever $s \ge x + \alpha/(2\pi)$. Letting $\alpha \downarrow 0$ proves optimality of $s = 1/p$.
\end{proof}

We next analyze the scalar kernel arising in \eqref{eq:radial-moment} to link $\delta$ with signed tests and total variation. With $X \sim \mathrm{Beta}(x, 1 - 2x)$, set $Y = 4X/(1 - X)^2$ and $\gamma = 1/2 - x$. Since $(1 + Y)^{-1} = ((1 - X)/(1 + X))^2 \sim \mathrm{Beta}(\gamma, x)$, $Y$ has density $p_x(y) = C y^{x-1}(1 + y)^{-1/2}$ on $(0,\infty)$. For $A \in (0, 1]$, define
\[
  m(A) = \mathbb{E}\log(1 + AY), \qquad
  d(A) = \frac{2A}{1 - x}\,\mathbb{E}\!\left[\frac{1 + Y}{1 + AY}\right],
\]
with $m(0) = d(0) = 0$. Since $-2[(1 + y)p_x(y)]' = [1 - 2x + 2(1 - x)/y]p_x(y)$, integrating against $\log(1 + Ay)$ (the boundary terms are $O(y^x)$ at $0$ and $O(y^{x-1/2}\log y)$ at $\infty$) yields, with $c_1 = (1 - 2x)/(1 - x)$,
\begin{equation}\label{eq:arc-kernel-remainder}
  0 \le d(A) - c_1 m(A) = 2\,\mathbb{E}\!\left[\frac{\log(1 + AY)}{Y}\right] \le 2A.
\end{equation}
Changing variables to $t = Ay$ gives
\[
  \frac{m(A)}{C A^\gamma} = \int_0^\infty \frac{t^{x-1}\log(1 + t)}{\sqrt{A + t}}\,dt, \qquad
  \frac{(1 - x)d(A)}{2C A^\gamma} = \int_0^\infty \frac{t^{x-1}\sqrt{A + t}}{1 + t}\,dt.
\]
The first integral decreases strictly in $A$ while the second increases strictly. Hence $m(A) \ge M_x A^\gamma$ with $M_x \coloneqq m(1)$, and the ratio $d(A)/m(A)$ increases strictly from $c_1$ (by \eqref{eq:arc-kernel-remainder}) to $2/[(1 - x)M_x]$ at $A = 1$. Evaluating the logarithmic moment of $(1 + Y)^{-1} \sim \mathrm{Beta}(\gamma, x)$ via Euler's gamma product gives
\[
  M_x = \sum_{j=0}^\infty \frac{x}{(j + 1/2)(j + 1/2 - x)} > \frac{4x}{1 - 2x}.
\]

Putting $b(z) = |1 - B(z)|$ and $\beta = 2x$, the averaged kernel in \eqref{eq:radial-moment} at $A = b(z)^2/4$ equals $2(1 - x)d(A) \ge 2(1 - x)c_1 M_x A^\gamma$. Consequently, for any continuous (possibly complex-valued) $f$ on $\T$ satisfying $|f| \le M_0 b^{1-\beta}$ and distributional bound $|\Lambda f| \lesssim b^{-\beta}$, we have $\Lambda f \in L^{s'}$ (with $s'$ as in \eqref{eq:anchored-transport-potential}), hence $f \in W^{1,s'}$ vanishing at the strictly convex prevertices, and Lemma~\ref{lem:boundary-phase-product} yields
\begin{equation}\label{eq:pinned-test}
  \int_\T |h|^2 \Lambda f\,dm = -2\int_\T f\|D_z h\|_2^2\,dm, \qquad
  \left|\int_\T |h|^2 \Lambda f\,dm\right| \le C_p M_0 \delta.
\end{equation}
Write $L_v(s) = -2\log|s - v|$ and $N(g) = \sup_{\eta \in \T} \int_\T |g(s)|\,|s - \eta|^{-\beta}\,dm(s)$. The logarithmic kernel satisfies
\begin{equation}\label{eq:log-shift}
  N(L_v - L_u) \lesssim_p |v - u|^{1-\beta}.
\end{equation}
Indeed, with $r = |u - v|$, $|L_v(s) - L_u(s)| \le 2\log(1 + r/|s - u|) + 2\log(1 + r/|s - v|)$; by the layer-cake representation on centered superlevel arcs, pairing either term with $|s - \eta|^{-\beta}$ is maximized when the centers coincide, and scaling bounds the integral by $C_p r^{1-\beta}\int_0^\infty t^{-\beta}\log(1 + 1/t)\,dt$.

\begin{cor}[Arc mass controls total variation of the transported measure]\label{cor:arc-order}
In the notation of Lemma~\ref{lem:radial-marker}, set $\omega = w_x\,dm$, $\eta = B_*(V\mu)/(\mu V)$, and $D = \delta/(\mu V) = \int_\T \Re\zeta\,d(\eta - \omega)$. For $I_t = \{e^{is} : |s| \le t\}$, write $a_t = \omega(I_t)$. Then
\begin{equation}\label{eq:arc-order}
  \frac{1}{2}\min(a_t, 1 - a_t)\,D
  \le \eta(I_t) - a_t
  \le \|\eta - \omega\|_{\mathrm{TV}}
  \le C_p D, \qquad 0 \le t \le \pi,
\end{equation}
where $\|\cdot\|_{\mathrm{TV}}$ denotes the supremum against real measurable tests bounded in modulus by $1$. In particular, for any fixed $t \in (0,\pi)$, equality $\eta(I_t) = a_t$ holds if and only if $P$ is convex and equiangular. When the convex prevertices are equally spaced, taking $B(z) = z^p$ and $V \equiv 1$ (after rotation) shows that $\eta(I_t)$ is the fraction of $|\partial P|$ whose preimage lies within angular distance $t/p$ of a convex prevertex, $a_t$ is that fraction for $R_p$, and equality characterizes $R_p$ up to similarity.
\end{cor}

\begin{proof}
For $v \in \T$, put $g_v = L_1 - L_v$, $m_v = \omega g_v$, and $d\lambda_v = (g_v - m_v)\,d\omega$. By \eqref{eq:log-shift},
\begin{equation}\label{eq:anchor-modulus}
  \|\lambda_v\|_{\mathrm{TV}} \le 2\int_\T |g_v|\,d\omega \lesssim_p |1 - v|^{1 - 2x}.
\end{equation}
Applying \eqref{eq:anchored-transport-potential} to $g = w_x(G - \omega G)$ for any bounded real $G$ gives $U(1) - U(v) = \int G\,d\lambda_v$ by logarithmic convolution. Since $w_x\,d\nu = C_p B_*(V\mu)$, this yields the signed-measure identity
\begin{equation}\label{eq:anchor-pairing}
  p(\mu V)(\eta - \omega) = 2\int_\T \lambda_{B(z)}\|D_z h\|_2^2\,dm(z).
\end{equation}
Combining \eqref{eq:anchor-modulus} with the scalar kernel bound preceding \eqref{eq:pinned-test} gives $\|\eta - \omega\|_{\mathrm{TV}} \le C_p D$.

For centered arcs and cosines, symmetrize $g_v$: writing $v = e^{i\theta}$ and $A = \sin^2(\theta/2)$, the even part of $g_v(e^{is})$ is $k(s) = \log|1 - A/\sin^2(s/2)|$. When $v \neq 1$, $\max(k(s), -A)$ decreases in $|s|$: it decreases on $|s| \le |\theta|$, and for $|s| > |\theta|$ we have $k(s) \le \log(1 - A) \le -A$. Set $T = \int (k - m_v)_+\,d\omega$ and $d_v = \int \Re\zeta\,d\lambda_v$; clearly $d_v \le 2T$. With $e_n = 1 - \cos(n\theta)$, summing $m_v = 2\sum_{n \ge 1} r_n e_n/n$ via the logarithm series and $d_v = \sum_{n \ge 1} a_{1n} e_n/n$ via the geometric sum of Lemma~\ref{lem:radial-marker} identifies $m_v = m(A)$ and $d_v = d(A)$ in \eqref{eq:arc-kernel-remainder}.

Now fix $t \in (0,\pi)$. If $k(t) \ge m_v$, monotonicity of $k$ on its positive set gives $\lambda_v(I_t) \ge a_t T$. If $-A < k(t) < m_v$, then on $I_t^c$ the average of $m_v - k$ is at least $T$ because removing the negative-density part of $I_t^c$ leaves larger values, so $\lambda_v(I_t) = -\lambda_v(I_t^c) \ge (1 - a_t)T$. Finally, if $k(t) \le -A$, then $m_v - k \ge m_v + A$ throughout $I_t^c$, giving $\lambda_v(I_t) \ge (1 - a_t)(m_v + A)$. Since $d(A) \le m(A) + 2A \le 2(m(A) + A)$ by \eqref{eq:arc-kernel-remainder} and $d_v \le 2T$, in every case we obtain
\[
  \lambda_v(I_t) \ge \min(a_t, 1 - a_t)\min(T, m_v + A) \ge \frac{1}{2}\min(a_t, 1 - a_t)\,d_v.
\]
The inequality is immediate for $v = 1$ or $t \in \{0, \pi\}$. Integrating \eqref{eq:anchor-pairing} proves the lower bound in \eqref{eq:arc-order}. Conversely, $D = \int_0^\pi (\eta(I_t) - a_t)\sin t\,dt$. By \eqref{eq:positive-tilt}, $D = 0$ holds if and only if $H$ is constant, i.e., $F \propto S^{-2/p}$, meaning that $P$ has turning weight $1/p$ at each of the $p$ convex prevertices and no other vertices. Thus $P$ is a convex equiangular $p$-gon, and a regular $p$-gon when the prevertices are equally spaced.
\end{proof}

\subsection{Reciprocal-density transport and canonical factorization}\label{sec:reciprocal-transport}
For $\beta = 1/p$, write $w = w_\beta$ and $r_j = r_j(\beta)$. Given coefficients $b_j \ge 0$ with $\sum_{j \ge 1} b_j < \infty$ and $\sum_{j \ge 1} j b_j^2 < \infty$, set $F_b(e^{i\phi}) = 2\sum_{j=1}^\infty b_j(r_j - \cos j\phi)$. Termwise integration of \eqref{eq:positive-tilt} and Parseval's identity yield
\begin{equation}\label{eq:regular-potential-transport}
  \int_\T u|A|^{2/p} F_b(B)\,d\mu \le 0, \qquad
  \int_\T F_b w\,dm = 0, \qquad
  \E(F_b) = 2\sum_{j=1}^\infty j b_j^2.
\end{equation}
Taking $b_j = r_j/j$ produces the Dirichlet potential
\begin{equation}\label{eq:regular-boundary-potential}
  F_p(e^{i\phi}) = 2S(\beta) - 2\sum_{j=1}^\infty \frac{r_j}{j}\cos(j\phi), \qquad \Lambda F_p = 1 - w,
\end{equation}
which achieves equality in the Dirichlet Cauchy--Schwarz estimate for $g \mapsto \int g(1 - w)\,dm$, and satisfies $F_p(e^{i\phi}) = 2S(\beta) + \mathbb{E}\log[(1 - X)^2 + 4X\sin^2(\phi/2)]$.

More generally, for $W = P_r(\zeta) = (1 - r^2)/|1 - r\zeta|^2$ or $W = w_\beta$, setting $U = \Lambda^{-1}(W - 1)$ and $K_W = U(1) - U \ge 0$ (specifically $K_{w_\beta} = F_p - F_p(1)$ and $K_{P_r}(\zeta) = 2\log(|1 - r\zeta|/(1 - r))$) in \eqref{eq:anchored-transport-potential} gives
\begin{equation}\label{eq:direct-poisson-defect}
  p\left(\int_\T W\,d\nu - \nu(\T)\right) = 2\int_\T K_W(B(z))\|D_z H\|_2^2\,dm(z) \ge 0.
\end{equation}
For the low-side comparisons in Section~\ref{sec:low-reflex}, we specialize \eqref{eq:regular-potential-transport} to the reciprocal-density test
\begin{equation}\label{eq:reciprocal-test}
  F = w_{-\beta} - D_\beta, \qquad D_\beta = \pi\beta\cot(\pi\beta), \qquad
  C \coloneqq \bar{F} = 1 - D_\beta, \qquad \E(F) = \frac{\beta}{2},
\end{equation}
corresponding to $b_j = -r_j(-\beta) > 0$. Here, for $X = |1 - z|^2$ with $z$ uniform on $\T$, the beta integral gives $\mathbb{E} X^t = 4^t \Gamma(t + 1/2)/[\sqrt{\pi}\,\Gamma(t + 1)]$ for $t > -1/2$, so Euler's reflection formula yields $w_\beta w_{-\beta} = (\mathbb{E} X^\beta\,\mathbb{E} X^{-\beta})^{-1} = \pi\beta\cot(\pi\beta) = D_\beta$. Writing $J = H^2$, $\omega = J/|J|$, $w_\beta = C_\beta|1 - \zeta|^{-2\beta}$, and $V = u|A|^{2\beta}$, combining $w_\beta w_{-\beta} = D_\beta$ and $w_\beta\,d\nu = C_\beta B_*(V\mu)$ with Lemma~\ref{lem:boundary-phase-product} gives
\[
  -\int_\T V F(B)\,d\mu
  = \frac{D_\beta}{p C_\beta}\iint_{\T \times \T} K_{w_\beta}(B(z))|J(w)|\left|\frac{\omega(w) - \omega(z)}{w - z}\right|^2 dm(w)\,dm(z).
\]
Furthermore, since $K_{w_x}(v) = m_v$, \eqref{eq:anchor-pairing} and \eqref{eq:direct-poisson-defect} express $p\delta$ and $p L/D_x$ (where $L = -\mu[V(w_{-x}(B) - D_x)]$) as integrals of $d_{B(z)}$ and $m_{B(z)}$ against $2\|D_z h\|_2^2\,dm(z)$, so the pointwise ratio bounds following \eqref{eq:arc-kernel-remainder} imply
\[
  \frac{1 - 2x}{(1 - x)D_x} L \le \delta \le \frac{2}{(1 - x)D_x M_x} L.
\]

The Dirichlet energy in \eqref{eq:reciprocal-test} is a special case of a general polarized identity obtained by telescoping the recurrence for $r_j$. For $a, b < 1/2$ with $a + b < 0$ and any integer $m \ge 0$, putting $V_j = (j - a)(j + m - b)r_j(a)r_{j+m}(b)$ gives $V_j \to 0$, $V_1 = a(m + b)r_m(b)$, and $V_{j+1} - V_j = 2[(a + b)j + am]r_j(a)r_{j+m}(b)$, whence
\begin{equation}\label{eq:polarized-energy}
  \sum_{j=1}^\infty \bigl[(a + b)j + am\bigr] r_j(a) r_{j+m}(b) = -\frac{1}{2}a(m + b)r_m(b).
\end{equation}
Setting $m = 0$, $a = -u$, and $b = -v$ ($u, v > 0$) yields the mixed energy formula
\begin{equation}\label{eq:reciprocal-mixed-energy}
  \E(w_{-u}, w_{-v}) = \frac{uv}{u + v} \quad (u, v > 0), \qquad \E(F) = \frac{\beta}{2}.
\end{equation}

Finally, for any $p$-tuple of distinct points $(\xi_1, \dots, \xi_p)$ on $\T$, there is a canonical factorization \eqref{eq:blaschke-factor} with $B(0) = 0$: set
\begin{equation}\label{eq:canonical-blaschke}
  S(z) = \prod_{i=1}^p (1 - \bar{\xi}_i z), \qquad
  A(z) = S(z) - \frac{z S'(z)}{p}, \qquad
  B(z) = 1 - \frac{S(z)}{A(z)}.
\end{equation}
Because
\[
  \frac{2A(z)}{S(z)} - 1 = \frac{1}{p}\sum_{i=1}^p \frac{1 + \bar{\xi}_i z}{1 - \bar{\xi}_i z}
\]
has strictly positive real part on $\D$ and purely imaginary boundary values away from the prevertices, while $A(\xi_i) = -\xi_i S'(\xi_i)/p \neq 0$, the polynomial $A$ is zero-free on $\overline{\D}$ and $B$ is a degree-$p$ Blaschke product satisfying $B(0) = 0$ and $B(\xi_i) = 1$. At the regular tuple of $p$-th roots of unity, $A \equiv u \equiv V \equiv 1$ and $B(z) = z^p$.

\subsection{Equally spaced prevertices and quadratic Hausdorff stability}\label{sec:hausdorff}
Define the regular spectral moment
\[
  D_p \coloneqq \sum_{l=1}^\infty \frac{2pl}{p^2 l^2 - 1} r_l(1/p)^2 = \|v\|^2,
\]
where $v$ is the coupling column of the regular $p$-gon Fourier matrix from Section~\ref{sec:regular}.

\begin{thm}[Comparison at equally spaced convex prevertices]\label{thm:symmetric}
Suppose the $p$ convex prevertices of $P$ are equally spaced on $\T$ (allowing any number of reflex vertices). Then $J(P) \le J(R_p)$, with equality precisely when $P$ is similar to $R_p$, in either of the following situations:
\begin{enumerate}[label=\textup{(\roman*)},leftmargin=*]
\item the disk coordinates are arc-length balanced and $p \ge 4$;
\item $p \ge 3$ and $P$ admits a nonidentity rotational symmetry.
\end{enumerate}
Moreover, if a sequence of simple polygons $P_k$ with $p$ strictly convex vertices has convex prevertices converging to an equally spaced $p$-tuple under either hypothesis, then $\liminf_{k \to \infty} d(P_k) \ge d(R_p)$. For any fixed $p \ge 3$, if the convex prevertices converge to an equally spaced $p$-tuple in arc-length balanced coordinates, then
\[
  \liminf_{k \to \infty} d(P_k) \ge d_p - \frac{\alpha_p}{2 + \alpha_p}(d_p - D_p),
  \qquad d_p = d(R_p), \quad \alpha_p = 1 - a_0^2,
\]
where $a_0$ is the zeroth coordinate of the positive unit maximizing eigenvector $(a_i)_{i \in \Z}$ of the regular $p$-gon Fourier matrix $M$.
\end{thm}

\begin{proof}
Write $x = 1/p$. Passing to a subsequence achieving a finite $\liminf$ $d$ of $d(P_k)$ and a further subsequence on which $\mu_k \rightharpoonup \mu$ weakly, rotate so that the limiting tuple is the $p$-th roots of unity. In the canonical factorization \eqref{eq:canonical-blaschke}, we have $A_k, V_k \to 1$ and $B_k \to z^p$ uniformly on $\T$. By \eqref{eq:positive-tilt}, $\mu_k[V_k(\Re B_k^l - r_l)] \ge 0$ for every $l \ge 1$; passing to the limit gives $\delta_l \coloneqq \Re\widehat{\mu}(lp) - r_l \ge 0$. For a single polygon with an exact regular tuple, Lemma~\ref{lem:hardy-marker-moments} gives $\delta_1 > 0$ unless $F(1 - z^p)^{2/p}$ is constant, i.e., unless $P$ is similar to $R_p$.

Set $j = 1 + d_p$. Write the unit-energy regular first eigenfield as $T(z) = z\psi(z^p)$, where $\psi(\zeta) = \sum_{i \in \Z} c_i \zeta^i$ and $c_i = a_i/\sqrt{\ell_i} > 0$. For the truncation $T_N$ at $|i| \le N$, positivity of $c_i c_l$ and $\delta_l \ge 0$ give $\mu|T_N|^2 \ge \mu_{R_p}|T_N|^2$. Since the variational inequality passes to the weak limit $\mu$ on trigonometric polynomials, it remains to control the mean $\mu T_N$ in $(1 + d)\E(T_N) \ge \Var_\mu(T_N)$.

In case~(ii), passing to a subsequence with a fixed nonzero cyclic shift $\ell_0 \in \{1, \dots, p-1\}$ induced by the rotational symmetry (which acts continuously on disk automorphisms determined by three boundary points), the limit measure $\mu$ is invariant under multiplication by $e^{2\pi i \ell_0/p}$. Hence $\mu T_N = 0$ for all $N$, yielding $d \ge d_p$. For an exact regular tuple of a nonregular polygon, $\delta_1 > 0$ leaves a fixed strictly positive excess from the modes $i = 0, 1$ as $N \to \infty$. Note also that a $p$-fold rotational symmetry automatically places the $p$ convex prevertices at an equally spaced orbit in coordinates centered at the fixed point.

In the balanced case, $\mu z = 0$. Allocating the frequencies $lp$ ($l \ge 1$) in \eqref{eq:spectralmom} (passed to the weak limit on finite frequency sets first) bounds the sum of $2n(n^2 - 1)^{-1}|\widehat{\mu}(n)|^2$ over all positive frequencies not divisible by $p$ by $d - D_p$. The frequencies $1 + pi$ ($i \neq 0$) have distinct absolute values $\ell_i \ge 2$, all coprime to $p$. Cauchy--Schwarz therefore yields
\[
  |\mu T_N|^2 \le \frac{d - D_p}{2}\sum_{0 < |i| \le N} a_i^2(1 - \ell_i^{-2})
  \le \frac{\alpha_p}{2}(d - D_p).
\]
Letting $N \to \infty$ gives $d \ge d_p - \frac{\alpha_p}{2}(d - D_p)$, which rearranges to the stated general lower bound for all $p \ge 3$.

To upgrade this to the sharp comparison $d \ge d_p$ for balanced coordinates when $p \ge 4$, assume $\lambda = 1 + d \le j$. The mean and variance bounds give (in any scalar normalization of $T$) $\mu|T|^2 \le (j + d/2)\E(T)$. Write $m_\psi = \min_\T |\psi|$; we will verify $m_\psi > 0$ below, and introduce the reciprocal test field $Y(z) = 1/\overline{\psi(z^p)}$. Because $T$ and $Y$ occupy disjoint Fourier residue classes modulo $p$, we have $\E(T, Y) = 0$, while $T\bar{Y} = z$ implies $\mu(T\bar{Y}) = 0$. By Cauchy--Schwarz and the Dirichlet difference kernel,
\[
  \mu|Y|^2 \ge (\mu|T|^2)^{-1}, \qquad \E(Y) \le \frac{p\E(\psi)}{m_\psi^4}.
\]
We verify that $\mu|Y|^2 > j\E(Y)$ separately for $p \ge 5$ and $p = 4$.

For $p \ge 5$, normalize $c_0 = 1$ and put $b_m = -r_m(-x)$ and $q_p = (p + 1)/(p - 1)^2 \le 3/8$. The Schur weight $h_i = \sqrt{\ell_i}\,b_{|i|}$ from Section~\ref{sec:regular} satisfies $M_* h = v \le q_p h$, so $M_* v \le q_p v$. By the Neumann series, $\sqrt{\ell_i}\,c_i = [(j I - M_*)^{-1} v]_i \le v_i/(j - q_p)$, which yields $c_i \le \kappa b_{|i|}$ with $\kappa = q_p/(j - q_p) \le 3/5$, and the zeroth eigenvector equation gives
\[
  U \coloneqq \sum_{i \neq 0} \ell_i c_i^2
  \le \frac{\sum_{i \neq 0} r_{|i|} c_i}{j - q_p} = \frac{d_p}{j - q_p} < \frac{4}{49},
\]
using $d_p \le d_5 < 1/18$ (Lemma~\ref{lem:regular-six-short}). Since $w_{-x}(1) = 0$, its absolutely convergent Fourier series gives $\sum_{i \neq 0} b_{|i|}(1 - \cos i\theta) = w_{-x}(e^{i\theta})$. Using $\cos t \ge -(1 - \cos t)/2$ and $\max w_{-x} \le 2^{2x} < 4/3$ (from Jensen's inequality and $\int \log|1 - z|\,dm = 0$), we obtain $\Re\psi \ge 1 - \kappa w_{-x}/2 > 3/5$. Because $p|i| \le \frac{5}{4}\ell_i$ and $\E(T) = 1 + U$, it follows that
\[
  j\E(Y) < \frac{19}{18}\cdot\frac{5}{4}\cdot\frac{4}{49}\cdot\Bigl(\frac{5}{3}\Bigr)^4 < \frac{5}{6}, \qquad
  \mu|T|^2 < \frac{13}{12}\cdot\frac{53}{49} < \frac{6}{5},
\]
and hence $\mu|Y|^2 > 5/6 > j\E(Y)$.

For $p = 4$, take the square eigenfield $T = H \circ f$ up to rotations, where $f : \D \to (-1,1)^2$ is conformal with $f(0) = 0$ and $H(w) = \Re\sin(\tau w) + i\Im\sinh(\tau w)$ with $\tau, t = \tanh\tau = \cot\tau$ as in Lemma~\ref{lem:reference-square}, so $j = \pi/(4\tau t)$. On the side $w = 1 + iv$ ($|v| \le 1$), $|H(w)|^2 = \sin^2\tau\cosh^2(\tau v) + \cosh^2\tau\sin^2(\tau v)$ has minimum $\sin^2\tau$ and average $\sin^2\tau\cosh^2\tau$. Thus $j\E(T) = (1 - t^4)^{-1}$ and $m_\psi^2 = (1 - t^2)j\E(T)$. Multiplying $T$ by $\bar{z}$ lowers the absolute frequency of each holomorphic mode by $1$ and raises that of each antiholomorphic mode by $1$; since the holomorphic and antiholomorphic parts of $T$ are $\frac{1}{2}(\sin + \sinh)(\tau f)$ and $\frac{1}{2}\overline{(\sin - \sinh)(\tau f)}$, we have
\[
  \E(T) - 4\E(\psi) = \int_\T \Re\bigl(\sin(\tau f)\overline{\sinh(\tau f)}\bigr)\,dm = 1.
\]
To check the last equality, put $w = x + iy$ and $L(s) = \cos s\cosh s - \sin s\sinh s$. Then
\[
  \Re\bigl[\cos(\tau w)\cosh(\tau w) - \sin(\tau w)\overline{\sinh(\tau w)}\bigr]
  = L(\tau x)L(\tau y),
\]
which vanishes on $\partial(-1,1)^2$ since $L(\pm\tau) = 0$. Hence the integrand agrees on $\T$ with the harmonic function $\Re(\cos(\tau f)\cosh(\tau f))$, whose value at $0$ is $1$. By Lemma~\ref{lem:reference-square}, $1.138 < j < 1.143$ and $t < 0.735$, giving
\[
  \bigl[1 - j(1 - t^4)\bigr]\Bigl(1 + \frac{d}{2j}\Bigr) < 0.195(1.063) < 0.21 < (1 - t^2)^2,
\]
which rearranges to $\mu|Y|^2 > j\E(Y)$.

Both $T$ and $Y$ are continuous with finite energy, so their variational inequalities pass to $\mu$ via Fej\'er approximation. Because $\mu(T\bar{Y}) = \E(T, Y) = 0$, both forms are diagonal on $\operatorname{span}\{T, Y\}$. Combining $\mu|T|^2 \ge j\E(T)$ with $\mu|Y|^2 > j\E(Y)$ and choosing a nonzero linear combination of $T$ and $Y$ with zero $\mu$-mean yields $\lambda \ge j$. For a single nonregular polygon with equally spaced prevertices, the strict inequality $\mu|T|^2 > j\E(T)$ forces $\lambda > j$, whereas $R_p$ attains equality.
\end{proof}

\begin{cor}[Quadratic Hausdorff stability and the marked logarithmic model]\label{cor:sharp-gap-stability}
Normalize $P$ to unit perimeter and work in arc-length balanced coordinates. Put
\[
  e^2 = \sum_{i=1}^p \Bigl(s_i - \frac{1}{p}\Bigr)^2, \qquad
  H_p(P) = \inf_g d_{\mathrm{H}}(\partial P, g\partial R_p),
\]
where $d_{\mathrm{H}}$ is Hausdorff distance, $R_p$ has perimeter $1$, and $g$ ranges over rigid motions.
\begin{enumerate}[label=\textup{(\roman*)},leftmargin=*]
\item For each fixed $p \ge 4$, there exist constants $c_p, K_p > 0$ such that
\[
  d(P) - d(R_p) \ge c_p H_p(P)^2 - K_p e^2.
\]
In particular, there is $\gamma_p > 0$ such that $d(P) \ge d(R_{p+1}) + \gamma_p$ whenever $p \ge 4$ and $e$ is sufficiently small.
\item At $p = 3$, if $P$ admits a nonidentity rotational symmetry or if $\delta = \mu[V(\Re B - r_1)]$ is sufficiently small, then the stronger linear bound $d(P) - d(R_3) \ge c_3 H_3(P) - K_3 e$ holds. For all balanced triples, $d(P) \ge D_3 - K_3 e^2$.
\item For convex prevertex tuples sufficiently close to the regular $p$-tuple, the $L^1$-distance of $\rho$ to the $(p-3)$-dimensional affine space $\mathcal{A}$ in \eqref{eq:finite-marker-model} is comparable to $\delta$. For an equally spaced balanced triple oriented at the cube roots of unity, $\delta \le \|\rho - w_{1/3}(z^3)\|_1 \le C\delta$.
\end{enumerate}
\end{cor}

\begin{proof}
Fix $p \ge 3$, put $x = 1/p$, $\beta = 2x$, and $j = 1 + d(R_p)$, and assume first that $e$ is small. Using the canonical factorization \eqref{eq:canonical-blaschke}, rotate so that $\xi_1 = 1$, and set
\[
  q = C_p |S|^{-\beta} = w(B)|A|^{-\beta}, \qquad
  h = (F S^\beta/C_p)^{1/2}, \qquad
  \delta = \mu[V(\Re B - r_1)],
\]
where $|F| = \rho = q|h|^2$. The factors $A, B, u, V$ depend smoothly on the gap lengths and reduce to $A = u = V = 1$, $B = z^p$ at $e = 0$; in particular, $\delta$ is uniformly bounded.

\paragraph{Weighted coercivity.}
Write $k = h - h(0)$. By \eqref{eq:radial-coercivity}, $\|k\|_{\dot{H}^x}^2 \lesssim_p \delta$. We claim that the corner weight $q$ is likewise controlled:
\begin{equation}\label{eq:weighted-tilt-coercivity}
  \int_\T q|k|^2\,dm \lesssim_p \sum_{n=1}^\infty n^{2x}|\widehat{h}(n)|^2 \lesssim_p \delta.
\end{equation}
For separated prevertices, $q \lesssim_p \sum_i w_x(\bar{\xi}_i z)$. Each rotated weight has Fourier matrix bounded entrywise in modulus by $K_{nm} = r_{|n-m|}(x)$ on $\Z$. Setting $a = (1 - 2x)/4$ and $v_n = r_{|n|}(a) > 0$, the identity $w_x w_a = D w_{x+a}$ ($D = \int w_x w_a\,dm$, justified in $L^1$ since $x + a < 1/2$) gives
\[
  \frac{(Kv)_n}{v_n} = D\frac{r_{|n|}(x + a)}{r_{|n|}(a)} \lesssim_p (1 + |n|)^{2x}
\]
by the rising factorial formula. Applying the quadratic form bound from Section~\ref{sec:regular} to $|\widehat{k}(n)|$ and passing from polynomials to $k$ proves \eqref{eq:weighted-tilt-coercivity}.

Let $\theta \mapsto z(\theta)$ be the circle diffeomorphism with $z(0) = 1$ and $B(z(\theta)) = e^{ip\theta}$; it varies smoothly with the gaps and equals $e^{i\theta}$ at $e = 0$. In the $\theta$-coordinate, the singular factor $w(e^{ip\theta})$ of $q$ is fixed, so $Q \coloneqq \int_\T q\,dm$ depends smoothly on the gaps; by cyclic symmetry on zero-sum gap variations, its first derivative vanishes at $e = 0$, giving $Q = 1 + O_p(e^2)$. Also $\int_\T qV\,dm = 1$ since $B_*(u\,dm) = dm$.

Put $G_0 = C_p S^{-\beta}/Q$ and pick $c_0 \in \C$ with $|c_0| = Q^{-1/2}$ and the same phase as $h(0)$. Since $\|h\|_{L^2(q)} = \|c_0\|_{L^2(q)} = 1$, the reverse triangle inequality gives $\|h - c_0\|_{L^2(q)} \le 2\|k\|_{L^2(q)}$. Factoring $h^2 - c_0^2$ and using Cauchy--Schwarz yields, with $\alpha = c_0^2/|c_0|^2$,
\[
  \|\rho - q/Q\|_1 \le \|F - \alpha G_0\|_1 \lesssim_p \sqrt{\delta}.
\]
Because $\|V - 1\|_\infty = O_p(e)$, integrating against $V - 1$ gives
\begin{equation}\label{eq:sharp-gap-mean}
  |\mu V - 1| \le |Q^{-1} - 1| + O_p(e\sqrt{\delta}) \lesssim_p e^2 + e\sqrt{\delta}.
\end{equation}

\paragraph{The marked logarithmic model.}
Let $T(z) = z\psi(z^p)$ be the unit-energy regular first eigenfield and define the transplanted field $T_e(z) = z\psi(B(z))$. Since $B(0) = 0$, $z B^k$ is holomorphic for $k \ge 0$ and antiholomorphic for $k < 0$; within either class its energy pairing with $z B^l$ equals $\pm\int_\T (1 + kpu)B^{k-l}\,dm = |1 + pk|\mathbf{1}_{k=l}$ (as $B_* dm = B_*(u\,dm) = dm$), while opposite classes are orthogonal. Thus $\E(T_e) = \E(T) = 1$.

Write $L_i = L_{\xi_i}$, and let $X$ be the real vector space of signed densities $q(a_0 + \sum_{i=1}^p a_i L_i)$ with $\sum_{i=1}^p a_i = 0$. Let $\mathcal{A} \subset X$ be the affine subspace of densities having total mass $1$ and zero first moment ($\int z\,m\,dm = 0$). We establish
\begin{equation}\label{eq:finite-marker-model}
  \dim \mathcal{A} = p - 3, \qquad \dist_{L^1}(\rho, \mathcal{A}) \asymp_p \delta.
\end{equation}
Let $\xi(z)$ denote a nearest convex prevertex to $z$, set $K(z) = \|D_z h\|_2^2$, and define
\[
  E(s) = -2\int_\T \bigl[L_z(s) - L_{\xi(z)}(s)\bigr] K(z)\,dm(z).
\]
By \eqref{eq:log-shift} and the kernel estimate preceding \eqref{eq:pinned-test}, $N(E) \lesssim_p \delta$. For smooth $f$, the linear functional $\int |h|^2 \Lambda f\,dm + 2\int [f(z) - f(\xi(z))]K(z)\,dm(z)$ vanishes whenever $f$ vanishes at all prevertices $\xi_i$ (by \eqref{eq:pinned-test}) and also at $f \equiv 1$, so it equals $\sum_{i=1}^p a_i f(\xi_i)$ for real $a_i$ summing to $0$. Applying $\Lambda^{-1}$ on nonconstant modes yields
\[
  |h|^2 = \int_\T |h|^2\,dm + \sum_{i=1}^p a_i L_i + E,
\]
with all Fubini interchanges justified by $N(E) < \infty$. Since $q \lesssim_p \sum_i |z - \xi_i|^{-\beta}$, multiplying by $q$ produces an element of $X$ within $O_p(\delta)$ of $\rho$ in $L^1$.

Distinct logarithmic singularities imply $\dim X = p$. At $e = 0$, we have $\int q\,dm = 1$, $\int zq\,dm = 0$, equal masses $\int q L_i\,dm$, and $\int zq L_i\,dm = \xi_i R(1)$ with $R(1) = \sum_{n \in \Z} r_{|n|}/|1 + pn| > 0$. Since the $p$-th roots of unity span $\R^2$ affinely, the three real moment constraints have full rank at $e = 0$ and vary by $O_p(e)$ in the $\theta$-coordinate, so they admit a uniformly bounded right inverse for small $e$. Adjusting the coefficients by $O_p(\delta)$ produces $m \in \mathcal{A}$ with $N(|h|^2 - m/q) \lesssim_p \delta$, proving $\dim \mathcal{A} = p - 3$ and the upper bound in \eqref{eq:finite-marker-model}.

Conversely, for every $m = q(a_0 + \sum_i a_i L_i) \in X$, the pushforward $B_*(Vm\,dm)$ is a constant multiple of $w\,dm$. Indeed, $qV = u w(B)$, and the roots of $(1 - \zeta)A - S$ (whose leading coefficient has modulus $1$ since $\deg A < p$) satisfy $\prod_{B(z) = \zeta} |z - \xi_i| = |1 - \zeta|\,|A(\xi_i)|$; averaging $L_i$ over the fiber and using $\sum_i a_i = 0$ gives $B_*(Vm\,dm) = [a_0 - (2/p)\sum_i a_i\log|A(\xi_i)|]w\,dm$. Consequently $\delta = \int V(\Re B - r_1)(\rho - m)\,dm$, which proves the lower bound in \eqref{eq:finite-marker-model}.

Because $\|m\|_1 \lesssim_p 1$ and the basis of $X$ stays linearly independent in the $\theta$-coordinate, the coefficients of $m$ are uniformly bounded. At $e = 0$, $K_T = \Lambda^{-1}(qT) = jT$ and $K_z = \Lambda^{-1}(qz)$ take values $j\xi_i\psi(1)$ and $\xi_i R(1)$ at the prevertices (here $T$ is continuous since $qT \in L^r$ for some $r > 1$ by Sobolev embedding, and its positive Fourier coefficients are summable by Abel convergence at $1$). With $\kappa = j\psi(1)/R(1)$, the function $T - \kappa z$ annihilates the regular space $X$. In the $\theta$-coordinate, $T_e = z(\theta)\psi(e^{ip\theta})$ differs uniformly from $T$ by $O_p(e)$, so the integrals of $T_e - \kappa z$ against the basis of $X$ are $O_p(e)$. Using $\mu z = 0$ and $\|\rho - m\|_1 \lesssim_p \delta$ gives
\[
  |\mu T_e| = \left|\int_\T (T_e - \kappa z)\rho\,dm\right| \lesssim_p e + \delta.
\]
At $p = 3$, the model $m_e \in \mathcal{A}$ is unique; correcting the first moment $\int zq/Q\,dm = O(e)$ shows $\|m_e - q/Q\|_1 = O(e)$ and hence $\|\rho - q/Q\|_1 \lesssim e + \delta$. At $e = 0$, $m_e = q = w_{1/3}(z^3)$, and testing against $\Re z^3$ gives the reverse bound $\delta \le \|\rho - w_{1/3}(z^3)\|_1$.

Retaining one adjacent mode pair in \eqref{eq:positive-tilt} and using positivity of the Fourier coefficients of $T$ yields $\mu[V(|T_e|^2 - j)] \ge c_p\delta$. The scalar integral $\int q(1 - V)(|\psi(B)|^2 - j)\,dm$ is $O_p(e^2)$ by smooth dependence and cyclic symmetry, and replacing $q/Q$ by $\rho$ incurs an error $O_p(e\sqrt{\delta})$. Since $\E(T_e) = 1$, the variational inequality $(1 + d(P))\E(T_e) \ge \mu|T_e|^2 - |\mu T_e|^2$ gives
\[
  d(P) - d(R_p) \ge c_p \delta - C_p(e^2 + \delta^2),
\]
whence $d(P) - d(R_p) \ge c'_p\delta - C'_p e^2$ for sufficiently small $\delta$. When $p \ge 4$, or $p = 3$ with a nonidentity rotational symmetry, compactness extends this bound to all $\delta$: otherwise, a sequence with $e \to 0$, $\delta \ge \delta_0 > 0$, and $\limsup d(P) \le d(R_p)$ would have a weak limit $\mu$ as in the proof of Theorem~\ref{thm:symmetric} carrying a strict mass excess for $T$ from its adjacent modes, contradicting the lower bound proved there (using the reciprocal field for $p \ge 4$ and $\mu T = 0$ by symmetry for $p = 3$). Since $\delta$ is bounded, a uniform margin follows.

\paragraph{Hausdorff distance bounds.}
In the $\theta$-coordinate, the unit-perimeter boundary velocity corresponding to $G_0$ and that of $R_p$ share the integrable factor $C_p(1 - e^{ip\theta})^{-\beta}/(2\pi)$, while their remaining factors $i z A(z)^{-\beta}(\arg z)'/Q$ and $i e^{i\theta}$ differ uniformly by $O_p(e)$. Combining this with $\|F - \alpha G_0\|_1 \lesssim_p \sqrt{\delta}$, matching basepoints at $\theta = 0$ by translation, and integrating around $\T$ yields $H_p(P) \lesssim_p e + \sqrt{\delta}$. Together with the spectral lower bound above, this proves (i) for small $e$; for larger $e$, increasing $K_p$ suffices since $d(P) \ge 0$ and $H_p(P) \le 1$. Strict monotonicity \eqref{eq:reg-monotonicity} immediately implies the neighborhood gap $d(P) \ge d(R_{p+1}) + \gamma_p$.

At $p = 3$, this improves to $H_3(P) \lesssim e + \delta$. For any real $g \in L^1(\T)$, define the Schwarz integral $\mathcal{S} g(z) = \int_\T (s + z)(s - z)^{-1} g(s)\,dm(s)$. Since $|S(r\eta)|^{-\beta} \lesssim_p (1 - r)^{-\beta}$ for separated prevertices, the primitive of $S^{-\beta}\mathcal{S} g$ starting from $0$ is uniformly bounded on $\D$ by $C_p N(g)$, because $\int_0^1 (1 - r)^{-\beta}(1 - r + |s - \eta|)^{-1}\,dr \lesssim_p |s - \eta|^{-\beta}$. By the rotated single-prevertex bound in \eqref{eq:weighted-tilt-coercivity}, $N(|k|^2) \lesssim \delta$, while $N(|h|^2 - 1/Q) \lesssim e + \delta$ because $m_e/q - 1/Q$ has logarithmic coefficients of size $O(e)$ and $N(|h|^2 - m_e/q) \lesssim \delta$. Writing $h(0) = |h(0)|e^{i\theta_0}$ (so $\alpha = e^{2i\theta_0}$), the function $e^{-2i\theta_0} h^2 \in H^1(\D)$ is real at $0$ and therefore equals the Schwarz integral of its boundary real part; since $\Im(e^{-i\theta_0} h) = \Im(e^{-i\theta_0} k)$, we have
\[
  h^2 - \frac{\alpha}{Q}
  = \alpha\,\mathcal{S}\!\left(|h|^2 - \frac{1}{Q} - 2\bigl[\Im(e^{-i\theta_0} k)\bigr]^2\right).
\]
As the boundary datum has $N$-norm $O(e + \delta)$, the primitive of $F - \alpha G_0 = C_3 S^{-\beta}(h^2 - \alpha/Q)$ has supremum norm $O(e + \delta)$. Comparing $G_0$ with $R_3$ adds $O(e)$, giving $H_3(P) \lesssim e + \delta$ and hence $d(P) - d(R_3) \ge c_3 H_3(P) - K_3 e$ in the stated three-vertex regimes.

Finally, for $R_3$ the series $B_\mu$ in \eqref{eq:spectralmom} equals $D_3$, and since its boundary density is nonconstant, $d(R_3) > D_3$. The asymptotic estimate of Theorem~\ref{thm:symmetric} therefore stays strictly above $D_3$ in a neighborhood of the regular triple, and enlarging $K_3$ covers the remaining triples.
\end{proof}

Note that Theorem~\ref{thm:intro-blaschke} is the direct union of Corollary~\ref{cor:arc-order}, Theorem~\ref{thm:symmetric}, and Corollary~\ref{cor:sharp-gap-stability}(i).

\section{One reflex vertex with six or more convex vertices}\label{sec:one-large}

We now establish the strict comparison $J(P) < J(R_{p+1})$ for all simple polygons with $p \ge 6$ strictly convex vertices and at most one reflex vertex. To reach the fourth-order precision needed to separate $P$ from $R_{p+1}$ uniformly down to $p = 6$, we combine three ingredients: a circle chain rule after removing the first harmonics, a cross-gap chord estimate, and a refined gap trial profile augmented by an interior plateau that captures the full within-gap phase contribution.

\subsection{First-harmonic-removed chain rule and cross-gap chord bounds}

\begin{lem}[Chain rule after removing the first harmonics]\label{lem:markov-square-chain}
Let $P_1$ be the orthogonal projection onto $\operatorname{span}\{z, \bar{z}\}$ on $\T$. Suppose $h \in H^{1/2}(\T)$ is real with $\bar{h} = 0$ and $P_1 h = 0$, and put $h_2 = \int_\T h^2\,dm$. For any real $C^1$ function $F$ with bounded derivative on the closed interval spanned by the essential range of $h$, let $G' = |F'|^2$. Then
\[
  \E(F(h)) - \Var_m\bigl(F(h)\bigr) + \frac{1}{2}\|P_1 F(h)\|_2^2
  \le \bigl\langle (\Lambda - I)h, G(h) \bigr\rangle.
\]
In particular, if $1 + 2h \ge r > 0$ almost everywhere, $A = \sqrt{1 + 2h}$, and $\psi(h) = h - \log A$, then
\begin{equation}\label{eq:phase8chain}
  \E(A) - \Var_m(A) + |\widehat{A}(1)|^2
  \le \E(h) - h_2 - \bigl\langle (\Lambda - I)h, \psi(h) \bigr\rangle,
\end{equation}
with equality if and only if $h = 0$ almost everywhere.
\end{lem}

\begin{proof}
Let $P_0$ be the projection onto constants. Adding $P_0 + \frac{1}{2}P_1$ to $\Lambda - I$ does not change its action on $h$, yet makes the bilinear difference kernel of $\Lambda - I + P_0 + \frac{1}{2}P_1$ pointwise nonnegative on $\T \times \T$:
\[
  \frac{1}{1 - \cos t} - 1 - \cos t = \frac{\cos^2 t}{1 - \cos t} \ge 0.
\]
Integrating the one-dimensional Cauchy--Schwarz inequality $|F(u) - F(v)|^2 \le (u - v)(G(u) - G(v))$ against this nonnegative kernel gives
\[
  \bigl\langle (\Lambda - I + P_0 + \tfrac{1}{2}P_1) F(h), F(h) \bigr\rangle
  \le \bigl\langle (\Lambda - I + P_0 + \tfrac{1}{2}P_1) h, G(h) \bigr\rangle.
\]
Because $P_0 h = P_1 h = 0$, the right-hand side reduces to $\langle(\Lambda - I)h, G(h)\rangle$, proving the first claim. Taking $F(x) = \sqrt{1 + 2x}$ gives $G(x) = \frac{1}{2}\log(1 + 2x) = \log A$ and yields \eqref{eq:phase8chain}. Since $F'$ is strictly decreasing, Cauchy--Schwarz is strict on every nondegenerate interval; as the difference kernel is strictly positive for almost every angular separation $t$, equality holds precisely when $h$ is constant almost everywhere, hence $h = 0$.
\end{proof}

To bound the phase interaction between distinct gaps, we first compare the distance to $p$ arbitrary points on $\T$ against equally spaced centers.

\begin{lem}[Comparison with equally spaced centers]\label{lem:capped-distance}
In arc-length balanced coordinates, let $D(t)$ be the distance in turns from $t \in \R/\Z$ to the nearest of $p$ prescribed points. Put $y = 1/p$ and $D_0(t) = \dist(t, y\Z)$. For every nondecreasing Lipschitz function $g$ on $[0, 1/2]$,
\begin{equation}\label{eq:equal-distance-test}
  \int_\T g(D)\,d\mu \ge \int_\T g(D_0)\,dm - \sqrt{B_\mu\,\mathcal{N}\bigl(g(D_0)\bigr)}.
\end{equation}
For the uniform measure ($B_\mu = 0$) and strictly increasing $g$, equality holds if and only if the $p$ points are equally spaced. For $p \ge 6$ and $B_\mu \le 3y^3$, set $D_c = \min(D, y/2)$ and $\gamma = \sqrt{B_\mu[\zeta(3) - y^2\zeta(5)]/(2\pi^4 y^3)} < 0.137$. Then the normalized chord functions satisfy
\begin{equation}\label{eq:relative-chord}
  \mu L_a(D_c) \ge 1 - 12\gamma y, \qquad
  L_a(t) = \frac{1 - \cos(2\pi a t/y)}{1 - \sin(\pi a)/(\pi a)} \quad (0 < a \le 1), \qquad
  L_0(t) = \frac{12 t^2}{y^2}.
\end{equation}
\end{lem}

\begin{proof}
Set $F(r) = m\{D \le r\}$. Increasing $r$ by $h > 0$ enlarges the union of the $p$ closed arcs by total length at most $2ph$, so $X \coloneqq F(D)/(2p)$ is $1$-Lipschitz and satisfies $X \le D$. Because the level sets of $D$ are finite, $F$ is continuous and $X$ has the same uniform law on $[0, y/2]$ as $D_0$.

On the subspace of even $y$-periodic functions, the rearrangement map $T : g(D_0) \mapsto g(X)$ preserves $L^2$ inner products and contracts $\E$: by reflection and periodicity, the harmonic extension of $g(D_0(t))$ to the half-cylinder has the form $U(s, D_0(t))$, and composing with $X$ gives an extension $U(s, X(t))$ whose Dirichlet energy is no larger because $X$ and $D_0$ are equidistributed while $|X'| \le 1 = |D_0'|$ almost everywhere. For any mean-zero $f$ in this subspace, setting $h = \Lambda^{-1} f$ (which lies in the same subspace) gives by duality
\[
  \|Tf\|_{\dot{H}^{-1/2}}^2 \ge \sup_{\widetilde{h}}\bigl\{2\langle Tf, T\widetilde{h}\rangle - \E(T\widetilde{h})\bigr\} \ge \|f\|_{\dot{H}^{-1/2}}^2.
\]
Thus $\mathcal{N}(Tf) \le \mathcal{N}(f)$. Applying \eqref{eq:duality} to $g(X) \le g(D)$ proves \eqref{eq:equal-distance-test}. If any gap exceeds $y$, then $D > y/2 \ge X$ on an open interval near its midpoint, giving strict inequality for $m$.

Each function $L_a(D_0)$ has mean $1$. For $0 < a < 1$, set $u = a^2$ and $P(u) = \prod_{j=2}^\infty (1 - u/j^2)^{-1}$. Integrating on one period and applying Euler's sine product gives, for every integer $n \neq 0$,
\[
  \frac{1}{|\widehat{L_a(D_0)}(pn)|}
  = \frac{(n^2 - u)(\pi a/\sin(\pi a) - 1)}{u}
  = \frac{n^2 - u}{1 - u}\left(1 + \frac{P(u) - 1}{u}\right).
\]
Both factors are nondecreasing in $u \in (0,1)$ (the second via its positive Taylor series), starting at $n^2\zeta(2)$ as $u \downarrow 0$ and extending by continuity to $a = 1$. Hence $\mathcal{N}(L_a(D_0)) \le \mathcal{N}(L_0(D_0)) = 2p[\zeta(3) - y^2\zeta(5)]/\zeta(2)^2$. Applying \eqref{eq:equal-distance-test} to $g(t) = L_a(\min(t, y/2))$ establishes \eqref{eq:relative-chord}.
\end{proof}

\begin{lem}[Cross-gap chord bound]\label{lem:cross-chord}
Suppose $p \ge 6$ and $P$ has at most one reflex vertex. Set $y = 1/p$ and assume $d(P) \le 3y^3$ and $u = \max_i s_i \le U \le 1/4$. Let $\mathcal{R}$ denote the contribution of all ordered pairs of distinct gaps $(i, k)$ ($i \neq k$) to the phase integral in \eqref{eq:phase-product} for $f|_{I_i} = s_i^3\Phi$, and set $v_\Phi = \int_0^1 (t - 1/2)^2\Phi(t)\,dt$. Then
\begin{align*}
  \frac{4Z\mathcal{R}}{y^4}
  \ge 4Z\pi^2 y\bigg[
    &v_\Phi\Bigl(1 - \frac{\pi^2 U^2}{12}\Bigr)(1 + 0.34 y)^{-5} \\
    &+ (1 - 2\pi^2 v_\Phi y^2)\Bigl(1 - \frac{\pi^2 y^2}{20}\Bigr)
      \Bigl(\frac{1}{12} - 0.137 y\Bigr)\frac{(1 - U)^4}{(1 - y)^3}\bigg].
\end{align*}
In particular, if $u \le 3y/2$, then $4Z\mathcal{R}/y^4 > 2.83 y$.
\end{lem}

\begin{proof}
Write $S_k = \mu(I_k)$, and let $A_i$ and $B_k$ be the average boundary phases $\int u$ on $I_i$ and $I_k$ with respect to the probability measures $\Phi(v)\,dv$ and $d\mu/S_k$, respectively. Since $|w - z|^2 \le 4$ on $\T$,
\[
  \mathcal{R} \ge \frac{1}{4}\sum_{i \neq k} s_i^4 S_k\bigl(2 - 2\Re(A_i \overline{B_k})\bigr)
  \ge \frac{1}{2}\sum_{i \neq k} s_i^4 S_k(1 - |A_i|\,|B_k|).
\]
Because each gap contains at most one reflex vertex, any two points $z, w$ in the same gap satisfy $\Re(u(z)\overline{u(w)}) \le \Re(z\bar{w})$ when $s_i \le 1/4$ (as in Lemma~\ref{lem:single-gap-loss}); averaging this pairwise comparison bounds $|A_i|$ and $|B_k|$ above by the moduli $a_i$ and $b_k$ of the corresponding averages of the circle coordinate $z$. By symmetry of $\Phi$, $a_i = \int_0^1 \cos(2\pi s_i(v - 1/2))\Phi(v)\,dv$, whence
\[
  1 - a_i \ge 2\pi^2 v_\Phi s_i^2\Bigl(1 - \frac{\pi^2 U^2}{12}\Bigr), \qquad
  a_i \ge 1 - 2\pi^2 v_\Phi s_i^2.
\]
Termwise integration of \eqref{eq:phi-profile} gives $v_\Phi = 1/12 - \zeta(5)/(2\pi^2 Z) \in (0.0395, 0.040)$.

Split $1 - a_i b_k = (1 - a_i) + a_i(1 - b_k)$. Since $s_i \le 1/4$, the cosine relative to the midpoint of $I_i$ is at least $\cos(\pi/4) > 2/3$, so balance $\mu z = 0$ forces $S_i \le 3/5$. On $[0, 3/5]$ we have $(1 - S)^{-1/5} \le 1 + 0.34 S$ (by convexity and $1.204^5 > 5/2$), so H\"older's inequality yields
\[
  \sum_{i=1}^p (1 - S_i)s_i^6
  \ge \biggl(\sum_{i=1}^p (1 - S_i)^{-1/5}\biggr)^{-5}
  \ge y^5(1 + 0.34 y)^{-5}.
\]
For the second part, $q \mapsto q^4(1 - 2\pi^2 v_\Phi q^2)$ is increasing and convex on $[0, 1/4]$; applying Jensen's inequality and $(1 - U)/(p - 1) \le y$ gives
\[
  \sum_{i \neq k} a_i s_i^4 \ge \frac{(1 - U)^4}{(p - 1)^3}(1 - 2\pi^2 v_\Phi y^2).
\]
Let $\alpha_k$ be the argument (in turns) of the $\mu$-average of $z$ on $I_k$, and let $D(t)$ be the distance from $t$ to the nearest $\alpha_k$. Monotonicity of $1 - \cos(2\pi t)$ gives $\sum_k S_k(1 - b_k) \ge \mu(1 - \cos(2\pi D)) \ge \mu(1 - \cos(2\pi D_c))$. Applying \eqref{eq:relative-chord} at $a = y$ and using $1 - \sin u/u \ge (u^2/6)(1 - u^2/20)$ at $u = \pi y$ completes the proof of the displayed inequality.

When $U = 3y/2$, $(1 - U)^4/(1 - y)^3$ decreases with $y$ and attains its minimum $(3/4)^4/(5/6)^3 > 0.546$ at $y = 1/6$, where every other factor in the bracket is likewise minimal; direct arithmetic gives
\[
  47.4\bigl[0.0395(0.948)(0.758) + (0.978)(0.977)(0.0605)(0.546)\bigr] > 2.83. \qedhere
\]
\end{proof}

\subsection{The combined trial profile and pointwise strip barriers}
We now augment the gap profile $f$ by an interior trapezoidal plateau $g$ that harvests the within-gap phase integral uniformly across each gap.

\begin{lem}[Combined profile for dispersion and energy estimates]\label{lem:kernel-dispersion}
In arc-length balanced coordinates, suppose $P$ has at most one reflex vertex. Set $y = 1/p$, $T_j = \sum_{i=1}^p s_i^j$, $M = T_4$, $u = \max_i s_i \le U < 1/2$, $w = u/y$, $\delta = (M - y^3)/y^4$, $\epsilon = 1/400$, and $\lambda = 1 - \epsilon$. For any $0 < k \le \kappa(U)$ (with $\kappa$ from Lemma~\ref{lem:single-gap-loss}), put
\[
  D_k(s) = c + 0.750 s + 0.026 k s^2 + (0.84 k - 0.70)s^3 + 0.77 k^2 s^4.
\]
Assign $X = s_i$ probability $s_i^4/M$, and set $m = \mathbb{E} X = T_5/M$, $\bar{D}_k = \mathbb{E} D_k(X)$, and $A_k = (1 + k\lambda m)^2/\bar{D}_k$. Then
\begin{equation}\label{eq:direct-dispersion}
  d(P) \ge M A_k + \frac{2 A_k}{1 + k\lambda m}\,\mathcal{R},
\end{equation}
where $\mathcal{R}$ is the distinct-gap phase contribution for $f|_{I_i} = s_i^3\Phi$.
\end{lem}

\begin{proof}
On gap $I_i$ with unit coordinate $v \in [0,1]$, put $h_i = s_i^3$, $f_i = h_i\Phi(v)$, and $g_i = h_i s_i B(v)$ extended by zero outside $I_i$, where $B(v) = \min\{1, v/\epsilon, (1 - v)/\epsilon\}$. Define $f = \sum_i f_i$, $g = \sum_i g_i$, and $G = (I + \Lambda)f + kg$. By Lemma~\ref{lem:single-gap-loss}, at every target point in $I_i$ the within-gap phase integral of $f_i$ is at least $\kappa(U) s_i h_i \ge k g_i$. Combining this with \eqref{eq:phase-product} and \eqref{eq:duality} yields
\begin{equation}\label{eq:profile-phase-dual}
  M(1 + k\lambda m) + \mathcal{R} \le \bar{G} - \mu G \le \sqrt{d(P)\,\mathcal{N}(G)}.
\end{equation}

We bound each term in the expansion of $\mathcal{N}(G)$. For $i \neq j$, the nonnegative supports of $f_i$ and $g_j$ are disjoint, so the energy difference kernel gives
\[
  \E(f_i, g_j) = -2\int_{I_i}\!\int_{I_j} \frac{f_i(z) g_j(w)}{|z - w|^2}\,dm(w)\,dm(z) \le 0.
\]
For a diagonal term $i = j$, extend $\Phi$ by zero from $(0,1)$ to $\R$ and set
\[
  J(v) = \frac{1}{2\pi^2}\,\pv\!\int_\R \frac{\Phi(v) - \Phi(w)}{(v - w)^2}\,dw.
\]
Adding nonnegative translated copies of $\Phi$ lowers $J$ on $(0,1)$: periodizing with period $1/s_i$ gives $\E(f_i, g_i) \le h_i^2 s_i \int_0^1 B J\,dv$ by scaling, while periodizing with period $1$ gives $J(v) \ge \Lambda_\T \Phi(v) = -Z^{-1}\sum_{n \ge 1} \cos(2\pi n v)/n^2 \ge -\zeta(2)/Z$. The internal principal values cancel in $\int_0^1 J(v)\,dv$, and since $\Phi(v) \le Z^{-1}\sum_{n \ge 1} (1 - \cos(2\pi n v))/n^2 = \pi^2 v(1 - v)/Z$, we obtain
\[
  \int_0^1 B J\,dv \le \frac{1}{2\pi^2}\int_0^1 \frac{\Phi(v)}{v(1 - v)}\,dv + \frac{\epsilon\zeta(2)}{Z}
  \le c + \frac{\epsilon\zeta(2)}{Z} < 0.420,
\]
whence $\E(f, g) \le 0.420\sum_i h_i^2 s_i$. By diagonal Cauchy--Schwarz, $\E(f) \le \sqrt{(b/2)(1 + e/2)}\sum_i h_i^2 < 0.714\sum_i h_i^2$, while \eqref{eq:phi-neumann-lower} gives $\E(f) \ge 2c^2\zeta(5)\sum_i h_i^2$. Comparing similarly with the zero extension on $\R$ gives $\E(g) \le \sum_i h_i^2 s_i^2 \E_\R(B)$ by scale invariance. On $\R$, integrating by parts expresses $\E_\R(B)$ via the kernel $-(2\pi^2)^{-1}\log|v - w|$ acting on $B'$ (the difference of uniform probability densities on $[0,\epsilon]$ and $[1-\epsilon, 1]$), so Jensen's inequality gives
\[
  \pi^2 \E_\R(B) = \int_0^1\!\int_0^1 \log\frac{1 - \epsilon + \epsilon(v - w)}{\epsilon|v - w|}\,dv\,dw
  \le \log\frac{1 - \epsilon}{\epsilon} + \frac{3}{2} < 0.77\pi^2,
\]
using $\int_0^1\int_0^1 (-\log|v - w|)\,dv\,dw = 3/2$. In angular coordinates, $f_i'' = c h_i q(v)/s_i^2$ with $q(v) = -2\log(2\sin\pi v)$, so integrating by parts gives $\langle f', g'\rangle \le 0.031 c\sum_i h_i^2$ because $\int_0^1 q B\,dv \ge 2\epsilon\log(4\epsilon) - 3\epsilon > -0.031$. Finally, $\Phi'$ vanishes at $0$ and $1/2$ and is concave on $(0, 1/2)$, so both $\Phi$ and $B$ increase toward the midpoint of each gap, while their gap averages $h_i$ and $\lambda s_i h_i$ increase together with $s_i$; hence $\Cov_m(f, g) \ge 0$.

Expanding $\mathcal{N}(G)$ and dropping $-2k\Cov_m(f, g) - \|f + kg\|_{\dot{H}^{-1/2}}^2 \le 0$, we find
\[
  \mathcal{N}(G) \le \|f\|_{\dot{H}^{3/2}}^2 + 2\|f'\|_2^2 + 2k\langle f', g'\rangle + 2k\E(f, g) + k^2\E(g) - 2\Var_m(f).
\]
Applying \eqref{eq:phi-energy-bounds} and the bounds just established yields $\mathcal{N}(G) \le M\bar{D}_k$. Substituting into \eqref{eq:profile-phase-dual} and expanding the square gives \eqref{eq:direct-dispersion}.
\end{proof}

We record three consequences of Lemma~\ref{lem:kernel-dispersion} that will be used both here and in Section~\ref{sec:low-reflex}:
\begin{enumerate}[label=\textup{(\alph*)},leftmargin=*]
\item \emph{Tangent lower bound on $A_k$.} Whenever $u/y \le 7/4$ and $ku \le 7/20$, using $X^2 \le uX$ and $\Var X \le u\mathbb{E}[(X - m)^2/X] = u(m^2\mathbb{E}(1/X) - m)$ gives $2\mathbb{E} X^2/m^2 \le 1 + u T_3/T_4 \le 1 + u/y \le 11/4$ (since $T_4 \ge y T_3$). Put $a = k\lambda$, $b_0 = 0.750/c \in (1.803, 1.805)$, and $\ell = 1 + (2a - b_0)m \in (0, 1.1)$. Since $D_k(X) \le c + 0.750 X + [0.026 + 0.14 ku + 0.77(ku)^2]X^2$ and the bracket is less than $0.17$, we have $\bar{D}_k/c \le 1 + b_0 m + (0.17/c)(11/8)m^2$, while $(1 + am)^2 - \ell(1 + b_0 m) = (b_0 - a)^2 m^2$. Because $(b_0 - a)^2 > 0.64 > 1.1(0.17)(2.406)(11/8)$, it follows that $A_k \ge \ell/c$; in particular,
\begin{equation}\label{eq:Ak-tangent}
  A_{0.95} \ge 2.4 + 0.19 m, \qquad A_1 \ge 2Z + 0.4 m.
\end{equation}
\item \emph{Independent increasing amplitudes.} Replacing $s_i^3$ by any nonnegative amplitudes $h_i$ increasing with $s_i$ preserves $\Cov_m(f, g) \ge 0$ and gives $\mathcal{N}(G) \le \sum_i h_i^2 D_k(s_i)/s_i^2$ with leading mean $\sum_i h_i s_i(1 + k\lambda s_i)$. Taking $h_i = s_i^3(1 + k\lambda s_i)/D_k(s_i)$ (which increases in $s_i$ since $3D_k - sD_k' = 3c + 1.50 s + 0.026 k s^2 - 0.77 k^2 s^4 > 0$) yields $d(P) \ge \sum_i s_i^4(1 + k\lambda s_i)^2/D_k(s_i) + 2\min_i(h_i/s_i^3)\mathcal{R}$.
\item \emph{Moment interpolation.} Put $r_i = s_i/y$ and $q_j(v) = \sum_{l=0}^{j-2}(j - 1 - l)v^l$. Whenever $w = u/y \le w_*$, setting $a_j = y q_j(w_*)/q_4(w_*)$ yields, for every integer $j \ge 4$,
\begin{equation}\label{eq:moment-interp}
  (1 + y\delta)^{(j-1)/3} \le \frac{T_j}{y^{j-1}} \le 1 + a_j\delta.
\end{equation}
The lower bound follows from power means with weights $s_i$, and the upper bound follows by summing $r_i^j - 1 - j(r_i - 1) = (r_i - 1)^2 q_j(r_i)$ and noting that $q_j/q_4$ is increasing on $[0,\infty)$ (the ratios of its first three coefficients increase and its higher coefficients are nonnegative).
\end{enumerate}

\begin{lem}[Initial dispersion and profile bounds]\label{lem:common-profile-bounds}
Assume $P$ has at most one reflex vertex, and retain the notation of Lemma~\ref{lem:kernel-dispersion}. For $p \ge 6$, the hypothesis $d(P) \le d(R_{p+1})$ implies
\begin{equation}\label{eq:initial-disp}
  u < \frac{1}{4}, \qquad \delta = \frac{T_4 - y^3}{y^4} < 0.85, \qquad \delta < 0.84 \quad (p \ge 7).
\end{equation}
For any $p$ with $u \le 1/4$, take $k = 1$ in $G = (I + \Lambda)f + g$, and set $\widetilde{h} = (\bar{G} - G)/c$, $h = \widetilde{h} - P_1\widetilde{h}$, $h_2 = \int_\T h^2\,dm$, $E = \E(h)$, and $\epsilon = 1/400$. Then
\begin{gather*}
  c\int_\T h\,d\mu \ge M + (1 - \epsilon)T_5 + \mathcal{R}, \\
  E \le 2Z M + 4W T_5 + K(u)T_6, \qquad
  h_2 \le c^{-2}\bigl[\sqrt{b/2} + (1.17 + u)u\bigr]^2 T_5,
\end{gather*}
where $K(u) = 2.406^2(0.714 + 0.84 u + 0.77 u^2) + 4.812(0.031)$. Moreover, $\|P_1\widetilde{h}\|_\infty \le (4|\widehat{f}(1)| + 2|\widehat{g}(1)|)/c$.
\end{lem}

\begin{proof}
By Corollary~\ref{cor:regular-uniform-nine} (for $p \ge 7$) and $216 d_7 < 216(0.013273) < 2.867$ (for $p = 6$), we have $d/y^3 \le d_* \coloneqq 2Z + 2.78 y$. Lemma~\ref{lem:gap} gives $M \le (c + bu)d$, whence $\delta < 1.16 + 2.16 w$. On the other hand, expanding $r_i^4$ about $1$ gives
\begin{equation}\label{eq:kernel-spread}
  \delta = \sum_{i=1}^p (r_i - 1)^2(r_i^2 + 2r_i + 3)
  \ge (w - 1)^2(w^2 + 2w + 3) \ge 6(w - 1)^2, \qquad r_i = \frac{s_i}{y},
\end{equation}
which forces $w < 7/4$. If $u \ge 1/4$, then necessarily $p = 6$, $u < 7/24 < 1/3$, and $\delta > 2$ by \eqref{eq:kernel-spread}; taking $k = 0.95$ (so $ku < 7/20$) in \eqref{eq:direct-dispersion} and \eqref{eq:Ak-tangent} gives $d/y^3 \ge 2.4(1 + y\delta) > 3.2$, contradicting $d/y^3 < 2.87$. Thus $u < 1/4$.

Now take $U = \min(7y/4, 1/4)$ and $k = 1$. The ratio $(1 - U)^4/(1 - y)^3$ attains its minimum at $y = 1/7$, where it exceeds $1/2$, so Lemma~\ref{lem:cross-chord} gives $\mathcal{R}/y^3 > y^2/2$ (since $\pi^2[0.039(0.94)(0.73) + (0.95)(0.97)(0.059)/2] > 1/2$). Since $m \le 1/4$, \eqref{eq:Ak-tangent} gives $A_1 \ge 2Z + 0.4 m \ge 2(1 + \lambda m)$. Substituting into \eqref{eq:direct-dispersion} and using $T_5 \ge y M$ yields
\[
  \frac{d}{y^3} \ge (1 + y\delta)(2Z + 0.4 y) + 2y^2,
\]
so $R \coloneqq (d - 2Z y^3)/y^4 \ge 2Z\delta + 0.4 + (2 + 0.4\delta)y$. At $\delta = 0.85$ this exceeds $2.78$ for $p = 6$; at $\delta = 0.84$ it exceeds $2.741$ for $p = 7$, $2.60$ for $p = 8, 9$, and $2.34$ for $p \ge 10$, contradicting Corollary~\ref{cor:regular-uniform-nine}. This proves \eqref{eq:initial-disp}.

For the profile bounds, the lower bound on $\int h\,d\mu$ follows directly from \eqref{eq:profile-phase-dual} and balance ($\mu P_1\widetilde{h} = 0$). Because removing $P_1$ does not change energy minus variance, expanding $G = (I + \Lambda)f + g$ yields
\begin{equation}\label{eq:kernel-core}
  c^2(E - h_2) = \|f\|_{\dot{H}^{3/2}}^2 + \|f'\|_2^2 - \E(f)
    + 2\langle f', g'\rangle + \E(g) - \Var_m(f + g).
\end{equation}
Similarly, $\E(G) \le \|f\|_{\dot{H}^{3/2}}^2 + 2\|f'\|_2^2 + 2\langle f', g'\rangle + \E(f) + 2\E(f, g) + \E(g)$, and substituting the estimates from the proof of Lemma~\ref{lem:kernel-dispersion} gives the bound for $E$. Centering $f + g$ and using $\|f + g\|_2 \le (1.17 + u)u\sqrt{T_5}$ bounds $h_2$. Finally, because $I + \Lambda$ multiplies modes $\pm 1$ by $2$, we have $\|P_1\widetilde{h}\|_\infty \le (4|\widehat{f}(1)| + 2|\widehat{g}(1)|)/c$.
\end{proof}

\begin{lem}[Pointwise barrier bounds for $\Lambda f$]\label{lem:strip-phase-pointwise}
Let $f|_{I_i} = s_i^3\Phi(v)$ ($v \in [0,1]$), and put $m = \min_i s_i$ and $u = \max_i s_i$. On each gap $I_i$,
\[
  \frac{\pi^2 s_i^2}{Z}v(1 - v) - \frac{\zeta(2)}{Z}u^2
  \le \Lambda f
  \le \frac{\pi^2 s_i^2}{Z}v(1 - v) - \frac{\zeta(2)}{Z}m^2.
\]
For equal gaps both bounds coincide, with $-\Lambda f = [\zeta(2)/Z]u^2$ at every prevertex. In general, $-\Lambda f \le [\zeta(2)/Z]u^2 < 1.37 u^2$, and in terms of $t = v - 1/2$,
\[
  \frac{\Lambda f}{s_i^2} \le C_\Phi(1 - 4t^2) - \frac{8t^2}{\pi^2}, \qquad
  C_\Phi \coloneqq \frac{2}{Z}\left(3\log 2 + \frac{17\pi^2}{72} - \frac{35}{9}\right) < 0.87 < \frac{7}{8}.
\]
\end{lem}

\begin{proof}
Work on the half-cylinder of angular circumference $2\pi$ with vertical depth $t > 0$. On a gap of length $s$ (in turns), consider the Neumann harmonic extension and its endpoint trace:
\[
  N_s(v, t) = s^3\left[1 - \frac{1}{Z}\sum_{n=1}^\infty \frac{e^{-nt/s}\cos(2\pi nv)}{n^3}\right], \qquad
  A_s(t) = N_s(0, t).
\]
Setting $V = N_s - A_s + A_u$ on each half-strip glues continuously across every vertical ray (where its value is $A_u(t)$ and its horizontal normal derivatives vanish). Moreover, $\Delta V = A_u'' - A_s'' \le 0$ because $-Z A_s''(t) = s[-\log(1 - e^{-t/s})]$ is a product of positive functions increasing in $s$. Thus $V$ is globally superharmonic, agrees with $f$ at $t = 0$, and tends to $u^3 \ge \bar{f}$ as $t \to \infty$, so $V$ lies above the harmonic extension of $f$ by the maximum principle. Replacing $u$ by $m$ produces a subharmonic lower barrier. Taking inward normal derivatives at $t = 0$ (where the harmonic extension has derivative $-\Lambda f$) and using $\sum_{n \ge 1}(1 - \cos(2\pi nv))/n^2 = \pi^2 v(1 - v)$ proves the two-sided inequality.

For the parabolic upper bound, extend $\phi(t) = \Phi(1/2 + t)$ by zero outside $[-1/2, 1/2]$. Dropping all other gaps and translated copies from the periodized difference kernel yields $\Lambda f/s_i^2 \le J(t)$, where $J(t) = (2\pi^2)^{-1}\,\pv\int_\R (\phi(t) - \phi(v))/(t - v)^2\,dv$. Restoring the unit-periodic copies gives
\[
  J(t) = \frac{\pi^2}{Z}\Bigl(\frac{1}{12} - t^2\Bigr)
    + \frac{1}{2\pi^2}\sum_{k \neq 0} \int_{-1/2}^{1/2} \frac{\phi(v)}{(k + v - t)^2}\,dv.
\]
Pairing $k + v$ with $-k - v$ expresses the sum as a power series in $t^2$ with nonnegative coefficients, so $J$ is convex as a function of $t^2$. At the endpoint $t = 1/2$, Jensen's inequality (since $\Phi$ is a symmetric probability density on $[0,1]$) gives $J(1/2) = -(2\pi^2)^{-1}\int_0^1 \Phi(v)v^{-2}\,dv \le -2/\pi^2$. At the center $t = 0$, integrating by parts once using $\phi'(0) = 0$ and $\phi''(v) = -(4\pi^2/Z)\log(2\cos\pi v)$ yields
\[
  J(0) = -\frac{1}{\pi^2}\int_0^{1/2} \frac{\phi'(v)}{v}\,dv
       = \frac{2}{Z}\int_0^1 (-\log x)\log\Bigl(2\cos\frac{\pi x}{2}\Bigr)\,dx.
\]
By Euler's cosine product and $1 - \theta \le e^{-\theta}$ on the factors past $1 - x^2$, we have $\cos(\pi x/2) \le (1 - x^2)\exp(-(\pi^2/8 - 1)x^2)$. Using $\int_0^1 (-\log x)\log(1 - x^2)\,dx = -\sum_{n \ge 1}[n(2n+1)^2]^{-1} = 2\log 2 + \pi^2/4 - 4$ gives $J(0) \le C_\Phi < 0.87$ (since $\log 2 < 0.6932$, $\pi^2 < 9.87$, and $Z > 1.20$). The chord joining $(0, C_\Phi)$ to $(1/4, -2/\pi^2)$ in the variable $t^2$ yields the stated parabolic bound.
\end{proof}

\subsection{Comparison for six or more convex vertices}
We now assemble the tools of Subsections~6.1--6.2 to rule out any polygon with $p \ge 6$ strictly convex vertices, at most one reflex vertex, and $d(P) \le d(R_{p+1})$. Retain the profile $f, g, h, B, \epsilon, K$ of Lemma~\ref{lem:common-profile-bounds}, and set $y = 1/p$, $M = T_4$, $w = u/y$, $\delta = (M - y^3)/y^4$, $s = \sqrt{\delta}$, $D = 2\int h\,d\mu - \E(h)$, and $h_2 = \int h^2\,dm$.

From \eqref{eq:kernel-core}, Lemma~\ref{lem:common-profile-bounds}, and the within-gap variance bound in \eqref{eq:phi-energy-bounds}, we have for $u \le U \le 1/4$:
\begin{align}
  D + h_2 &\ge 2Z M + \bigl[4Z(1 - \epsilon) - 2W\bigr]T_5 + B_6 T_6 + 2\zeta(6)T_7 - B_8 T_8 + 4Z\mathcal{R}, \label{eq:phase8core} \\
  \E(h) &\le 2Z M + 4W T_5 + K(U)T_6, \notag
\end{align}
where $B_6 \coloneqq 2.072 - 4.812(0.031)$ and $B_8 \coloneqq 2.406^2(0.77)$. Note that $W = \pi^4/90 < 1.0824$. Note also that since the $p - 1$ deviations $r_i - 1$ other than the maximum sum to $1 - w$, \eqref{eq:kernel-spread} strengthens to $\delta \ge (w - 1)^2[w^2 + 2w + 3 + 3/(p - 1)]$.

\paragraph{Applying the chain rule and second-harmonic compensation.}
Suppose $-a \le h \le b_+$ with $r = 1 - 2a > 0$, and put $A = \sqrt{1 + 2h}$ and $\psi(h) = h - \log A$, so that \eqref{eq:phase8chain} applies. Using the barrier bound of Lemma~\ref{lem:strip-phase-pointwise}, we have $(\Lambda - I)h + c^{-1}\Lambda g = c^{-1}(f'' + f + g - \bar{G}) \ge -b_0$, where $b_0 = 1.387 u + \bar{G}/c$. For any $t \ge 0$, let $R_t(x) = \psi(x) - \psi(t) - \psi'(t)(x - t) \ge 0$ be the nonnegative convex tangent remainder of $\psi$ at $t$. The chord bound for $(1 + 2sx)^{-1}$ on $[-a, \infty)$ gives
\[
  \psi(x) = 2x^2\int_0^1 \frac{s'\,ds'}{1 + 2s' x}
  \le 2x^2\int_0^1 s'\bigl[1 + s'(r^{-1} - 1)\bigr]\,ds' = v_r x^2, \qquad v_r \coloneqq \frac{1 + 2/r}{3}.
\]
Because $(t\psi'(t) - \psi(t))' = 2t/(1 + 2t)^2 \le 2t$ and $\bar{h} = 0$, we have $\overline{R_t(h)} \le v_r h_2 + t^2$. Moreover, $R_t'(x) = 2(x - t)/[(1 + 2t)(1 + 2x)]$, so $\Lip(R_t) \le 2\max(a + t, b_+ - t)/[r(1 + 2t)]$. Combining these estimates with $\E(g) \le 0.77 T_8$ yields
\begin{equation}\label{eq:phase8tangent}
  \bigl\langle (\Lambda - I)h, \psi(h) \bigr\rangle
  \ge \psi'(t)\bigl(\E(h) - h_2\bigr) - b_0(v_r h_2 + t^2)
    - c^{-1}\Lip(R_t)\sqrt{0.77 T_8\,\E(h)}.
\end{equation}
We also note that \eqref{eq:duality} gives $\E(h) \ge [(M + \bar{g})/c]^2/d \ge 2Z y^3$ (since $4Z(1 - \epsilon) > 2.78$).

Now set $j = 1 + d$, $A_1 = \widehat{A}(1)$, and $k = C(A - \bar{A} - P_1 A)$ (with $C$ from Lemma~\ref{lem:circle-gram}), and write $X^2 = \E(k) = \sum_{|n| \ge 2} |\widehat{A}(n)|^2/|n|$. For $0 \le t \le 1$, Lemma~\ref{lem:circle-gram} gives the energy of $U = z(A + itk)$ as $\E(U) = (\bar{A})^2 + \E(A) - (2t - t^2)X^2$. By balance and \eqref{eq:duality} applied to $U - \bar{A} z$ and $zk$,
\[
  |m| \coloneqq |\mu U| \le M_A \coloneqq |A_1| + \sqrt{d\E(h)/r}, \qquad
  |\mu(zk)| \le \sqrt{d}\,X,
\]
using $\E(A) \le \E(h)/r$ and $|A_1| \le 2h_2/(1 + \sqrt{r})^2$ (from $2(A - 1) = 2h - (A - 1)^2$ and $P_1 h = 0$).

To absorb the nonzero mean $m = \mu U$, add the real tangential perturbation $q = 2\Im(\bar{m} z)/j$, forming $W = z(A + i(tk + q)) = U + (\bar{m} z^2 - m)/j$. Writing $\beta = \mu z^2$, we have
\[
  \mu W = \frac{d m + \bar{m}\beta}{j}, \qquad
  \E(W) = \E(U) + \frac{2|m|^2}{j^2} + \frac{4\Re(m A_1)}{j}.
\]
Testing $z$ and $\bar{z}$ gives $|\beta| \le d$, so $\mu q^2 \ge 2(1 - d)|m|^2/j^2$ and $|\mu(kq)| \le 2M_A\sqrt{d}\,X/j$. Substituting into $j\E(W) \ge \mu A^2 + \mu(tk + q)^2 - |\mu W|^2$, using $\overline{A^2} = 1$ and $2j - 2(1 - d) + 4d^2 = 4dj$, and maximizing $j(2t - t^2)X^2 - 4t M_A\sqrt{d}\,X/j$ over $t \in [0,1]$ yields
\[
  d \ge 2\mu h - j\bigl(\E(A) - \Var_m(A)\bigr) + j\Bigl(X - \frac{2M_A\sqrt{d}}{j^2}\Bigr)_+^2
    - \frac{4d M_A^2}{j} - 4|A_1|M_A,
\]
which holds also when $X = 0$. Setting $f_\perp = f - M - P_1 f$ and $v = -f_\perp/c$, Fourier Cauchy--Schwarz gives $X^2 \ge L \coloneqq 2\langle A, v\rangle - \E(v)$. Applying $(X - b)_+^2 \ge (1 - \eta)X^2 - (\eta^{-1} - 1)b^2$ for any $\eta \in (0,1)$ and using $j \ge 1$, we obtain
\begin{equation}\label{eq:phase8schur}
  d \ge 2\mu h - j\bigl(\E(A) - \Var_m(A)\bigr) + (1 - \eta)L
    - \frac{4d M_A^2}{\eta} - 4|\widehat{A}(1)|M_A.
\end{equation}
Writing $F_1 = |\widehat{f}(1)|$ and $G_1 = |\widehat{g}(1)|$, inserting $A = 1 + h - (A - 1)^2/2$ into $L$, and using $f_\perp \ge -M - 2F_1$ together with $\Cov_m(f, g) \ge 0$ gives
\[
  L \ge c^{-2}\bigl[\E(f) - 2F_1^2 + 2\|f_\perp\|_2^2 - 4F_1 G_1\bigr]
    - \frac{4(M + 2F_1)h_2}{c(1 + \sqrt{r})^2}.
\]

\begin{prop}[One reflex vertex with six or more convex vertices]\label{prop:phase64}
If $p \ge 6$ and $P$ has at most one reflex vertex, then $J(P) < J(R_{p+1})$.
\end{prop}

\begin{proof}
Suppose for contradiction that $d = d(P) \le d(R_{p+1})$. Put $Y = 1/6$, $\tau = y/Y \in (0,1]$, $s = \sqrt{\delta}$, and $R = (d - 2Z y^3)/y^4$. By Corollary~\ref{cor:regular-uniform-nine}, $R$ is strictly less than $2.78$, $2.741$, $2.60$, and $2.34$ for $p = 6$, $p = 7$, $p \in \{8, 9\}$, and $p \ge 10$, respectively. Set $d_* = 2.406 + 2.78 y$, so $d < d_* y^3$.

Define $S = 0.925$, $W_* = 1.34$, $U = y W_*$, $\omega = \sqrt{(1 - y)/6}$, and $w_0 = 1 + \omega s$. By \eqref{eq:initial-disp}, $s < S$. Moreover, $u/y \le \min(W_*, w_0)$: putting $l = u/y - 1$, the remaining $p - 1$ normalized gaps have mean $1 - l/(p - 1)$, so Jensen's inequality gives
\[
  \delta \ge (1 + l)^4 + (p - 1)\Bigl(1 - \frac{l}{p - 1}\Bigr)^4 - p \ge \frac{6 l^2}{1 - y}
\]
(all higher terms in $l$ are nonnegative), whence $l \le \omega s$; on the other hand, \eqref{eq:kernel-spread} gives $l < 0.34$ since $0.34^2(6 + 4\cdot 0.34 + 0.34^2) > 0.85$.

Set $a_j = y q_j(W_*)/q_4(W_*)$, $Q_j = 1 + a_j s^2$, and $M_0 = 1 + y s^2$, so that \eqref{eq:moment-interp} gives $(1 + y\delta)^{(j-1)/3} \le T_j/y^{j-1} \le Q_j$ for $j \ge 4$. Both profiles $\Phi$ and $B$ have the form $V(t) = \sum_{j \ge 1} b_j(1 - \cos(2\pi j t))$ with $b_j \ge 0$ and $\sum_j b_j = \int_0^1 V \le 1$ (specifically $b_j = (Z j^3)^{-1}$ for $\Phi$ and $b_j = 2\epsilon(\sin(\pi j\epsilon)/(\pi j\epsilon))^2$ for $B$, since $1 - B$ is a periodic triangular tent). On a gap of length $x_0$, the first Fourier coefficient of $x_0^m V$ coincides with that of the constant $x_0^m b_V(x_0)$ on the same gap, where
\[
  \frac{\int_0^1 (1 - \cos(2\pi j t))e^{-2\pi i x_0 t}\,dt}{\int_0^1 e^{-2\pi i x_0 t}\,dt}
  = \frac{j^2}{j^2 - x_0^2}, \qquad
  b_V(x_0) = \sum_{j=1}^\infty \frac{b_j}{1 - x_0^2/j^2}.
\]
For $m \in \{3, 4\}$, since $U^2 < 0.05$, we have
\[
  0 \le \frac{(x_0^m b_V(x_0))'}{m x_0^{m-1}}
  \le \frac{1}{1 - U^2} + \frac{2 U^2}{3(1 - U^2)^2} < 1.09.
\]
Assign $\ell = s_i$ probability $s_i$. By \eqref{eq:moment-interp}, $\Var(\ell^m)/y^{2m} \le (a_{2m+1} - 2my/3)s^2$ for $m \in \{3, 4\}$ (squaring the lower bound for $T_{m+1}$). Twice the squared first Fourier coefficient of any real step function is bounded above by its variance; applying the Lipschitz bound to two independent copies of $\ell$ and noting $d(\ell^4)/d(\ell^3) \le 4U/3 < 0.30$ yields $2\max\{F_1^2, (3G_1/(4U))^2\} \le 1.09^2\Var(\ell^3)$. Since $2y = \tau/3$, $a_7 < 0.846\tau$, and $1.09\sqrt{(0.846 - 1/3)/2} < 0.553$, we obtain
\[
  F_1 \le y^3\chi, \qquad G_1 \le 0.30 y^3\chi, \qquad \chi \coloneqq 0.553 s\sqrt{\tau}.
\]
Furthermore, $\Var_m(f) \ge e T_7/2 + \Var(\ell^3)$ by conditional variance, so $2\|f_\perp\|_2^2 = 2\Var_m(f) - 4F_1^2 \ge e T_7 - 0.64 F_1^2$ (since $4(1 - 1/1.09^2) < 0.64$).

Set $c_* = 2.406$, $j_* = 1 + d_* y^3$, $\lambda = 1.25$, $\eta = 0.25$, $k_0 = 2\lambda/(1 + 2\lambda y^2)$, $z_0^2 = 1.017/(2\cdot 1.203^2)$, and $\alpha = 1 - \tau/20$. Replacing $w$ by $w_0$ for brevity in the envelopes below, define
\begin{align*}
  a &= \alpha\Bigl\{c_*\bigl[\tfrac{7}{8}w_0^2 + 1.75 y w_0^3 + y^2 w_0^4 - y\bigl(1 + (1 - \epsilon)y\bigr)\bigr] + 4.6 c_*\chi y\Bigr\}, \\
  r &= 1 - 2a y^2, \qquad d_1 = (1 + \sqrt{r})^2, \\
  H &= \alpha^2 c_*^2\left[\sqrt{\frac{1.0824 Q_5}{2\cdot 1.202^2}}
    + y\sqrt{0.353 + (1.353 a_7 - 2y)s^2}
    + y^2\sqrt{0.002 + (1.002 a_9 - 8y/3)s^2}\right]^2, \\
  E &= \alpha^2\bigl[c_* M_0 + 4(1.0824)y Q_5 + K(U)y^2 Q_6\bigr], \qquad
  m = \frac{2Hy}{d_1} + \sqrt{\frac{d_* E}{r}}, \\
  h_0 &= \alpha\bigl[1.387 w_0 + c_* y^2(M_0 + y Q_5)\bigr], \qquad
  v_r = 1 + \frac{2(1 - r)}{3r}.
\end{align*}
Scaling both $h$ and the tangential trial $-f_\perp/c$ by $\alpha$, Lemma~\ref{lem:strip-phase-pointwise} and Lemma~\ref{lem:common-profile-bounds} ensure that $-a y^2 \le \alpha h \le (a + 2\lambda)y^2$, $h_2(\alpha h) \le H y^4$, $\E(\alpha h) \le E y^3$, and $M_A \le m y^3$. Indeed, $M_0, Q_5 \le w_0$ (as $a_5 \le Y q_5(W_*)/q_4(W_*) < 0.313$ and $\omega > 0.37$), and the lower barrier of Lemma~\ref{lem:strip-phase-pointwise} gives
\[
  \frac{\alpha h/y^2 - a}{\alpha c_*}
  \le 0.495 w_0^2 + y(1 + M_0 - 1.75 w_0^3) + y^2(1 - \epsilon + Q_5 - w_0^4)
  \le 0.495(1.38)^2 + \frac{y}{4} + y^2 < \frac{2\lambda}{c_*},
\]
so $|\alpha h - \lambda y^2| \le (a + \lambda)y^2$ for the tangent estimate \eqref{eq:phase8tangent} at $t = \lambda y^2$. The bound for $H$ follows from centering $f$ and $g$ before applying the triangle inequality to $G$, noting that the within-gap profile variances are bounded by $0.353$ and $0.002$ while the remaining contributions are the gap-power variances bounded above.

Combining \eqref{eq:phase8core}--\eqref{eq:phase8schur} for the $\alpha$-scaled trial gives
\[
  R \ge P(y, \delta) - \mathcal{L}(y, s) + \frac{4Z\alpha\mathcal{R}}{y^4}.
\]
Here the leading term $2Z(2\alpha - \alpha^2)M$ leaves $-2Z(1 - \alpha)^2/y$ after subtracting $2Z y^3$ and dividing by $y^4$. Setting $g_\alpha = 4(1.202)\alpha(1 - \epsilon) - 2(1.0824)\alpha^2$, $a_0 = 1 + (1 - \epsilon)y$, and $a_1 = y(1 + 4(1 - \epsilon)y/3)$, and using $\E(h) \ge (M + \bar{g})^2/(c^2 d)$ before scaling, we obtain
\begin{align*}
  P \ge{}& -c_*\frac{(1 - \alpha)^2}{y} + 2.404(2\alpha - \alpha^2)\delta + g_\alpha\Bigl(1 + \frac{4}{3}y\delta\Bigr) \\
  &+ \alpha^2\Bigl[2.072(2 - \eta)y\Bigl(1 + \frac{5}{3}y\delta\Bigr) + 2.034 y^2(1 + 2y\delta)
    + 2(1 - \eta)(2.404)^2 z_0^2 y^2\Bigr] \\
  &+ \alpha^2 k_0(2.404)^2 y\,\frac{a_0^2 + 2a_0 a_1\delta}{d_*}.
\end{align*}
With $c_\alpha = \alpha c_*$ and $\E(g) \le 0.77 T_8$, the accumulated error term is
\begin{align*}
  \mathcal{L}(y, s) ={}& 2\alpha c_\alpha(0.031)y Q_6 + c_\alpha^2(0.77)y^3 Q_8 + d_* E y^2 + k_0 H y^2
    + j_* h_0(v_r H + \lambda^2)y \\
  &+ \frac{2j_* c_\alpha(a + \lambda)\sqrt{0.77 Q_8 E}}{r(1 + 2\lambda y^2)}y^3
    + \frac{4}{\eta}d_* m^2 y^5 + \frac{8H m y^3}{d_1} \\
  &+ (1 - \eta)\left[3.84 c_\alpha^2 \chi^2 y^2 + \frac{4c_\alpha(M_0 + 2\chi)H y^3}{d_1}\right].
\end{align*}

Denote evaluation at $Y = 1/6$ (with $s$ fixed) by a superscript $\circ$, and set $k = \alpha/0.95 = 1 + (1 - \tau)/19$. We claim that
\[
  \mathcal{L}(y, s) \le \tau k^3 \mathcal{L}^\circ(s).
\]
Among the unscaled envelopes, the only decreasing input in $y$ is $w_0 = 1 + \omega s$, which satisfies $-w_0' \le s/4$, $w_0 < 1.38$, and $w_0 - 1 > 0.37 s$ for $s > 0$. Consequently, $y(w_0^j - 1)$ increases with $y$ for $1 \le j \le 4$ (as $w_0^j - 1 \ge (w_0 - 1)w_0^{j-1}$ and $jy|w_0'| \le s/6$). Since $\frac{7}{8}w_0^2 + 0.75 y$ also increases in $y$, grouping the terms of $a/\alpha$ yields $a \le k a^\circ$; for $h_0$, we have $1.387(w_0 - w_0^\circ) \le 0.33(Y - y) \le c_*(Y^2 - y^2)$, while all other unscaled numerators increase directly in $y$. Because $1 - \tau k = (1 - \tau)(1 - \tau/19) \ge 0$, we have $r \ge r^\circ$. Thus $H$ and $E$ grow by at most $k^2$, while $a + \lambda$, $m$, and $h_0$ grow by at most $k$. Factoring out one power of $\tau$ and noting that $y/(1 + 2\lambda y^2)$ increases on $(0, 1/6]$, every summand of $\mathcal{L}$ grows by at most $k^3$, establishing the claim.

Moreover, $\mathcal{L}^\circ(s)$ is increasing and convex on $[0, S]$: every radical in $\sqrt{H}$, $\sqrt{E}$, and $\sqrt{Q_8}$ has the form $\sqrt{A + B s^2}$ with $A, B \ge 0$; $r^{-1}$, $r^{-1/2}$, and $d_1^{-1}$ are increasing convex functions of $1 - r$; and nonnegative products and increasing convex compositions preserve convexity. The moment recurrence yields $(a_5, a_6, a_7, a_8, a_9) < (0.313, 0.531, 0.846, 1.291, 1.909)$ and $K(U) < 7.41$, while $\omega < 0.390$ at $y = Y$. Direct evaluation at the three nodes $s = 0$, $7/20$, and $S$ gives the strict upper bounds
\begin{align*}
  s = 0:\quad & H < 2.71, & E &< 3.013, & 1/r &< 1.147, & \mathcal{L}^\circ &< 1.85; \\
  s = 7/20:\quad & H < 2.93, & E &< 3.096, & 1/r &< 1.251, & \mathcal{L}^\circ &< 2.27; \\
  s = S:\quad & H < 4.242, & E &< 3.588, & 1/r &< 1.522, & \mathcal{L}^\circ &< 4.10.
\end{align*}
(For example, at $s = 7/20$ one has $h_0 < 1.577$ and $m < 3.62$, and the first four summands, the two tangent terms, the two mean terms, and the final bracket of $\mathcal{L}^\circ$ are bounded by $0.484$, $1.330$, $0.205$, $0.189$, and $0.059$, respectively.) Set $Q(s) = 1.99 - 0.12 s + 3.48 s^2$. At the three nodes $s \in \{0, 7/20, S\}$, the margins $Q(s) - \mathcal{L}^\circ(s)$ exceed $0.14$, $0.10$, and $0.75$. On each subinterval of length $L_0$, the endpoint margins $e_0, e_1$ satisfy $\sqrt{e_0} + \sqrt{e_1} > \sqrt{3.48}\,L_0$; since $(1 - z)e_0 + z e_1 \ge z(1 - z)(\sqrt{e_0} + \sqrt{e_1})^2$ dominates the quadratic sag $3.48 L_0^2 z(1 - z)$ of $Q$ below its secant while $\mathcal{L}^\circ$ is convex, we conclude that $\mathcal{L}^\circ(s) < Q(s)$ for all $s \in [0, S]$.

Now set $t = 1 - \tau \in [0, 1)$, $c_1 = (1 - \epsilon)/6$, and $b(\tau) = 0.95^2 k_0(2.404)^2/(6 d_*)$ evaluated at $y = \tau/6$. Then $k = 1 + t/19$, $b(1) > 0.708$ (since $k_0(1) = 180/77$), and
\[
  \frac{b(\tau)}{b(1)}
  = \Bigl(1 + \frac{2.78 t}{6 d_*}\Bigr)\Bigl(1 + \frac{5(1 + \tau)t}{72 + 5\tau^2}\Bigr) > 1 + 0.225 t \qquad (t > 0),
\]
because the two linear coefficients in $t$ exceed $0.161$ and $0.064$. Writing $P_0 + P_2 s^2$ for the displayed lower bound on $P$, expansion gives $g_\alpha - 6c_*(1 - \alpha)^2/\tau > 2.631 - 0.0595\tau - 0.0055\tau^2$, while the remaining constant terms exceed $k^2\tau[0.545 + 0.127\tau + b(\tau)(1 + c_1\tau)^2]$ (using $z_0^2 > 0.35136$). Since $k^2 \ge 1 + 0.105 t$ and $0.105(1.253 + 0.362\tau) + 0.225(0.708)(1 + 2c_1\tau) > 0.290 + 0.090\tau$, we find
\[
  P_0 > 2.631 - 0.065\tau + \tau(1.253 + 0.362\tau + 0.019\tau^2) + \tau t(0.290 + 0.090\tau).
\]
Similarly, using $g_\alpha \ge g_{0.95} > 2.6024$, $1 + c_1\tau \ge \tau(1 + c_1)$, and $0.545(5/18) + 0.95^2(2.034)/108 + (0.708/3)(1 + c_1)(1 + 4c_1/3) > 0.504$ gives $P_2 > 2.3979 + 0.5783\tau + 0.504\tau^4$. Subtracting $3.48\tau k^3$ produces a function decreasing on $[0,1]$ (since $(\tau k^3)' = k^2(k - 3\tau/19) \ge 16/19$ and $3.48(16/19) > 0.5783 + 4(0.504)$) whose value at $\tau = 1$ is $0.0002 > 0$.

For the cross-gap term $C = 4Z\alpha\mathcal{R}/y^4$, apply Lemma~\ref{lem:cross-chord} with $U = y\min(1.34, w_0)$. Since $4y\omega/(1 - y) \le 0.30\tau$, Bernoulli's inequality gives $(1 - U)^4/(1 - y)^3 \ge (1 - y)(1 - 0.30\tau s)$. With $4Z\pi^2(0.95)/6 > 7.5$, the two bracketed summands at $\tau = 1$ (prior to the factor $1 - 0.30\tau s$) exceed $7.5(0.0395)(0.958)(0.758) > 0.215$ and $7.5(0.978)(0.977)(0.0605)(5/6) > 0.361$. For the first summand, $1 - 0.042\tau^2$ and $(1 + 0.34\tau/6)^{-5}$ have ratios to their values at $\tau = 1$ at least $1 + 0.042 t$ and $1 + 0.268 t$; for the second, $1 - 0.022\tau^2$, $1 - 0.023\tau^2$, $1/12 - 0.137\tau/6$, and $1 - \tau/6$ give linear coefficients in $t$ summing to at least $0.022 + 0.023 + (0.137/6)/0.0605 + 0.20 > 0.622$. Because $\tau(1 + 0.622 t) \le 1 - 0.378 t$ implies $k\tau(1 + 0.622 t) \le 1$, it follows that
\[
  \frac{C}{\tau} > k\bigl[0.215(1 + 0.310 t) + 0.361(1 + 0.622 t)(1 - 0.30\tau s)\bigr]
  > 0.576 + 0.321 t - 0.109 s.
\]

Finally, the chord bound $1.99 k^3 \le 1.99 + 0.333 t$ substituted into $R \ge P - \tau k^3 Q(s) + C$ yields
\[
  R > 2.631 + 0.052\tau + 0.174\tau^2 - 0.071\tau^3 + 0.011\tau s.
\]
The cubic polynomial is strictly increasing on $[0, 1]$ (its derivative is at least $0.052 + 0.135\tau > 0$), exceeding $2.78$ at $\tau = 1$ ($p = 6$), $2.75$ at $\tau = 6/7$ ($p = 7$), and $2.60$ throughout $(0, 1]$, contradicting Corollary~\ref{cor:regular-uniform-nine} for every $p \ge 6$.
\end{proof}

\section{One reflex vertex with three, four, or five convex vertices}\label{sec:low-reflex}

It remains to establish the strict comparison $J(P) < J(R_{p+1})$ for polygons with a single reflex vertex and $p \in \{3, 4, 5\}$ strictly convex vertices. When $p = 3$, an elementary convex-hull pocket decomposition bounds $J(P)$ strictly below $3\sqrt{3} < J(R_4)$, independently of the number of reflex vertices. For $p = 4$ and $p = 5$, gap dispersion forces the convex prevertices into a concrete neighborhood of the regular $p$-tuple, where the pulled-back reciprocal regular density $w_{-1/p}(B)$ from Subsection~\ref{sec:reciprocal-transport} contradicts $d(P) \le d(R_{p+1})$.

\subsection{Three strictly convex vertices}

\begin{prop}[Three strictly convex vertices]\label{prop:hull}
Every simple planar polygon $P$ with exactly three strictly convex vertices satisfies
\[
  J(P) = |\partial P|\,\sigma_1(P) < 3\sqrt{3} < J(R_4).
\]
\end{prop}

\begin{proof}
For any polygon $D$, write $I_D = \min_{z \in \R^2} \int_{\partial D} |X - z|^2\,ds$. Testing \eqref{eq:variational-physical} with the two Cartesian coordinates centered at the boundary centroid yields $J(D) I_D \le 2|D|\,|\partial D|$. We first prove the geometric inequality
\[
  I_P \ge b|P|\,|\partial P|, \qquad b \coloneqq \frac{2}{3\sqrt{3}},
\]
with equality only for an equilateral triangle.

The three vertices of the convex hull $T$ of $P$ are necessarily strictly convex vertices of $P$, and therefore exhaust all strictly convex vertices of $P$. They partition $\partial P$ into three simple arcs; the arc opposite vertex $i$ together with the hull edge of length $a_i$ encloses a (possibly degenerate) convex pocket $K_i$, since every interior vertex along the arc is reflex and the hull edge supports $K_i$. Rescale so that $|\partial T| = 2$, set $\Delta = |T|$, and write $|K_i| = \Delta\alpha_i$ and $\lambda = \sum_{i=1}^3 \alpha_i < 1$, so that $|P| = (1 - \lambda)\Delta$.

Let $h_i$ be the signed distance from hull edge $i$, positive inside $T$, and put $M_i = \max_{K_i} h_i$. Because $K_i$ contains a triangle of base $a_i$ and altitude $M_i$, we have $M_i/(2\Delta/a_i) \le \alpha_i$. On the other hand, $K_i$ is contained in the trapezoid $T \cap \{h_i \le M_i\}$; monotonicity of convex perimeter and similarity of the excised triangle $T \cap \{h_i \ge M_i\}$ yield
\[
  |\partial P| \le L_+ \coloneqq 2 + \sum_{i=1}^3 (2 - 2a_i)\alpha_i.
\]

Over the open base interval $|t| < a_i/2$ measured from the midpoint of hull edge $i$, convexity implies that the normal cross-section of $K_i$ is $[0, f_i(t)]$, where the graph of $f_i$ lies on the replacing arc of $\partial P$ (with $f_i \equiv 0$ if $K_i$ is degenerate). Along these graphs $ds \ge dt$ and $\sum_{i=1}^3 \int_{-a_i/2}^{a_i/2} f_i(t)\,dt \le \lambda\Delta$. For any center $z \in \R^2$, the tangential component of $|X - z|^2$ contributes at least $\frac{1}{12}\sum_{i=1}^3 a_i^3$ to $I_P$. For the normal component, Cauchy--Schwarz over the union of the three base intervals (of total length $|\partial T|$), combined with the signed-area identity $\sum_i a_i h_i(z) = 2|T|$ and the subtraction of at most $|T| - |P|$ from the pockets, leaves an integral of at least $|T| + |P| = (2 - \lambda)\Delta$. Hence
\[
  I_P \ge \frac{1}{12}\sum_{i=1}^3 a_i^3 + \frac{(|T| + |P|)^2}{|\partial T|}.
\]
(Note that AM--GM and $(\sum_i a_i^2)^2 \le 2\sum_i a_i^3$ also give $2\sqrt{3}\,I_P \ge (|T| + |P|)\sum_i a_i^2$.)

Set $v_i = a_i - 2/3$, $v = \sum_{i=1}^3 v_i^2$, and $t = 3\sqrt{3}\,\Delta$. Heron's identity $\Delta^2 = \prod_{i=1}^3 (1 - a_i)$ implies $t \le 1$ by AM--GM and $\sum_{i=1}^3 a_i^3 = 1 + \frac{3}{2}v - 3\Delta^2$ by expansion. Substituting into the lower bound for $I_P$ gives
\[
  108 I_P \ge 9 + \frac{27}{2}v + (7 - 8\lambda + 2\lambda^2)t^2
  \ge 2(1 - t) + \frac{27}{2}v + 8t(2 - \lambda),
\]
because the difference between the two right-hand sides equals $7(1 - t)^2 + 8t\lambda(1 - t) + 2t^2\lambda^2 \ge 0$. Since $\sum_i v_i = 0$, we have $|v_i| \le \sqrt{2v/3}$, so the perimeter bound $L_+$ yields $108 b|P|\,|\partial P| \le 16 t(1 - \lambda)(1 + \lambda/3 + \lambda\sqrt{2v/3})$. Subtracting and applying $2\lambda\sqrt{2v/3} \le \frac{2}{3}\lambda^2 + v$ gives
\[
  108\bigl(I_P - b|P|\,|\partial P|\bigr)
  \ge 2(1 - t) + \frac{8}{3}t\lambda + \Bigl(\frac{27}{2} - 8t\Bigr)v
  \ge 2(v + \lambda),
\]
using $(1 - t) + t\lambda \ge \lambda$. Thus $I_P \ge b|P|\,|\partial P|$ holds with equality only for an equilateral triangle.

Finally, the summed coordinate trial $J(P) I_P \le 2|P|\,|\partial P|$ is strict on every polygon: equality would force $\nu = \sigma_1(P)(X - z)$ along every edge, which is impossible since $\nu$ is constant on an edge while $X - z$ varies. Therefore $J(P) < 2/b = 3\sqrt{3}$. By Lemma~\ref{lem:reference-square}, $J(R_4) > 5.5 > 3\sqrt{3}$.
\end{proof}

\subsection{Gap dispersion and Blaschke multiplier estimates for \texorpdfstring{$p = 4, 5$}{p = 4, 5}}
We first narrow the balanced convex gaps for $p = 4$ and $p = 5$ using Lemma~\ref{lem:kernel-dispersion}.

\begin{lem}[Gap dispersion for four and five convex vertices]\label{lem:marked-gap-masses}
Let $P$ be a simple polygon with $p \in \{4, 5\}$ strictly convex vertices, one reflex vertex, and balanced convex gaps $s_1, \dots, s_p$.
\begin{enumerate}[label=\textup{(\roman*)},leftmargin=*]
\item If $p = 4$ and $d(P) \le 0.05215$, then $s_i \in (0.08, 3/8)$ for every $i$ and $\sum_{i=1}^4 (s_i - 1/4)^2 < 0.021125$.
\item If $p = 5$ and $d(P) \le 0.02420$, then $s_i \in (0.08, 0.30)$ for every $i$ and $\sum_{i=1}^5 (s_i - 1/5)^2 < 0.00908$.
\end{enumerate}
\end{lem}

\begin{proof}
Retain $y = 1/p$, $T_j$, $M = T_4$, $u = \max_i s_i$, $w = u/y$, and $\delta = (M - y^3)/y^4$ from Lemma~\ref{lem:kernel-dispersion}. First, $w < 5/3$ in both cases: by convexity in \eqref{eq:gap}, setting one gap equal to $5y/3$ and the remaining $p - 1$ gaps equal minimizes $\sum_i \mathcal{H}(s_i)$ subject to $u \ge 5y/3$; using $D(s) < 0.416 + 0.750 s - 0.70 s^3$ yields $0.03014/0.6779 + 3(0.001429)/0.557 > 0.05215$ for $p = 4$ and $1/(81\cdot 0.641) + 4/(1296\cdot 0.54) > 0.02420$ for $p = 5$.

Taking $k = 9/11 \le \kappa(5/12)$ in the tangent bound following Lemma~\ref{lem:kernel-dispersion} gives $A_k > 2.22$. This forces $u < 3y/2$ in both cases: otherwise, convexity gives $M \ge u^4 + (1 - u)^4/(p - 1)^3 > 0.025$ (for $p = 4$) and $> 0.011$ (for $p = 5$), so $2.22 M$ already exceeds the assumed upper bound on $d(P)$.

With $u < 3/8$ ($p = 4$) and $u < 3/10$ ($p = 5$), Lemma~\ref{lem:single-gap-loss} gives $\kappa(3/8) > 8/9$ and $\kappa(3/10) > 63/64$, because the corrections to $1$ in their reciprocals are $\pi/12 - 9/(16\pi Z) < 1/8$ and $\pi/60 - 3/(20\pi Z) < 1/63$, respectively. For these choices of $k$, the tangent bound gives $A_k > 2.37$ ($p = 4$) and $A_k > 2.404 + 0.38 m$ ($p = 5$). Using $T_5 \ge y M$ in \eqref{eq:direct-dispersion} yields $d(P) \ge 2.37 M$ ($p = 4$) and $d(P) \ge 2.48 M$ ($p = 5$), forcing $\delta < 1.64$ and $\delta < 1.10$, respectively.

Set $x_i = s_i/y - 1$, $V = \sum_{i=1}^p x_i^2$, $A = \sqrt{(p - 1)V/p}$, and $b_0 = A/(p - 1)$. Since $\sum_i x_i = 0$, every $x_i$ satisfies $|x_i| \le A$; moreover, $3V \le \delta$ implies $A < 1$. For $g(x) = x^4 + 4x^3 + 6x^2$, the polynomial
\[
  Q(x) = g(x) - (x + A)(x - b_0)^2(x + 4 - A + 2b_0)
\]
has degree $2$ and satisfies $Q(x) \le g(x)$ for all $x \ge -A$, with equality at $x = -A$ and $x = b_0$. Because the vector $(-A, b_0, \dots, b_0)$ shares the same first two moments as $(x_1, \dots, x_p)$, we obtain
\[
  \delta = \sum_{i=1}^p g(x_i) \ge g(-A) + (p - 1)g(b_0)
  = V\bigl[6 - 4(A - b_0) + (A - b_0)^2 + A b_0\bigr].
\]
Writing $K = (p - 2)/\sqrt{p(p - 1)}$, we have $A - b_0 = K\sqrt{V}$ and $A b_0 = V/p$. This lower bound increases on $V \in [0, 1]$ and exceeds $1.64$ at $V = 0.338$ ($p = 4$, using $\sqrt{0.338/3} < 0.336$) and $1.10$ at $V = 0.227$ ($p = 5$, using $\sqrt{0.227} < 0.477$). Thus $\sum_i (s_i - y)^2 = y^2 V < 0.021125$ and $0.00908$, respectively, and $s_i \ge y(1 - A) > 0.08$ in both cases.
\end{proof}

Next, we derive an exact Euler integral representation and an energy product estimate for the Blaschke multipliers appearing in the centered $p = 4$ and $p = 5$ models.

\begin{lem}[Euler integral for the Blaschke multiplier]\label{lem:marker-multiplier}
Let $p \ge 3$ be an integer, and put $\beta = 1/p$, $x = z^2$, $t \in (0,1)$, $q = t^2$, $k = (1 - q)^{-4\beta^2}$, and $D = |1 - tx|^2$. Define
\[
  \phi_t(x) = \frac{t - x}{1 - tx}, \qquad
  B = \eta z^{p-2}\phi_t(x), \qquad
  u = 1 - 2\beta + \frac{2\beta(1 - q)}{D}, \qquad
  V = k u D^\beta,
\]
with $|\eta| = 1$. For integers $j \ge 0$ and $n$ such that $a = ((p - 2)j - n)/2$ is a positive integer,
\begin{equation}\label{eq:five-euler-mode}
  \widehat{V B^j}(n)
  = -\frac{\eta^j k\beta\sin(\pi\beta)\,t^{a+j}}{\pi(j + \beta)}
    \int_0^1 v^{a - \beta - 1}(1 - v)^{j + \beta}(1 - qv)^{\beta - j}
    \Bigl(n + 1 + \frac{4\beta qv}{1 - qv}\Bigr)\,dv,
\end{equation}
while Fourier coefficients of the opposite parity vanish. Writing $v_m = \widehat{V}(2m)$, for $p = 5$ and $t \le 0.605$ one has
\[
  1 + 0.04 q \le v_0 \le k(1 - 0.12 q), \qquad
  0 \le v_1 \le 0.2 k t, \qquad
  0 \le v_m \le 0.24 k t^m \quad (m \ge 2).
\]
\end{lem}

\begin{proof}
Set $L = (1 - t z^2)^\beta$ and $X = z L'/L$, so $V = k u|L|^2$ and $u = 1 - X - \bar{X}$. Integrating the angular derivative of $W = z^{-n} B^j |L|^2$ around $\T$ yields $0 = \int_\T W(-n + pju + X - \bar{X})\,dm$, which rearranges to
\begin{equation}\label{eq:marker-transfer}
  (1 + pj)\widehat{V B^j}(n)
  = k\int_\T z^{-n} B^j |L|^2\Bigl(n + 1 - \frac{2z L'}{L}\Bigr)\,dm.
\end{equation}
For any function $f$ holomorphic on a neighborhood of $\overline{\D}$ except possibly for a pole at $t$, shrinking the Cauchy contour for $f(x)(1 - t/x)^\gamma/x$ onto the two sides of $[0, t]$ yields
\begin{equation}\label{eq:five-euler-integral}
  \int_\T f(x)(1 - t\bar{x})^\gamma\,dm
  = -\frac{\sin(\pi\gamma)}{\pi}\int_0^1 f(tv) v^{-\gamma - 1}(1 - v)^\gamma\,dv
\end{equation}
whenever $\operatorname{ord}_0 f > \gamma$ and $\operatorname{ord}_t f > -1 - \gamma$. In \eqref{eq:marker-transfer}, substituting $x = z^2$ gives $z^{-n} B^j = \eta^j x^a \phi_t(x)^j$ in the stated parity (and zero integral for opposite parity). Applying \eqref{eq:five-euler-integral} with $\gamma = \beta$ and $-2z L'/L = 4\beta tx/(1 - tx)$ establishes \eqref{eq:five-euler-mode}.

Now fix $p = 5$ and $t \le 0.605$. For $j = 0$ and $n = -2m$ ($m \ge 1$), the integrand in \eqref{eq:five-euler-mode} is strictly negative since $4\beta q/(1 - q) < 0.48$. Taking $f(x) = x^m$ in \eqref{eq:five-euler-integral} identifies the coefficient of $t^m$ on the left as $-d_m$, where $(1 - x)^\beta = 1 - \sum_{m=1}^\infty d_m x^m$; symmetry and \eqref{eq:five-euler-mode} thus give $0 \le v_m = v_{-m} \le k t^m (2m - 1)d_m$. From $d_1 = \beta$ and $d_{m+1}/d_m = (m - \beta)/(m + 1)$, we have $(2m - 1)d_m \le 3d_2 = 0.24$ for $m \ge 2$. Expanding $V/k = |L|^2 - 2\Re(z L'\overline{L})$ via $L = 1 - \sum_{m \ge 1} d_m t^m z^{2m}$ gives $v_0/k = 1 - \sum_{m \ge 1} (4m - 1)d_m^2 q^m \le 1 - 0.12 q$. Since $(4m - 1)d_m^2$ decreases with $m$, we also obtain $v_0/k \ge 1 - 0.12 q - 0.0448 q^2/(1 - q) \ge (1 - q)^{0.16}(1 + 0.04 q)$ by expanding $(1 - q)^{0.16}$ to second order.
\end{proof}

\begin{lem}[Energy of the four- and five-vertex multipliers]\label{lem:marker-gamma}
In the setting of Lemma~\ref{lem:marker-multiplier}, take $(p, \ell) = (4, 1)$ or $(5, 2)$, allow $t = 0$, assume $4\beta q < 1 - q$, and set
\[
  R = \max_{0 \le v \le 1} (tv)^{(p-2)\ell/2}\left(\frac{t(1 - v)}{1 - qv}\right)^\ell, \qquad
  \Gamma(V)(z) = \int_\T \frac{|V(z) - V(w)|^2}{|z - w|^2}\,dm(w).
\]
Then for every $F \in L^2(\T)$,
\begin{equation}\label{eq:marker-gamma-push}
  \int_\T \Gamma(V)|F(B)|^2\,dm
  \le \left(1 + \frac{2R}{(1 + p\ell)(1 - R)}\right)\E(V)\int_\T |F|^2\,dm,
\end{equation}
and for every real $f = F(B)$ of finite energy,
\begin{equation}\label{eq:marker-product}
  \E(Vf) \le \E(f, V^2 f) + \int_\T \Gamma(V)f^2\,dm.
\end{equation}
In particular, if $t \le 0.525$ for $p = 4$ or $t \le 0.605$ for $p = 5$, then
\begin{equation}\label{eq:covariant-marker-product}
  \E(Vf) \le \E(f, V^2 f) + 0.575 q\int_\T F^2\,dm.
\end{equation}
\end{lem}

\begin{proof}
Write $v_m = \widehat{V}(2m) = v_{-m}$ and $U_j = V B^j$. For $m \ge 1$, formula \eqref{eq:five-euler-mode} at $n = -2m$, $j = 0$ expresses $v_m$ as a positive multiple of the integral of the strictly positive weight
\[
  v^{m - \beta - 1}(1 - v)^\beta (1 - qv)^\beta \Bigl(2m - 1 - \frac{4\beta qv}{1 - qv}\Bigr).
\]
At $j = \ell s$ ($s \ge 1$), the integrand picks up the $s$-th power of the expression maximized in $R$, while the prefactor multiplies by $\eta^{\ell s}/(1 + p\ell s)$; thus $|\widehat{U}_{\ell s}(-2m)| \le v_m R^s/(1 + p\ell s)$.

From $\Gamma(V) = V\Lambda V - \frac{1}{2}\Lambda(V^2)$, we have $\widehat{\Gamma(V)}(2n) = 4\sum_{m=1}^\infty m v_m v_{m+n}$ for $n \ge 0$. Because $B$ is holomorphic,
\[
  \left|\int_\T \Gamma(V) B^{\ell s}\,dm\right|
  = \left|4\sum_{m=1}^\infty m v_m \widehat{U}_{\ell s}(-2m)\right|
  \le \frac{R^s}{1 + p\ell s}\E(V).
\]
(In fact, for $t > 0$, $(B^\ell)_*(\Gamma(V)\,dm)/\E(V)$ is a mixture of Poisson kernels centered at $\eta^\ell a^{p\ell} b(v)$, obtained by sampling $m$ with weight $4m v_m^2/\E(V)$, $v$ from its normalized Euler weight, and $a \sim \mathrm{Unif}[0,1]$ independently.) At $p = 5$, odd $B$-moments vanish under $z \mapsto -z$. Summing the absolute Fourier coefficients of $B_*(\Gamma(V)\,dm)$ proves \eqref{eq:marker-gamma-push}. Inequality \eqref{eq:marker-product} follows immediately from the difference-kernel representation of $\E(Vf) - \E(f, V^2 f)$, whose integrand contains $f(z)f(w)|V(z) - V(w)|^2 \le \frac{1}{2}(f(z)^2 + f(w)^2)|V(z) - V(w)|^2$.

When $t = 0$, $V \equiv 1$ and \eqref{eq:covariant-marker-product} is immediate. For $t > 0$, expanding $V/k = |L|^2 - 2\Re(z L'\overline{L})$ yields only negative cross terms in positive modes, so retaining the first two terms gives
\[
  0 \le \frac{v_m}{k t^m} \le (2m - 1)d_m - (2m + 3)\beta d_{m+1} q.
\]
For $m \ge 2$, $(2m - 1)d_m$ decreases from $M_2 = 3d_2$. At $m = 1, 2$, the ratios of the $q$-correction coefficients to the leading coefficients are $5d_2 \ge 0.4$ and $7\beta(2 - \beta)/9 \ge 0.28$, so since $q < 0.367$, squaring improves the $m = 1, 2$ bounds by factors $1 - 0.7q$ and $1 - 0.5q$. Hence
\[
  \frac{\E(V)}{4k^2 q} = \sum_{m=1}^\infty \frac{m v_m^2}{k^2 q}
  \le \beta^2(1 - 0.7q) + M_2^2\bigl((1 - q)^{-2} - 1 - q^2\bigr).
\]
Evaluating at $q < 0.276$, $k^2 < 1.176$ ($p = 4$) and $q < 0.367$, $k^2 < 1.166$ ($p = 5$) gives $\E(V)/q < 0.550$ and $< 0.506$, respectively. Using $v(1 - v) \le 1/4$ bounds $R \le q/[4(1 - q)]$ ($p = 4$) and $R \le t(q/[4(1 - q)])^2$ ($p = 5$), both below $0.1$. Substituting into \eqref{eq:marker-gamma-push}--\eqref{eq:marker-product} (with $(47/45)(0.550) < 0.575$) proves \eqref{eq:covariant-marker-product}.
\end{proof}

We also collect the harmonic-analytic properties of the signed measure $F\Lambda F$ associated with the reciprocal test $F = w_{-\beta} - D_\beta$ from \eqref{eq:reciprocal-test}.

\begin{lem}[The reciprocal-density energy measure $F\Lambda F$]\label{lem:reciprocal-measure}
Fix $\beta \in (0, 1/2)$, set $X = |1 - z|^2$ for $z$ uniform on $\T$, and define $w = w_{-\beta} = X^\beta/\mathbb{E} X^\beta$, $F = w - D_\beta$, $C = \bar{F} = 1 - D_\beta$, and $b_j = -r_j(-\beta)$.
\begin{enumerate}[label=\textup{(\roman*)},leftmargin=*]
\item Let $T$ be independent of $X$ with density $\beta t^{\beta-1}$ on $(0,1)$. Then $\Lambda F = \Lambda w$ is represented by the finite even signed measure satisfying
\[
  \int_\T J(X)\,d\Lambda w = \beta\,\mathbb{E}\bigl[J(X) - J(TX)\bigr]
\]
for every continuous $J$ on $[0, 4]$. Its density with respect to $dm$ is increasing in $X$ and has mean zero, so for every bounded nonnegative increasing function $\psi$ on $[0, 4]$,
\begin{equation}\label{eq:chord-monotone-pairing}
  \int_\T \psi(X) F\,d\Lambda F \ge C\int_\T \psi(X)\,d\Lambda F \ge 0.
\end{equation}
\item The negative total variation of the signed measure $F\Lambda F$ satisfies
\begin{equation}\label{eq:reciprocal-negative-mass}
  \|(F\Lambda F)_-\| \le \frac{\beta}{4}\bigl(3\|F\|_1 - C\bigr),
\end{equation}
and its positive mass equals $\|(F\Lambda F)_-\| + \beta/2$.
\item The Fourier coefficients $h_j = \widehat{F\Lambda F}(j)$ are given for $j = 0, 1, 2$ by
\begin{equation}\label{eq:small-reciprocal-product}
  h_0 = \frac{\beta}{2}, \qquad
  h_1 = D_\beta b_1 - \frac{\beta}{2}, \qquad
  h_2 = 2D_\beta b_2 + \frac{\beta}{2} - \frac{\beta}{(1 + \beta)^2}.
\end{equation}
Whenever $C < (1 - \beta)/2$, every $h_n$ ($n \ge 1$) satisfies $0 < h_n \le D_\beta b_1$.
\end{enumerate}
\end{lem}

\begin{proof}
For $J(X) = X^k$ ($k \ge 1$), both sides of the representation identity equal $\beta k\,\mathbb{E} X^k/(\beta + k)$ by \eqref{eq:reciprocal-mixed-energy}, and both vanish on constants. Since every $\cos(n\theta)$ is a polynomial in $X$ and odd modes vanish by symmetry, the representation holds on all even continuous functions; uniform and energy approximation identifies $F\Lambda F$ with the energy pairing $\E(F, JF)$ on smooth tests.

The density of $X$ on $(0,4)$ is $[\pi\sqrt{x(4-x)}]^{-1}$, so the ratio of the densities of $TX$ and $X$ at $x \in (0,4)$ is
\[
  R(x) = \beta x^{\beta - 1/2}\sqrt{4 - x}\int_x^4 \frac{dy}{y^{\beta + 1/2}\sqrt{4 - y}}.
\]
Since $\beta < 1/2$, every factor decreases in $x$, so $d\Lambda F = \beta(1 - R(X))\,dm$ has an increasing density of mean zero. Under the reweighted probability measure $\psi\,dm/\int \psi\,dm$ (for $\int \psi\,dm > 0$), positive covariance of increasing functions ensures that $F$ and $1 - R(X)$ have means at least $C$ and $0$, respectively, and applying positive covariance once more to their product yields \eqref{eq:chord-monotone-pairing}.

Writing $F = F_+ - F_-$, the bound $F_+(TX) \le T^\beta F_+(X)$ and $\mathbb{E} T^\beta = 1/2$ give $\|(F\Lambda F)_-\| \le \beta\,\mathbb{E}(F_- + F_+/2) = \frac{\beta}{4}(3\|F\|_1 - C)$, proving \eqref{eq:reciprocal-negative-mass}; the positive mass exceeds this by $\int F\,d\Lambda F = \E(F) = \beta/2$.

For (iii), $\mathbb{E}[w X^k] = 4^k(\beta + 1/2)_k/(\beta + 1)_k$ and \eqref{eq:reciprocal-mixed-energy} give $\int X^k w\,d\Lambda w = \beta(\beta + k)\mathbb{E}[w X^k]/(2\beta + k)$ for $k = 0, 1, 2$; expanding $\cos\theta = 1 - X/2$ and $\cos(2\theta) = 1 - 2X + X^2/2$ produces \eqref{eq:small-reciprocal-product}. For general $n \ge 1$, applying \eqref{eq:polarized-energy} at $a = b = -\beta$, $m = n$ to the two tails of the Fourier convolution yields $h_n = ((n - \beta)/2 - nC)b_n + \frac{n}{2}\sum_{j=1}^{n-1} b_j b_{n-j}$. When $C < (1 - \beta)/2$, this is strictly positive; replacing $n/2$ by $j$ in the symmetric sum and using $j b_j \le b_1$ with $\sum_{j=1}^\infty b_j = 1/2$ (from $w(1) = 0$) gives $h_n \le D_\beta b_1$.
\end{proof}

\subsection{Four convex vertices and one reflex vertex}
We now rule out any simple polygon with $p = 4$ strictly convex vertices and $q = 1$ reflex vertex achieving $d(P) \le 0.05215$. In balanced disk coordinates, the two hyperbolic geodesics joining opposite pairs of the four convex prevertices intersect at a unique interior point. Applying a disk automorphism that sends this intersection point to the origin turns both geodesics into diameters; after a rotation, the four prevertices become $\pm e^{i\alpha}, \pm e^{-i\alpha}$, with marker polynomial $1 - 2t z^2 + z^4$, where $t = \cos(2\alpha) \in [0, 1)$.

\begin{lem}[Antipodal normalization for four prevertices]\label{lem:four-hand-dispersion}
Under the hypotheses of Lemma~\ref{lem:marked-gap-masses} with $p = 4$, set $W_0 = \sum_{i=1}^4 (s_i - 1/4)^2 < 0.021125$. In the antipodal coordinates sending the intersection of the two opposite hyperbolic diagonals to the origin, let $t \in [0, 1)$ be the parameter of $1 - 2t z^2 + z^4$ and let $r$ be the modulus of the image of the balanced origin. Then $t \le \sin(\pi\sqrt{W_0}) < 0.46$ and $r < 1/4$.
\end{lem}

\begin{proof}
Write the two pairs of opposite balanced gaps as $1/4 + h \pm x$ and $1/4 - h \pm y$, and set $s = 2\pi|h|$, $u = 2\pi^2(x^2 + y^2)$, and $A = s^2 + u = \pi^2 W_0 < 0.209$. The products of opposite chord lengths are proportional to $\sin^2(\pi(1/4 \pm h)) - \sin^2(\pi x_\pm)$ (where $x_+ = x$, $x_- = y$). By cross-ratio invariance, $t$ equals the absolute difference of these two products divided by their sum, which yields $t \le \sin s + (1 + t)u/2$. Putting $b = \sqrt{A} < 1/2$ and $T = \sin b$, we have $(1 + T)b < \cos b$, so the mean-value theorem gives $T - \sin s \ge (b - s)\cos b \ge (1 + T)u/2$. Rearranging yields $(1 - u/2)(t - T) \le 0$, whence $t \le T = \sin(\pi\sqrt{W_0}) < 0.46$.

In balanced coordinates, the Klein-model image $K$ of the diagonal intersection has modulus $2r/(1 + r^2)$. The two diagonal chords have unit normals at angle $\pi/2 \pm s$ and signed distances from the origin of moduli $|\sin\pi(x + y)|$ and $|\sin\pi(x - y)|$. Summing the squares of the two normal projections gives $(1 - \sin s)|K|^2 \le \sin^2\pi(x + y) + \sin^2\pi(x - y) \le u$. Setting $d_0 = 2(1 - \sqrt{1 - A})$, we have $u = A - s^2 = d_0(1 - s) - (s - d_0/2)^2 \le d_0(1 - \sin s)$, so $|K|^2 \le d_0$. Because $A < 0.209 < 64/289 - (32/289)^2$, we obtain $d_0 < 64/289$; since $4r^2/(1 + r^2)^2$ increases with $r$ and equals $64/289$ at $r = 1/4$, it follows that $r < 1/4$.
\end{proof}

\begin{prop}[Four convex vertices and one reflex vertex]\label{prop:four-one-reflex}
Every simple polygon $P$ with $p = 4$ strictly convex vertices and $q = 1$ reflex vertex satisfies $d(P) > 0.05215 > d(R_5)$, and in particular $J(P) < J(R_5)$.
\end{prop}

\begin{proof}
Suppose for contradiction that $d(P) \le 0.05215$. By Lemmas~\ref{lem:marked-gap-masses} and~\ref{lem:four-hand-dispersion}, in the antipodal coordinates the parameter $t$ and the displaced balanced origin $a \in \D$ satisfy $q \coloneqq t^2 \le \pi^2 W_0 < 0.209$ (hence $t < 0.46 < 0.525$) and $r = |a| < 1/4$.

Set $x = z^2$, $\kappa = q/(1 - q) < 0.265$, $B = x(t - x)/(1 - tx)$, $h = \log(|1 - tx|/\sqrt{1 - q})$, $u = (1 + e^{-2h})/2$, and $V = u e^{h/2}$. Taking the reciprocal test $F = w_{-1/4} - D_{1/4}$ from \eqref{eq:reciprocal-test} at $\beta = 1/4$, put $G = V F(B)$. Since $V$ is proportional to $u|1 - tx|^{1/2}$, inequality \eqref{eq:regular-potential-transport} gives $\int_\T G\,d\mu \le 0$. Meanwhile, the mean of $G$ in balanced coordinates equals the value $M(a)$ of the harmonic extension of $G$ at $a$.

At $\beta = 1/4$, we have $D_\beta = \pi/4$, $0.2145 < C < 0.215$, $b_1 = 1/5$, and $\E(F) = 1/8$. Parseval's identity yields $\int_\T F^2\,dm \le C^2 + \E(F)/2 + b_1^2 < 0.149$. The two branches of $B = x\phi_t(x)$ are swapped by the disk involution $\phi_t$, which maps $P = e^{-2h}$ to $P^{-1}$. Direct algebra gives $P + P^{-1} - 2 = \kappa|1 - B|^2$, so $\sinh^2 h = \kappa v$ with $v = \sin^2(\phi/2) = X/4$. Since $|B'| = 4u$, the measure $u\,dm$ averages the two branches equally; writing $u = e^{-h}\cosh h$, the fiber averages of $e^{h/2}$ and $V^2$ are $c_0 = \cosh(h/2)$ and $J(v) = \cosh^3 h = (1 + \kappa v)^{3/2}$. Because $J$ is convex on $[0,1]$, $1 + J'(1)v - J(v)$ is nonnegative and increasing; applying \eqref{eq:chord-monotone-pairing} with $\int v F\,d\Lambda F = (h_0 - h_1)/2$ (where $h_0 - h_1 = 1/4 - \pi/20 < 0.093$ by \eqref{eq:small-reciprocal-product}), together with $\Lambda(F(B)) = 4u(\Lambda F)(B)$ and \eqref{eq:covariant-marker-product}, gives
\begin{align*}
  \E(G) &\le 4h_0 + 3\kappa\sqrt{1 + \kappa}\,(h_0 - h_1) + 0.575 q\int_\T F^2\,dm \\
  &< 0.5 + 3(0.265)(1.125)(0.093) + 0.575(0.209)(0.149) < 0.602.
\end{align*}

To bound $M(a)$ from below, write $x' = \phi_t(x)$, $\zeta = x x' = e^{i\phi}$, and denote the average over the two branches by $\av{\cdot}$. From $x + x' = t(1 + \zeta)$ and $(x - x')\sinh h = t(\zeta - 1)\cosh h$ (via $\tanh h = t(t - \Re x)/(1 - t\Re x)$), averaging at $\pm h$ gives
\[
  \av{x e^{h/2}} = t\Bigl(c_0\zeta + \frac{1 - \zeta}{4c_0}\Bigr).
\]
Using $0 \le c_0 - 1 \le \kappa/8$ and $\int \zeta F\,dm = -b_1$ yields $2|\widehat{G}(2)| \le t(|C - 3b_1|/2 + 3\kappa\sqrt{0.149}/8) \le \sqrt{\kappa}/4$. Furthermore, $c_0'(v) = \kappa/(8c_0\sqrt{1 + \kappa v}) \ge \kappa/10$, so $c_0(v) - \kappa v/10$ is increasing and at least $1$; positive covariance with $F$ gives $M(0) \ge C + \kappa(C + b_1)/20 \ge C + \kappa/50$. For $r = |a| < 1/4$, Cauchy--Schwarz bounds the contribution of all even harmonics $|2j| \ge 4$ to the Poisson integral by $\sqrt{\E(G)\sum_{j \ge 2} r^{4j}/j} \le 9 r^4/16$. Completing the square in $\sqrt{\kappa}$ yields
\[
  M(a) \ge C + \frac{\kappa}{50} - \frac{r^2\sqrt{\kappa}}{4} - \frac{9 r^4}{16}
  \ge C - \frac{43 r^4}{32} > 0.2092.
\]
By \eqref{eq:duality} in balanced coordinates, $d(P) \ge M(a)^2/\E(G) > 0.2092^2/0.602 > 0.072$, contradicting $d(P) \le 0.05215$. By Lemma~\ref{lem:regular-six-short}, $d(R_5) < 0.0518 < 0.05215$, so $J(P) < J(R_5)$.
\end{proof}

\subsection{Five convex vertices and one reflex vertex}
When $p = 5$ and $q = 1$, we balance the five convex prevertices themselves (via Lemma~\ref{lem:balanced-coords}) to place them in a normal form with vanishing linear symmetric function, and then test with the pulled-back reciprocal density at $\beta = 1/5$.

\begin{lem}[Centering five boundary points]\label{lem:five-center-hand}
Let five cyclically ordered points on $\T$ have successive angular gaps $g_j$ satisfying $W_0 = \sum_{j=1}^5 (g_j/(2\pi) - 1/5)^2 \le 0.01$. A disk automorphism centers the five points at the origin and maps the original origin to a point of modulus $h < 0.205$, and the centered images $z_1, \dots, z_5$ satisfy $t \coloneqq \frac{1}{2}|\sum_{j=1}^5 z_j^2| < 0.590$.
\end{lem}

\begin{proof}
Up to rotation, write the original points as $w e^{i\gamma}$ on $\{w^5 = 1\}$ with $[\gamma] = 0$, where brackets denote the five-point average and subscripts denote discrete Fourier indices modulo $5$. Set $a_0 = |\gamma_1|$ and $d_0 = |\gamma_2|$. Parseval's identity on successive differences gives $10(a_0^2\sin^2(\pi/5) + d_0^2\sin^2(2\pi/5)) = \pi^2 W_0$, whence $d_0^2 + 0.381 a_0^2 < 0.105^2$. This ellipse bound persists along the continuous path $w e^{i\tau\gamma}$ ($\tau \in [0,1]$); below $\gamma$ refers to the phase at a point of this path. By Lemma~\ref{lem:balanced-coords}, the centering point varies continuously along the path. Write the centered points as $z = w e^{ic}$ with $[c] = 0$ after rotation, and put $a = |c_1|$, $d = |c_2|$, and $L = a^2 + d^2$.

Suppose on a closed subset of $[0,1]$ that $h \le 0.205$ and $L \le 0.0139$, so $d < 0.118$. The centering condition $(e^{ic})_{-1} = 0$ and Taylor expansion first give $a \le [c^2]/2 = L$, and then $(1 - d)a \le d^2/2 + [|c|^3]/6$. Using Cauchy--Schwarz and Minkowski on the two Fourier pairs gives $\|c\|_4^2 \le \sqrt{6}(a + d)^2$ and $[|c|^3]/6 \le (a + d)^3/\sqrt{3}$, yielding $a < 0.0094$.

Write $x = w e^{i(\gamma - \lambda)} = (z + b)/(1 + \bar{b} z)$ with $|b| = h$, and set $m = (\gamma + c)/2$ and $B = \bar{b} e^{-i\lambda/2}$. Dividing $x - z = b - \bar{b} x z$ by $2i w e^{i(m - \lambda/2)}$ gives $\sin((\gamma - \lambda - c)/2) = -\Im(B w e^{im})$. Since the continuous lift satisfies $(\gamma - \lambda - c)/2 = -\arg(1 + \bar{b} z) \in [-\arcsin h, \arcsin h]$, expanding $\arcsin$ gives (up to an additive constant)
\[
  \gamma = c - 2\Im(B w e^{im}) + R, \qquad \|R\|_\infty \le \frac{h^3}{3(1 - h^2)^{3/2}} \le 0.36 h^3.
\]
Put $A = |m_1|$, $D = |m_2|$, and $K = [m^2]/2$. From the ellipse bound, $\sqrt{a_0^2 + d_0^2} < 0.105/\sqrt{0.381} < 0.1702$ and $a_0 + d_0 < 0.105\sqrt{1 + 1/0.381} < 0.200$, giving $K < ((0.1702 + 0.118)/2)^2 < 0.0208$, $A + D < (0.200 + 0.0094 + 0.118)/2 < 0.164$, and $D \le (d_0 + d)/2$. Every Fourier coefficient of $e^{im} - 1 - im$ has modulus at most $K$; reading off the first two discrete Fourier modes of $2\Im(B w e^{im})$ (involving modes $(0,-2)$ and $(1,2)$ of $e^{im}$) gives
\[
  a_0 \ge h(1 - 2K - D) - a - 0.36 h^3, \qquad
  d_0 \ge d - h(A + D + 2K) - 0.36 h^3.
\]
The first inequality yields $0.8842 h \le a_0 + h d_0/2 + 0.0094 < 0.105\sqrt{1/0.381 + 0.205^2/4} + 0.0094 < 0.1804$, so $h < 0.2045$. Using $D \le (0.105 + 0.118)/2$, the two bounds also give $a_0 \ge 0.831 h - 0.0094$ and $d_0 \ge d - 0.221 h$. Completing the square in the ellipse gives $d_0 + 2a_0^2 < 2(0.105)^2/0.381 + 0.381/8 < 0.1056$, whence (using $(u - v)_+^2 \ge u^2 - 2uv$ for $u, v \ge 0$)
\[
  d < 0.1056 + 0.253 h - 1.38 h^2 \le 0.1056 + \frac{0.253^2}{4(1.38)} < 0.1174.
\]
Thus $L < 0.0094^2 + 0.1174^2 < 0.0139$. By continuity along $\tau \in [0,1]$, both strict bounds $h < 0.205$ and $L < 0.0139$ hold at $\tau = 1$.

Finally, after rotating so that $\prod_j z_j = 1$, the elementary symmetric polynomials of $z_1, \dots, z_5$ satisfy $e_1 = e_4 = 0$, $e_5 = 1$, $e_3 = \bar{e}_2$, and $|e_2| = t$. Hence $z_j^5 - 1 = z_j^2(\bar{e}_2 - e_2 z_j)$ for each $j$. Averaging the squared modulus over $j$ and using $[z] = 0$ gives $2t^2 = [|z^5 - 1|^2] = [|e^{5ic} - 1|^2] \le 25[c^2] = 50 L$, so $t < 0.590$.
\end{proof}

\begin{lem}[Energy and displaced mean for the five-vertex reciprocal test]\label{lem:pure-five-energy}
In centered five-vertex coordinates, put $t = \frac{1}{2}|\sum_{j=1}^5 z_j^2|$. Whenever $t \le 0.605$, the test function $G$ defined below has harmonic extension $M(a) > 0$ satisfying $M(a)^2 > 0.02420\,\E(G)$ throughout $\{|a| \le 0.22\}$.
\end{lem}

\begin{proof}
Rotate so that the quadratic coefficient of $S(z) = \prod_{j=1}^5 (1 - \bar{z}_j z)$ equals $-t$. Vanishing of the linear coefficient and conjugate-reciprocal symmetry give
\[
  S(z) = 1 - t z^2 + \eta(z^5 - t z^3), \qquad |\eta| = 1.
\]
Factor $S = A(1 - B)$ with outer factor $A(z) = 1 - t z^2$ and degree-five Blaschke product $B(z) = -\eta z^3(z^2 - t)/(1 - t z^2)$, mapping each $z_j$ to $1$. Put $x = z^2$, $q = t^2$, $T = 0.605$, $Q = T^2$, $D = |1 - tx|^2$, $k = (1 - q)^{-4/25}$, $\tau = q/Q$, and
\[
  u = \frac{|B'|}{5} = \frac{3}{5} + \frac{2(1 - q)}{5 D}, \qquad V = k u D^{1/5}.
\]
Take $F = w_{-1/5} - D_{1/5}$ from \eqref{eq:reciprocal-test} at $\beta = 1/5$, for which $b_1 = 1/6$, $b_2 = 2/33$, $\E(F) = 1/10$, and $0.13517 < C = 2\sum_{m=1}^\infty \zeta(2m)/25^m < 0.136$ (using the first three terms with $1.6449 < \zeta(2) < 1.645$, $1.082 < W < 1.083$, and $\zeta(6) > 1$). Define $G = V F(B)$.

Set $P = (1 - q)/D$, $y = q/(1 - q)$, $e = P - 1$, $S_0 = P + P^{-1} - 2$, and $L_c = 1 + e/5 + c S_0$, with $c_- = 1/35$ and $c_+ = 1/25$. We first verify the two-sided barrier
\begin{equation}\label{eq:five-both-barriers}
  L_{c_-} \le u P^{-1/5} \le L_{c_+}.
\end{equation}
In terms of $v = P^{1/5} \in [3/4, 4/3]$, the ratio $(u P^{-1/5} - 1 - e/5)/S_0$ equals
\[
  \frac{3 + 2v + v^2 - v^4}{5(v^{-2} + v^{-1} + 1 + v + v^2)^2},
\]
which increases strictly to $1/25$ at $v = 1$ and decreases strictly beyond $v = 1$ by opposite monotonicities of numerator and denominator. At both endpoints $v = 3/4, 4/3$, the numerator exceeds $4.28$ while $v^{-2} + \cdots + v^2 < 5.43$, keeping the ratio above $1/35$.

\paragraph{Energy bound.}
Put $U_c = u L_c^2$. The pushforward density $B_*(U_c\,dm)$ has constant coefficient $I_c$ and only one nonconstant frequency pair, at $\pm 2$, of modulus $J_c$, given by
\[
  I_c = 1 + (0.4 + 4c)y + (0.096 + 0.96 c + 6 c^2)y^2, \qquad
  J_c = 0.4(0.2 + c)^2 t y^2.
\]
To verify this, expand $U_c$ in powers $P^{-2}, \dots, P^3$. Because the involution $\phi_t$ interchanges $P$ and $P^{-1}$ with Jacobian $P$, we have $\int P = 1$, $\int P^2 = \int P^{-1} = 1 + 2y$, and $\int P^3 = \int P^{-2} = 1 + 6y + 6y^2$, giving $I_c$. Odd $B$-moments vanish by parity. In the variable $x = z^2$, $B^{2j}$ has zeros of orders $3j$ and $2j$ at $0$ and $t$, so $P^l B^{2j}$ is holomorphic and vanishes at $0$ for $-2 \le l \le 3$, $j \ge 1$ except when $(l, j) = (3, 1)$, in which case $P^3 B^2\,dx/x$ has residue $\eta^2 t y^2$ at $t$ and the coefficient of $P^3$ in $U_c$ is $0.4(0.2 + c)^2$. The same formulas hold at $t = 0$.

By Parseval, $\int_\T F^2\,dm \le C^2 + \E(F)/3 + \frac{4}{3}b_1^2 + \frac{2}{3}b_2^2 < 0.092$, so \eqref{eq:reciprocal-negative-mass} gives $\|(F\Lambda F)_-\| \le (3\sqrt{0.092} - C)/20 < 0.039$. Both barriers in \eqref{eq:five-both-barriers} are positive, and the pushforward density of $U_{c_+} - U_{c_-}$ is bounded above by $I_{c_+} - I_{c_-} + 2(J_{c_+} - J_{c_-}) \le 0.046 y + 0.016 y^2 + 0.0043 t y^2 < 0.033\tau$ (using $y < 0.578\tau$). We also have $J_{c_+} < 0.0047\tau$ and $(1 - q)^{0.08} I_{c_+} \le (1 - 0.08 q)(1 + 0.56 y + 0.144 y^2) \le 1 + 0.333\tau$ (all nonlinear Taylor coefficients in $q$ are nonnegative, and at $Q$ the value is below $0.971(1 + 0.56(0.578) + 0.144(0.578)^2) < 1.333$). Pairing the smooth pushforward density of $u V^2\,dm$ with $F\Lambda F$, using $h_2 = 4D_\beta/33 + 1/10 - 5/36 < 0.067$ from \eqref{eq:small-reciprocal-product}, $\Lambda(F(B)) = 5u(\Lambda F)(B)$, and \eqref{eq:covariant-marker-product}, we obtain
\[
  \E(G) \le 5\bigl[h_0(1 + 0.333\tau) + 2(0.0047)h_2\tau + 0.033(0.039)\tau\bigr]
    + 0.575 q\int_\T F^2\,dm \le 0.5 + 0.196\tau.
\]

\paragraph{Displaced mean bound.}
Let $r = |a| \le 0.22$, $s = t r^2 < 0.0293$, $\mathcal{P}_a$ denote harmonic extension at $a$, and $L_0 = L_{c_-} = (26 + 8P + P^{-1})/35$. By \eqref{eq:five-both-barriers}, $V/(1 - q)^{0.04}$ lies between $L_0$ and $L_0 + \frac{2}{175}S_0$. On $\T$, the holomorphic identities
\[
  P B = -\eta(1 - q)\frac{z^5}{A^2}, \qquad P^{-1} B = -\eta\frac{z(z^2 - t)^2}{1 - q}
\]
imply that $S_0 B^j$ is holomorphic and vanishes at $0$ for every $j \ge 1$, so $B_*(S_0\,dm) = 2y\,dm$. Moreover, the holomorphic series $L_0\sum_{j \ge 1} b_j B^j$ can be evaluated directly at $a$: since $|B(a)| \le r^3(t + r^2)/(1 - s) < 0.0072$, twice its modulus at $a$ is bounded above by
\[
  \frac{2b_1}{35(1 - 0.0072)}\left[26(0.0072) + \frac{8 r^5}{(1 - s)^2} + \frac{r(t + r^2)^2}{1 - q}\right] < 0.0033.
\]
The remaining error from replacing $V/(1 - q)^{0.04}$ by $L_0$ is bounded above by $\frac{4}{175}y\|P_a\|_\infty\int F_-\,dm \le 0.023 C y$, using $\|P_a\|_\infty \le (1 + r)/(1 - r) < 1.565$, $\int F_-\,dm \le (\sqrt{0.092} - C)/2 < 0.085$, and $C > 0.13517$. With $w = t a^2$, the harmonic extensions of $P - 1$ and $P^{-1} - 1$ at $a$ are $2\Re(w/(1 - w))$ and $2y - 2\Re(w)/(1 - q)$, yielding
\[
  \mathcal{P}_a L_0 \ge 1 + \frac{2y}{35} - 0.4 s - \frac{16 s^2}{35(1 - s)} > 1 + \frac{2y}{35} - 0.0122.
\]
The tangent to $\log(1 + y)$ at $y = 1/4$ gives $\log(1 + y) \le 0.025 + 0.8 y$; since $1 - q = (1 + y)^{-1}$ and $y < 0.578$, we obtain $(1 - q)^{0.04}(1 + 2y/35) \ge (0.999 - 0.032 y)(1 + 2y/35) \ge 0.999 + 0.023 y$. The linear terms in $y$ cancel, leaving
\[
  \frac{M(a)}{C} > 0.999 - 0.0122 - \frac{0.0033}{0.13517} > 0.962.
\]
Finally, $(0.962 C)^2 > 0.02420(0.5 + 0.196) \ge 0.02420\,\E(G)$ throughout $\{|a| \le 0.22\}$.
\end{proof}

\begin{prop}[Five convex vertices and one reflex vertex]\label{prop:five-one-reflex}
Every simple polygon $P$ with $p = 5$ strictly convex vertices and $q = 1$ reflex vertex satisfies $d(P) > 0.02420 > d(R_6)$, and in particular $J(P) < J(R_6)$.
\end{prop}

\begin{proof}
Assume for contradiction that $d(P) \le 0.02420$. By Lemmas~\ref{lem:marked-gap-masses} and~\ref{lem:five-center-hand}, centering the five convex prevertices yields $t < 0.590 < 0.605$ and a displaced balanced origin $a \in \D$ with $|a| < 0.205 < 0.22$. For the test function $G$ of Lemma~\ref{lem:pure-five-energy}, \eqref{eq:regular-potential-transport} gives $\int_\T G\,d\mu \le 0$ (since $V$ is proportional to $u|A|^{2/5}$), while its mean in balanced coordinates is $M(a) > 0$. Applying \eqref{eq:duality} gives $d(P) \ge M(a)^2/\E(G) > 0.02420$, a contradiction. Lemma~\ref{lem:regular-six-short} gives $d(R_6) < 0.02420$, whence $J(P) < J(R_6)$.
\end{proof}

\section{Proof of the main theorem and the spectral gap}\label{sec:assembly}

We now combine the comparisons established in Sections~\ref{sec:regular}--\ref{sec:low-reflex} to complete the proofs of Theorems~\ref{thm:main} and~\ref{cor:defective-gap}.

\begin{proof}[Proof of Theorem~\ref{thm:main}]
Let $P$ be a simple planar polygon with $p$ strictly convex vertices and $q$ reflex vertices after merging consecutive collinear sides, and set $n = p + \min(q, 2)$.
\begin{itemize}[leftmargin=*]
\item If $q \ge 2$, then $n = p + 2$, and Proposition~\ref{prop:pure-two-finite} yields $J(P) < J(R_{p+2})$.
\item If $q = 1$, then $n = p + 1$, and strict inequality $J(P) < J(R_{p+1})$ holds by Proposition~\ref{prop:hull} ($p = 3$), Proposition~\ref{prop:four-one-reflex} ($p = 4$), Proposition~\ref{prop:five-one-reflex} ($p = 5$), and Proposition~\ref{prop:phase64} ($p \ge 6$).
\item If $q = 0$, then $P$ is a convex $p$-gon and $n = p$; Theorem~\ref{thm:pure-triangles} ($p = 3$), Proposition~\ref{prop:convex-quad} ($p = 4$), Theorem~\ref{thm:convex-five-nine} ($5 \le p \le 9$), and Theorem~\ref{thm:convex-ten} ($p \ge 10$) establish $J(P) \le J(R_p)$, with equality if and only if $P$ is similar to $R_p$.
\end{itemize}
Finally, for any simple polygon $P$ with at most $N$ sides ($N \ge 3$), we have $n \le p + q \le N$. By the strict monotonicity \eqref{eq:reg-monotonicity} of $J(R_k)$, it follows that $J(P) \le J(R_n) \le J(R_N)$, with equality if and only if $n = N$, $q = 0$, and $P$ is similar to $R_N$.
\end{proof}

\begin{proof}[Proof of Theorem~\ref{cor:defective-gap}]
Every polygon $P \in \mathcal{D}_N$ has at most $p \le N - 1$ strictly convex vertices. We first determine the asymptotic infimum of $J(R_N) - J(P)$ as $N \to \infty$. If $p \ge N/2$, Theorem~\ref{thm:unrestricted-reflex-asymptotic} and the monotonicity \eqref{eq:reg-monotonicity} give $d(P) \ge d(R_{N-1}) - O(N^{-5})$. If $p < N/2$, Lemma~\ref{lem:gap} yields the stronger bound $d(P) \ge (c p^3 + b p^2 - e)^{-1} > d(R_{N-1})$ for all sufficiently large $N$, in view of \eqref{eq:regular-two}. Since $R_{N-1} \in \mathcal{D}_N$, we conclude that
\[
  \inf_{P \in \mathcal{D}_N} d(P) = d(R_{N-1}) + O(N^{-5}) \qquad (N \to \infty).
\]
Using $J(P) = 2\pi/(1 + d(P))$ together with $d(R_{N-1}) - d(R_N) = 6Z N^{-4} + O(N^{-5})$ from \eqref{eq:regular-two} yields
\[
  \inf_{P \in \mathcal{D}_N}\bigl(J(R_N) - J(P)\bigr) = 12\pi Z N^{-4} + O(N^{-5}).
\]

To establish the uniform lower bound $J(R_N) - J(P) \ge c N^{-4}$ for all $N \ge 4$, it suffices to verify a strictly positive margin $\inf_{P \in \mathcal{D}_N}(J(R_N) - J(P)) > 0$ for each fixed $N \ge 4$.
\begin{itemize}[leftmargin=*]
\item If $P$ has fewer than $N$ genuine sides, then $p + q \le N - 1$, so Theorem~\ref{thm:main} and \eqref{eq:reg-monotonicity} give $J(P) \le J(R_{N-1}) < J(R_N)$.
\item If $P$ has $N$ genuine sides and $q \ge 2$ reflex vertices, then $p \le N - 2$, so Proposition~\ref{prop:pure-two-finite} gives $d(P) > d(R_N) + (N - 1)^{-4}$.
\item If $P$ has $N$ genuine sides and $q = 1$ reflex vertex, then $p = N - 1$. For $p = 3$, Proposition~\ref{prop:hull} gives $J(P) < 3\sqrt{3} < J(R_4)$. For $p = 4$ and $p = 5$, Propositions~\ref{prop:four-one-reflex} and~\ref{prop:five-one-reflex} exclude $d(P) \le 0.05215$ and $d(P) \le 0.02420$, both strictly above $d(R_5)$ and $d(R_6)$ by Lemma~\ref{lem:regular-six-short}. For $p \ge 6$, the proof of Proposition~\ref{prop:phase64} uses $d(R_{p+1})$ exclusively through the strict inequalities of Corollary~\ref{cor:regular-uniform-nine} (and, at $p = 6$, the strict inequality $d_7 < 0.013273$); because these strict inequalities remain valid in a right neighborhood of $d(R_{p+1})$, the same argument excludes a positive neighborhood of $d(R_N)$.
\end{itemize}
Combining the fixed-$N$ margins with the large-$N$ asymptotic limit yields the uniform bound $J(R_N) - J(P) \ge c N^{-4}$ for all $N \ge 4$.
\end{proof}

\section*{Statement on the use of AI tools}
AI tools (including ChatGPT, OpenAI, and Gemini, Google) were used during the preparation of this manuscript to assist with mathematical exposition, manuscript restructuring, and language polishing, and to explore and refine proof strategies through interactive discussion. All proofs in the present version are purely analytic and self-contained; the authors take full responsibility for the independent verification of all mathematical results, proofs, and explicit bounds, as well as for the final content of the article.

\end{document}